\documentclass[12pt]{amsart}
\usepackage[utf8]{inputenc}
\usepackage[margin=1in]{geometry}
\usepackage[dvipsnames]{xcolor}

\usepackage{bbm}
\usepackage{tikz}
\usetikzlibrary{knots}
\usetikzlibrary{arrows}
\usepackage{tikz-cd}
\usepackage{pgffor}
\usepackage{amssymb,amsmath,amsfonts,mathrsfs,amsthm}
\usepackage{enumitem}
\usepackage{adjustbox}
\usepackage{hyperref}
\usepackage{relsize}
\usepackage{graphicx}
\usepackage[font={small}]{caption}
\usepackage{float}
\usepackage{subcaption}
\usepackage{changepage}
\usepackage{changes}
\usepackage[toc]{appendix}
\usepackage{pgfplots}
\pgfplotsset{compat=1.7}
\usepackage{adjustbox}
\usepackage{varwidth}
\usepackage{chngcntr}
\usepackage{tcolorbox}
\usepackage{mathtools}
\usepackage{changes}
\usepackage{chngcntr}
\usepackage{standalone}
\usepackage{relsize}
\usetikzlibrary{calc}
\usetikzlibrary{patterns}
\usetikzlibrary{patterns.meta}

\newcommand{\RR}{\mathbb{R}}
\renewcommand{\SS}{\mathbb{S}}

\newcommand{\Z}{\mathbb{Z}}

\newcommand{\A}{\mathbb{A}}

\renewcommand{\l}{\ell}

\renewcommand{\d}{\partial}

\newcommand{\til}[1]{\widetilde{#1}}
\renewcommand{\b}[1]{\overline{#1}}

\newcommand{\ggmod}{{\rm -}\operatorname{ggmod}}
\newcommand{\gggmod}{{\rm -}\operatorname{g}_3\operatorname{mod}}

\newcommand{\kk}{\Bbbk}
\newcommand{\qdeg}{\operatorname{qdeg}}
\newcommand{\adeg}{\operatorname{adeg}}
\newcommand{\lr}[1]{\vert {#1} \vert}

\newcommand{\lra}{\longrightarrow}
\newcommand{\mcD}{\mathcal{D}}

\newcommand{\rar}[1]{\xrightarrow{#1}}

\newcommand{\al}{\alpha}

\renewcommand{\o}{\otimes}
\newcommand{\eps}{\varepsilon}

\renewcommand{\:}{\colon}

\newcommand{\sqrtLF}{\sqrt{f'|_{L\cap F}}}

\newcommand{\Ra}{R_\alpha} 
\newcommand{\Aa}{A_\alpha} 
\newcommand{\Fa}{\mathcal{F}_{\alpha}} 
\newcommand{\Ga}{\mathcal{G}_{\alpha}} 
\renewcommand{\P}{\mathcal{P}} 
\newcommand{\Fr}{\operatorname{Fr}} 

\newcommand{\ACob}{\operatorname{ACob}}
\newcommand{\AFoam}{\operatorname{AFoam}}
\newcommand{\AFoamor}{\operatorname{AFoam}_{\operatorname{or}}}
\newcommand{\grank}{\operatorname{grank}}
\newcommand{\rank}{\operatorname{rank}}
\newcommand{\MV}{\operatorname{MV}}
\newcommand{\anch}{\mathrm{an}}
\newcommand{\forget}{\operatorname{forget}}

\newcommand{\Comp}{\operatorname{Comp}} 

\newcommand{\brak}[1]{\ensuremath{\left\langle #1\right\rangle}}
\newcommand{\adm}{\operatorname{adm}}

\newtheorem{theorem}{Theorem}[section]
\newtheorem{lemma}[theorem]{Lemma}
\newtheorem{proposition}[theorem]{Proposition}
\newtheorem{corollary}[theorem]{Corollary}

\newcommand*\circled[1]{\tikz[baseline=(char.base)]{
            \node[shape=circle,draw,inner sep=1pt] (char) {${#1}$};}} 
            
\pgfdeclarepattern{name=hatch,parameters={\hatchsize,\hatchangle,\hatchlinewidth},bottom left={\pgfpoint{-.1pt}{-.1pt}},top right={\pgfpoint{\hatchsize+.1pt}{\hatchsize+.1pt}},tile size={\pgfpoint{\hatchsize}{\hatchsize}},tile transformation={\pgftransformrotate{\hatchangle}},code={\pgfsetlinewidth{\hatchlinewidth}\pgfpathmoveto{\pgfpoint{-.1pt}{-.1pt}}\pgfpathlineto{\pgfpoint{\hatchsize+.1pt}{\hatchsize+.1pt}}\pgfpathmoveto{\pgfpoint{-.1pt}{\hatchsize+.1pt}}\pgfpathlineto{\pgfpoint{\hatchsize+.1pt}{-.1pt}}\pgfusepath{stroke}}} 

\tikzset{hatch size/.store in=\hatchsize,hatch angle/.store in=\hatchangle,hatch line width/.store in=\hatchlinewidth,hatch size=5pt,hatch angle=0pt,hatch line width=.5pt,} 

\theoremstyle{definition}

\theoremstyle{definition}

\theoremstyle{remark}
\newtheorem{example}[theorem]{Example}

\theoremstyle{remark}
\newtheorem{remark}[theorem]{Remark}

\theoremstyle{definition}
\newtheorem{definition}{Definition}[section]

\theoremstyle{remark}
\newtheorem{rmk/}{Remark}
\theoremstyle{definition}
\newtheorem{case/}{Case}

\date{May 3, 2021} 

\title{Anchored foams and annular homology} 

\begin{document}

\author[R. Akhmechet]{Rostislav Akhmechet}
\address{Department of Mathematics, University of Virginia, Charlottesville VA 22904-4137}
\email{\href{mailto:ra5aq@virginia.edu}{ra5aq@virginia.edu}}

\author[M. Khovanov]{Mikhail Khovanov}
\address{Department of Mathematics, Columbia University, New York NY 10027}
\email{\href{mailto:khovanov@math.columbia.edu}{khovanov@math.columbia.edu}}

\begin{abstract} We describe equivariant SL(2) and SL(3) homology for links in the solid torus via foam evaluation. The solid torus is replaced by 3-space with a distinguished line in it. Generators of state spaces for annular webs are represented by foams with boundary that may intersect the distinguished line; intersection points, called anchor points, contribute additional terms, reminiscent of square roots of the Hessian, to the foam evaluation. Both oriented and unoriented SL(3) foams are treated in the paper. 
\end{abstract}

\maketitle
\tableofcontents

\section{Introduction} 
Asaeda-Przytycki-Sikora~\cite{APS} homology of links in the solid torus has led to a number of interesting developments  \cite{Roberts,GN, BG,  GLWsl2, GLWsinvt, BPW,Akh} and  extensions of their work to $SL(N)$ and $GL(N)$ link homology in the solid torus~\cite{QR, QW, QRS}. 

$GL(N)$ and $SL(N)$ link homology theories are closely related to foam evaluation. This connection was made the most transparent by the work of 
Robert and Wagner~\cite{RWfoamev}, who wrote down a combinatorial formula for $GL(N)$ closed foam evaluation that allows to build $GL(N)$ link homology from the ground up, bypassing categorical approaches to the latter. A variation of their formula was used to evaluate unoriented $SL(3)$ foams~\cite{KRfoamev}, giving a combinatorial approach to some of the structures discovered by Kronheimer and Mrowka~\cite{KM1}. 

In this paper we extend foam evaluation framework to build equivariant $SL(2)$ and $SL(3)$ state spaces for annular webs and, consequently, equivariant $SL(2)$ and $SL(3)$ homology for links in the solid torus. Our construction complements earlier work~\cite{QR, QW} on the subject. 
The same approach allows to define state spaces for unoriented $SL(3)$ annular webs, extending the construction in~\cite{KRfoamev}. 

\vspace{0.1in} 

In the APS (Asaeda-Przytycki-Sikora) annular homology and its equivariant and $SL(N)$ generalizations, one first defines state spaces for annular $SL(2)$ and $SL(N)$ webs, where annular $SL(2)$ webs are just collections of embedded circles in an annulus. 

Our idea is to think of an open solid torus as the complement to a line $L$ in $\RR^3$, chosen for convenience to be the $z$-axis. An annular $SL(N)$ web $\Gamma$ is then placed into the $xy$-plane with $(0,0)$ removed. To define its state space $\brak{\Gamma}$, we consider $SL(N)$ foams $F$ in the half-space $\RR^3_-$ bounded by the $xy$-plane such that $\Gamma$ is the boundary of $F$. These foams may intersect the $z$-axis, and we refer to the intersection points as \emph{anchor points} and to such foams as \emph{anchored foams}. Anchor points additionally carry a label from $1$ to $N$, and we modify foam evaluation by adding a new type of factors associated to anchor points. 

In this paper we treat $N=2$ and $N=3$ cases, with modified evaluations given by formulas (\ref{eq:evaluation for coloring}) and (\ref{eq:oriented anchored sl3 eval color}), respectively, also see (\ref{eq:unoriented eval color}) for the unoriented $SL(3)$ anchored foam evaluation. 

Anchored foam evaluation take values in the ring of polynomials rather than the ring of symmetric polynomials. One starts with an admissible coloring $c$ of facets of a foam $F$, as usual. 
An anchor point labeled $i$ lying on a facet of color $j$ contributes $\delta_{i,j}\sqrt{\pm f'(x_i)}$ to the evaluation $\brak{F,c}$, 
where, in the $SL(3)$ case as an example, $f(x)=(x-x_1)(x-x_2)(x-x_3)$ is the polynomial of degree three with roots $x_1,x_2,x_3$. 
The full evaluation $\brak{F}$ is given by summing over $\brak{F,c}$ for all admissible colorings $c$.  
We check integrality property of these evaluations, with $\brak{F}$ a polynomial in $x_1,x_2,x_3$, in the $SL(3)$ case.

\vspace{0.1in} 

Given evaluations of anchored closed foams, one can form state spaces for annular webs. We show that this modified evaluation, with anchor points contributing $\delta_{i,j}\sqrt{\pm f'(x_i)}$,  perfectly matches the structure of state spaces of annular homology, in $SL(2)$ and $SL(3)$ cases. The construction also allows us to define unoriented $SL(3)$  homology for annular trivalent graphs, extending~\cite{KRfoamev} to the annular framework.

With state spaces at hand, it is straightforward to define annular $SL(2)$ and $SL(3)$ link homology, by analogy with \cite{Kh,APS, BNtangles,Akh} in the $SL(2)$ setting, with \cite{KRfoamev} in the unoriented $SL(3)$ setting, and with  \cite{Khsl3, MV, RWfoamev} in the oriented $SL(3)$ setting.  State spaces and link homology carry additional gradings coming from intersection points of foams with the $z$-axis. We show  that the result matches equivariant $SL(2)$ homology~\cite{Akh} of the first author. A simple modification of the construction (truncating the ground ring by sending $x_i$'s to $0$ upon evaluation) gives a foam approach to the original APS homology. We expect that the non-equivariant variant of our $SL(3)$ construction recovers $N=3$ case of the homology in~\cite{QR}. It seems that the equivariant annular $SL(3)$ homology, as described in the present paper, is new.

\vspace{0.1in} 

Section~\ref{sec:anchored surfaces} describes $SL(2)$ homology via anchored foams. 
The evaluation is defined in Section~\ref{subsec_anch_sl2}, 
which also contains the skein relations for anchored $SL(2)$ foams.  The state spaces are studied in Section~\ref{subsec_state_sp_2}.
The state space of $n$ circles in the annulus is a free module of rank $2^n$ over the ground ring $R_{\alpha}$ of polynomials in two variables, see Theorem~\ref{thm:state space sl2}. The numbers of contractible and essential circles control the bigraded rank.  This section also discusses categories of \emph{anchored} and \emph{annular} cobordisms. Annular cobordisms between annular $SL(2)$ webs are disjoint from the $z$-axis, while anchored cobordism may intersect it. 

Theorem~\ref{thm_functor_iso} identifies the annular cobordism functor with that constructed in~\cite{Akh}. Consequently, equivariant annular $SL(2)$ link homology~\cite{Akh}  can be rederived via anchored foams. To obtain the original APS homology, one can use anchored foam evaluation, combined with the homomorphism $R_{\alpha}\lra \Z$  taking $\alpha_1,\alpha_2$  to $0$ to get state spaces and cobordism maps in the APS theory.

\vspace{0.1in}

Section~\ref{sec:unoriented sl3} constructs the state spaces for the annular unoriented $SL(3)$ foam theory, extending the construction of~\cite{KRfoamev}. We start with the evaluation (Section~\ref{subsec_unor_anch}), followed by skein relations on annular foams (Section~\ref{subsec_skein}) and properties of state spaces (Section~\ref{sec:state spaces unoriented}).  Section~\ref{sec:Remark on Lee's theory} describes similarities between anchor points contributions and Lee's theory, given by inverting the discriminant in the ground ring. Similar to the planar case~\cite{KRfoamev}, we don't know a way  to describe  the state space of an annular web when regions of valency at most four, allowing an inductive simplification, are absent.

\vspace{0.1in} 

In Section~\ref{sec_oriented_hom} we describe annular equivariant $SL(3)$ link homology, based on anchored (annular) oriented $SL(3)$ foams. This homology extends Mackaay-Vaz~\cite{MV} equivariant $SL(3)$ homology of links in $\RR^3$, also see~\cite{Khsl3,MN,Clark,Robert} for the non-equivariant homology in $\RR^3$. 
We start with a review of oriented $SL(3)$ foams in Section~\ref{sec:oriented sl3 foams} and then follow a similar route to that of the earlier sections. 

\vspace{0.1in} 

We expect that our construction admits a  generalization to $SL(N)$ homology for all $N$ via an extension of the Robert-Wagner formula~\cite{RWfoamev} to the anchored case.

\vspace{0.1in} 

{\bf Acknowledgments:} M.K. was partially supported by NSF grant DMS-1807425  while working on the paper. R.A. was supported by the Jefferson Scholars Foundation. R.A. would like to thank his advisor Slava Krushkal for encouraging him to pursue this project. 

\vspace{0.1in} 

%
%

\section{\texorpdfstring{$SL(2)$}{SL(2)} anchored homology}
\label{sec:anchored surfaces}

\subsection{Anchored surfaces and their evaluations}
\label{subsec_anch_sl2}

Consider the integral polynomial ring $\Ra = \Z[\al_1, \al_2]$ in two variables $\al_1, \al_2$. Define a grading on $\Ra$ by setting \begin{equation}
\label{eq:grading of ground ring}
\deg(\al_1) = \deg(\al_2) = 2.
\end{equation}
Denote by $\tau$ the nontrivial involution of $\{1,2\}$. It is given by $\tau(i)=3-i$, for $i\in \{1,2\}$. Also denote by $\tau$ the induced involution of $\Ra$ which permutes $\al_1,\al_2$, so that $\tau(\alpha_i)=\alpha_{\tau(i)}=\alpha_{3-i}$. Let $R$ be the $\tau$-invariant subring of $\Ra$, which consists of symmetric polynomials in $\al_1, \al_2$. The subring $R$ is itself a polynomial ring, $R= \Z[E_1, E_2]$, where $E_1, E_2$ are elementary symmetric polynomials in $\al_1, \al_2$, 
\[
E_1 = \al_1 + \al_2, \hskip1em E_2 = \al_1 \al_2. 
\]
Degrees of $E_1$ and $E_2$ are $2$ and $4$, respectively.

Let $L\subset \RR^3$ denote the $z$-axis, $L= (0,0) \times \RR$. Let $S\subset \RR^3$ be a closed, smoothly embedded surface which intersects $L$ transversely. The surface $S$ may be decorated by dots, disjoint from $L$, that can otherwise float freely on components of $S$. The intersection points $S\cap L$ are called \emph{anchor points}. Fix a labeling $\l$, which is a map  from the set of anchor points to $\{1, 2\}$, 
\[
\l \: S \cap L \to \{1, 2\}.
\]
Order the anchor points by $1, \ldots, 2k$, read from bottom to top, so that the labeling $\l$ consists of a choice $\l(j) \in \{1, 2\}$ for each $1\leq j \leq 2k$. We will define an evaluation 
\[
\brak{S} \in \Ra. 
\]
for $S$ with the fixed labeling $\l$, which is omitted from the notation. 

Let $\Comp(S)$ denote the set of connected components of $S$. A \emph{coloring} of $S$  is a function $c \: \Comp(S) \to \{1,2\}$, and we denote by $\mathrm{adm}(S)$ the set of colorings of $S$. Surface $S$ has $2^{|\Comp(S)|}$ colorings. Fix a coloring $c$. For $i=1,2$, let $d_i(c)$ denote the number of dots on components colored $i$. Let $S_2$ denote the union of the $2$-colored components. For $1\leq j \leq 2k$, let $c(j)$ denote the color of the $j$-th anchor point, induced by $c$, which may in general be different from the fixed label $\l(j)$. Define 
\begin{equation}
\label{eq:evaluation for coloring}
\brak{S,c} = (-1)^{\chi(S_2)/2} \frac{\al_1^{d_1(c)} \al_2^{d_2(c)} \left( \prod_{j=1}^{2k} (\al_{c(j)} - \al_{\l(j)})
\right)^{1/2}}{(\al_1-\al_2)^{\chi(S)/2}}.
\end{equation}

Note that $\chi(S_2)$ is even since $S_2$ is a closed surface in $\RR^3$. Let us explain the square root in the above equation.

Each component $S'$ of $S$ intersects $L$ at an even number of points $p_1,\dots, p_{2m}$, which can be ordered as encountered along $L$, from bottom to top.  Suppose $S'$  is colored by $c(S') = j$, and moreover $S'$ contains an anchor point labeled $j$. Then the product $\prod_{j=1}^{2m} (\al_{c(j)} - \al_{\l(j)})=0$, since it contains a term $\alpha_j-\alpha_j=0$, and the entire evaluation $\brak{S,c}=0$. Thus, the evaluation  \eqref{eq:evaluation for coloring} is only nonzero when the anchor points on a component $S'$ colored $j$ are all labeled by the complementary color $\tau(j)$. In this case, each component contributes an even number of factors of either $\al_1 - \al_2$ or $\al_2-\al_1$ to the product $\prod_{j=1}^{2m} (\al_{c(j)} - \al_{\l(j)})$, and we define the square root to be $(\al_1 - \al_2)^m$ or $(\al_2-\al_1)^m$, respectively. If $S'$ has no anchor points, this term is $1$ and can be removed from the product. 

Note that the evaluation is the product of evaluations of individual components, 
\begin{equation}
    \brak{S,c} \ = \ \prod_{S'\in \Comp{(S)}} \brak{S',c(S')}. 
\end{equation}
Thus, if $S'$ is colored $1$ by $c'=c(S')$, has $2k$ anchor points all labeled $2$ and carries $d$ dots, then
\begin{equation}
    \brak{S',c'}  = \alpha_1^d (\alpha_1-\alpha_2)^{k-\chi(S')/2}. 
\end{equation}
If 
$S'$ is colored $2$ by $c'=c(S')$, has $2k$ anchor points all labeled $1$ and carries $d$ dots, then
\begin{equation}
    \brak{S',c'}  = (-1)^{\chi(S')/2+k}\alpha_2^d (\alpha_1-\alpha_2)^{k-\chi(S')/2}=
    \alpha_2^d (\alpha_2-\alpha_1)^{k-\chi(S')/2}. 
\end{equation}
Otherwise, if one of the anchor points has the same label as the color of $S'$, the evaluation $\brak{S',c'}=0$ and $\brak{S,c}=0$. 

Define the evaluation of $S$ by 
\begin{equation}
\label{eq:evaluation} 
\brak{S} = \sum_c \brak{S,c},
\end{equation}
where the sum is over all colorings of $S$. Note that if $S \cap L = \varnothing$, then $\brak{S}$ agrees with the evaluation in \cite{RWfoamev, KRfrobext}. Also note that $\brak{S} = 0$ if a component of $S$ has two anchor points with different labels $1$, $2$. 

We have 
\begin{equation}
    \brak{S} \ = \ \prod_{S'\in \Comp{S}} \brak{S'}, 
\end{equation}
that is, evaluation of $S$ is the product of evaluations over connected components of $S$. 

We can rewrite $\brak{S}$ as follows. First, suppose $S$ is connected, carrying $d$ dots, with $2k\geq 0$ anchor points. For $i=1,2$, let $c_i$ denote the coloring of $S$ by $i$. Define 
\begin{align}
\brak{S, c_1} &= \frac{\al_1^{d} ((\al_1-\al_{\l(1)}) \cdots (\al_1 - \al_{\l(2k)}) )^{1/2}}{(\al_1 - \al_2)^{\chi(S)/2}} ,\label{eq:color 1} \\
 \brak{S, c_2} &= (-1)^{\chi(S)/2} \frac{\al_2^{d} ((\al_2-\al_{\l(1)}) \cdots (\al_2 - \al_{\l(2k)}) )^{1/2}}{(\al_1 - \al_2)^{\chi(S)/2}} . \label{eq:color 2} 
\end{align}

Again, square roots in the above equations are taken in the natural way. If $S$ has oppositely labeled anchor points then both \eqref{eq:color 1} and \eqref{eq:color 2} are zero. If all anchor points are labeled $1$, then \eqref{eq:color 1} is zero, whereas \eqref{eq:color 2} is equal to 
\[
\brak{S, c_2} = (-1)^{\chi(S)/2} \frac{\al_2^d (\al_2 - \al_1)^k}{(\al_1 - \al_2)^{\chi(S)/2}}.
\]
On the other hand, if all anchor points are labeled by $2$ then \eqref{eq:color 2} is zero and \eqref{eq:color 1} equals 
\[
\frac{\al_1^d (\al_1 - \al_2)^k}{(\al_1 - \al_2)^{\chi(S)/2}}.
\]
Then for connected $S$ with anchor points, we have 
\[
\brak{S} = \brak{S,c_1} + \brak{S,c_2},
\]
where at most one of the summands $\brak{S, c_i}$ is nonzero.

Clearly the evaluation is multiplicative under disjoint union. That is, if $S = S_1 \sqcup \cdots \sqcup S_n$, then 
\[
\brak{S} = \brak{S_1} \cdots \brak{S_n}.
\]

\begin{remark}
Unlike closed foam evaluations appearing elsewhere \cite{RWfoamev, KRfoamev, KK, KRfrobext, RoseWedrich}, our evaluation  does not in general produce a symmetric function. The following examples illustrate this. 
\end{remark}

\begin{example}
Let $S$ be a sphere intersecting $L$ in two points with labels $i$ and $j$ and carrying $d$ dots. If $i\neq j$, then each coloring $c$ yields $\brak{S, c}=0$. If both anchor points are labeled $1$, then only coloring $S$ by $2$ contributes to the sum, and we have 
\[
\brak{S} = \brak{S, c_2} = - \frac{\al_2^d (\al_2-\al_1)}{\al_1 - \al_2} = \al_2^d.
\]
On the other hand, if both anchor points are labeled $2$, then 
\[
\brak{S} = \brak{S, c_1} = \al_1^d.
\]
This is summarized pictorially in \eqref{eq:sphere}. Both signs are positive since $k+\chi(\SS^2)/2=1+1=2$ is even. 

\begin{equation}
\begin{aligned}\label{eq:sphere}
\includestandalone{sphere}
\end{aligned}
\end{equation}
\end{example} 
Note that these evaluations are not symmetric in $\alpha_1,\alpha_2$.

\begin{example}
\label{ex:evaluation for surface disjoint from L}
More generally, let $S$ be a genus $g$ surface with $d$ dots and $2k$ anchor points. If $k=0$ (that is, if $S$ is disjoint from $L$) then the evaluation is 
\[
\brak{S} = \frac{\al_1^d + (-1)^{g-1} \al_2^d}{(\al_1 - \al_2)^{1-g}}.
\]
On the other hand, if $k>0$, then 
\begin{equation}
\label{eq:general formula}
\brak{S} = 
\begin{cases}
 \al_2^d (\al_2 - \al_1)^{k+g-1} & \text{ if } \l(1) = \cdots = \l(2k) = 1, \\
\al_1^d (\al_1 - \al_2)^{k + g-1} & \text{ if } \l(1) = \cdots = \l(2k) = 2 , \\
0 & \text{ otherwise.}
\end{cases}
\end{equation}
\end{example}

\begin{proposition}
\label{prop:evaluation is polynomial}
For any anchored surface $S\subset \RR^3$ with $d$ dots and $2k$ anchor points, its evaluation $\brak{S}$ is a homogeneous polynomial in $\al_1$ and $\al_2$ of degree $-\chi(S) + 2d +2k$. 
\end{proposition}

\begin{proof}
If $S$ does not intersect $L$, then this follows from Example \ref{ex:evaluation for surface disjoint from L}. Suppose that $S$ intersects $L$. It suffices to verify the statement for connected surfaces. If $S$ intersects $L$, then the statement follows from \eqref{eq:general formula}, since $k>0$. 
\end{proof}

We recall the following notation from \cite{KRfrobext}. For $i=1,2$, we allow surfaces to carry decorations $\circled{i}$ consisting of $i$ inscribed in a small circle. They must be disjoint from $L$ and are allowed to float freely along the connected component on which they appear. We call these \emph{shifted} dots. Diagrammatically, a shifted dot $\circled{i}$ is the difference between a dot and $\al_i$, 
\begin{equation}
\label{eq:shifted dot}
\begin{aligned}
\includestandalone{shifted_dot}.
\end{aligned}
\end{equation}

\begin{lemma}
\label{lem:shifted dots}
Let $S$ be an anchored foam and let $S\cup \circled{i}$ denote the anchored foam obtained by placing a shifted dot $\circled{i}$ on some component $S'$ of $S$. Then 
\[
\brak{S\cup \circled{i}} =
\begin{cases}
0 & \text{ if } S' \text{ has an anchor point labeled } \tau(i) \\
(-1)^i (\al_1 - \al_2) \brak{S} & \text{ if all anchor points on } S' \text{ are labeled } i. 
\end{cases}
\]
\end{lemma}

\begin{proof}
This is clear from the definitions.
\end{proof}

Lemma \ref{lem:shifted dots} is summarized diagrammatically in the relations \eqref{eq:shifted dots relations}.

\begin{equation}
\label{eq:shifted dots relations}
\begin{aligned}
\includestandalone{shifted_dots_relations1} \\
\includestandalone{shifted_dots_relations2} \\
\includestandalone{shifted_dots_relations3} 
\end{aligned}
\end{equation}
Alternatively, the skein relations \eqref{eq:shifted dots relations} may written compactly as in \eqref{eq:dot relations}. 

\begin{equation}
\label{eq:dot relations}
\begin{aligned}
\includestandalone{dot_relation}
\end{aligned}
\end{equation}

\begin{lemma}
The local relations \eqref{eq:two dots}, \eqref{eq:neck cutting}, and \eqref{eq:neck cutting line} hold. 

\begin{equation}
\label{eq:two dots}
\begin{aligned}
\includestandalone{two_dots}
\end{aligned}
\end{equation}

\begin{equation}
\label{eq:neck cutting}
\begin{aligned}
\includestandalone{neck_cutting}
\end{aligned}
\end{equation}

\begin{equation}
\label{eq:neck cutting line}
\begin{aligned}
\includestandalone{neck_cutting_line}
\end{aligned}
\end{equation}

\end{lemma}

\begin{proof}
The relation \eqref{eq:two dots} is straightforward. Let us now verify equation \eqref{eq:neck cutting}, which is proved in the  same way as for  non-anchored foams, see~\cite[Lemma 3.5]{KRfrobext}. Let $S$ denote the surface on the left, and let $F$ denote the surface obtained by surgering $S$ as shown on the right. Denote by $F^t$ (resp. $F^b$) the surface obtained from $F$ by placing an additional dot on the top (resp. bottom) depicted disk. Note that anchor points, as well as their labels, are the same for $F^t, F^b$, and $F$. Colorings of $F, F^t$, and $F^b$ are in a canonical bijection. There are four local models for a coloring of $F$, illustrated in Figure \ref{fig:neck cutting possibilities}. 

\begin{figure}
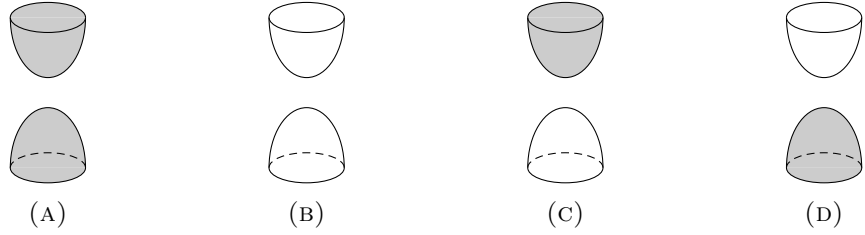

\centering
\subcaptionbox{
\label{fig:neck cutting 1}}[.2\linewidth]
{\includestandalone{neck_cutting_color1}
}
\subcaptionbox{
\label{fig:neck cutting 2}}[.2\linewidth]
{\includestandalone{neck_cutting_color2}
}
\subcaptionbox{
\label{fig:neck cutting 3}}[.2\linewidth]
{\includestandalone{neck_cutting_color3}
}
\subcaptionbox{
\label{fig:neck cutting 4}}[.2\linewidth]
{\includestandalone{neck_cutting_color4}
}

\caption{Local models for colorings of $F$. Shaded indicates color $1$ and solid white indicates color $2$.}\label{fig:neck cutting possibilities}
\end{figure}

Let $c$ be a coloring of $F$ of the type shown in Figure \ref{fig:neck cutting 3}, with the corresponding coloring of $F^t$ and $F^b$ still denoted by $c$. We have 
\[
\brak{F^t,c} = \al_1 \brak{F,c} \hskip2em \brak{F^b,c} = \al_2 \brak{F,c}, 
\]
hence $\brak{F^t,c} + \brak{F^b,c} - E_1 \brak{F,c} = 0$. A similar calculation holds for a coloring $c$ of Figure~\ref{fig:neck cutting 4} type. 

There is a natural bijection between colorings of $S$ and colorings of $F$ of Figures \ref{fig:neck cutting 1} and \ref{fig:neck cutting 2} types. Let $c$ be a coloring of $F$ of Figure \ref{fig:neck cutting 1} type, and continue to denote by $c$ the corresponding coloring of $S$. Then 
\begin{align*}
\chi(F) &= \chi(S) + 2 & \chi(F_2(c)) &= \chi(S_2(c)),\\
\brak{F^t,c} &= \al_1 \brak{F,c} & \brak{F^b,c} &= \al_1 \brak{F,c},
\end{align*}
so we have 
\[
\brak{F^t, c} + \brak{F^b,c} - E_1 \brak{F,c} = (\al_1 - \al_2) \brak{F,c} = \brak{S,c}. 
\] 
Finally, if $c$ is a coloring of $F$ of the Figure \ref{fig:neck cutting 2} type, then 
\begin{align*}
\chi(F) &= \chi(S) + 2  & \brak{F^t,c} &= \al_2 \brak{F,c},\\
 \chi(F_2(c)) &= \chi(S_2(c)) + 2  & \brak{F^b,c} &= \al_2 \brak{F,c},
\end{align*}
which yields
\[
\brak{F^t, c} + \brak{F^b,c} - E_1 \brak{F,c} = (\al_2 - \al_1) \brak{F,c} = - (\al_2-\al_1) \frac{\brak{S,c}}{\al_1 - \al_2} = \brak{S,c}.
\]

We now address equation \eqref{eq:neck cutting line}, where anchor points are present. Let $S$ denote the surface on the left-hand side of the equality.  Let $F^1, F^2$ denote the two anchored foams obtained by surgery on $S$ in which the new anchor points are both labeled $1$ or  $2$, respectively, so that \eqref{eq:neck cutting line} reads $\brak{S} = \brak{F^1} + \brak{F^2}$. For each $i=1,2$ there are four local models for a coloring of $F^i$, shown in Figure \ref{fig:neck cutting line possibilities}. Colorings $c$ in Figure~\ref{fig:neck cutting line 3} and Figure~\ref{fig:neck cutting line 4} evaluate to zero for both $i=1,2$, 
\[ \brak{F^1,c}=\brak{F^2,c}=0 
\] 
and they don't correspond to any colorings of $S$. There is a natural bijection between colorings of $S$ and colorings of $F^i$ of the types in Figures \ref{fig:neck cutting line 1} and \ref{fig:neck cutting line 2}.

\begin{figure}
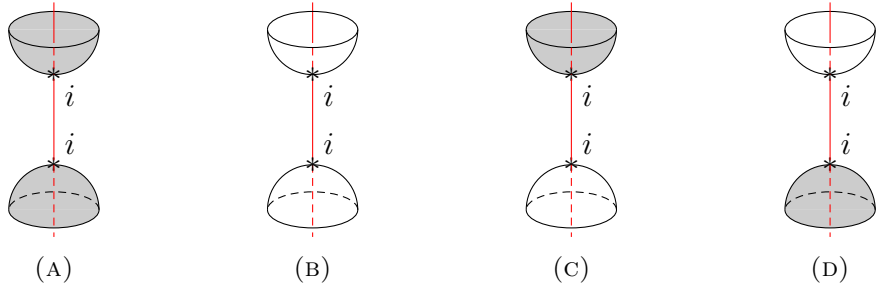

\centering
\subcaptionbox{
\label{fig:neck cutting line 1}}[.2\linewidth]
{\includestandalone{neck_cutting_line_color1}
}
\subcaptionbox{
\label{fig:neck cutting line 2}}[.2\linewidth]
{\includestandalone{neck_cutting_line_color2}
}
\subcaptionbox{
\label{fig:neck cutting line 3}}[.2\linewidth]
{\includestandalone{neck_cutting_line_color3}
}
\subcaptionbox{
\label{fig:neck cutting line 4}}[.2\linewidth]
{\includestandalone{neck_cutting_line_color4}
}
\caption{Local models for colorings of $F^i$. Shaded indicates color $1$ and solid white indicates color $2$.}\label{fig:neck cutting line possibilities}
\end{figure}

Let $c$ be a coloring of $S$ in which the depicted region of $S$ in \eqref{eq:neck cutting line} is colored $1$, with the corresponding colorings of $F^1$ and $F^2$ still denoted $c$. We have immediately that $\brak{F^1,c} = 0$. On the other hand, 
\[
\chi(F^2) = \chi(S) + 2,  \hskip2em \chi(F^2_2(c)) =\chi(S_2(c)), 
\]
and $F^2$ has two additional anchor points compared to $S$, both labeled $2$ and their regions colored $1$. Therefore
\[ \brak{F^1,c}+\brak{F^2,c} =
\brak{F^2,c} = (\al_1 - \al_2) \frac{\brak{S,c}}{\al_1 - \al_2 } = \brak{S,c}.
\]

Now let $c$ be a coloring of $S$ in which the depicted region of \eqref{eq:neck cutting line} is colored $2$, and continue to denote by $c$ the corresponding 
colorings of $F^1$ and $F^2$. Then $\brak{F^2, c} = 0$. Since
\[
\chi(F^1) = \chi(S) + 2,  \hskip2em \chi(F^1_2(c)) =\chi(S_2(c))+2, 
\]
and $F^1$ contains two more anchor points labeled $1$ and colored $2$ than $S$ does, we obtain
\[  \brak{F^1,c} + \brak{F^2,c} =
\brak{F^1,c} = -(\al_2 - \al_1) \frac{\brak{S,c}}{\al_1 - \al_2 } = \brak{S,c}.
\]
Relation  $\brak{S} = \brak{F^1} + \brak{F^2}$ in Figure~\eqref{eq:neck cutting line} follows.
\end{proof}

Equation \eqref{eq:neck cutting} can also be written using shifted dots,

\begin{equation}
\label{eq:neck cutting shifted dots}
\begin{aligned}
\includestandalone{neck_cutting_shifted_dots}.
\end{aligned}
\end{equation}

\begin{corollary}
The local relation \eqref{eq:moving cup off line} holds.

\begin{equation}
\label{eq:moving cup off line}
\begin{aligned}
\includestandalone{moving_cup_off_line}
\end{aligned}
\end{equation}
\end{corollary}

\begin{proof}
This can be seen by applying the neck-cutting relation \eqref{eq:neck cutting} near the depicted contractible circle and evaluating the resulting anchored sphere according to the formula \eqref{eq:sphere}. 
\end{proof}

\subsection{State spaces}\label{subsec_state_sp_2}

Following \cite{BHMV,KRfrobext}, we can apply the universal construction to the evaluation described above. Let $\P = \RR^2 \setminus {(0,0)}$ denote the punctured plane. Given a collection $C$ of disjoint simple closed curves in $\P$, let $\Fr(C)$ denote the free $\Ra$-module with a basis consisting of properly embedded compact surfaces $S\subset \RR^2 \times (-\infty, 0]$ with $\d S = C$ and which are transverse to the ray $L_-:= (0,0) \times (-\infty, 0]$. The intersection $S\cap  L_-$ is a $0$-submanifold of $L_-$ and consists of finitely many points. 
 Moreover, each such surface $S$ must carry a labeling, a map 
\[
\l=\l_S \ : \ S\cap L_- \lra \{1,2\}
\] 
from the set of its intersection points with the ray $L_-$ (its \emph{anchor points}) to $\{1,2\}$. For a basis element $S\in \Fr(C)$, let 
$\b{S} \subset \RR^2 \times [0,\infty)$ denote its reflection through the plane $\RR^2$. Labels of anchor points do not change upon reflection.  For two basis elements $S, S' \in \Fr(C)$ denote by $\b{S} S'$ the closed anchored surface obtained by gluing $\b{S}$ to $S'$ along their common boundary $C$.

Define a bilinear form 
\begin{equation}
\label{eq:bilinear form}
(-,-) \: \Fr(C) \times \Fr(C) \to \Ra
\end{equation}
by setting $(S, S') = \brak{\b{S}S'}$. A direct computation shows that the form is symmetric, since for a closed surface $T$ the evaluation satisfies $\brak{\b{T}}=\brak{T}$. 

Define the \emph{state space of $C$}, denoted $\brak{C}$, to be the quotient of $\Fr(C)$ by the kernel 
\[
\{ x\in \Fr(C) \mid (x,y) = 0 \text{ for all } y\in \Fr(C)\}
\]
of this bilinear form. For a basis element $S\in \Fr(C)$, we will write $[S]$ to denote its equivalence class in $\brak{C}$. 

Equip the ground ring $\Ra$ with a bigrading by placing $\al_1, \al_2$ in bidegree $(2,0)$. We extend this  bigrading $(\qdeg, \adeg)$ to $\Fr(C)$ as follows. For a basis element $S\in \Fr(C)$ with $d$ dots and $m$ anchor points, set the quantum grading $\qdeg(S)\in \Z$ to be 
\begin{equation}
\label{eq:qdeg}
\qdeg(S) = -\chi(S) + 2d +m.
\end{equation}
Note that if $S$ is a closed surface, then $\brak{S} \in \Ra$ is a homogeneous polynomial of degree $\qdeg(S)$, following the degree convention \eqref{eq:grading of ground ring}. 

Next, let $\l(1), \ldots, \l(m)$ denote the labels of the anchor points of $S$, ordered from bottom to top, and define the annular grading $\adeg(S)\in \Z$ by setting 
\begin{equation}
\label{eq:adeg}
\adeg(S) = \sum_{i=1}^m (-1)^{i + \l(i) }.
\end{equation}
In other words, if the $i$-th anchor point $p_i$ is labeled $1$, then it contributes $1$ to $\adeg$ if $i$ is odd and $-1$ if $i$ is even. Likewise, if $p_i$ has label $2$ then it contributes $-1$ if $i$ is odd and $1$ if $i$ is even, see also Figure \ref{fig:adeg}. Multiplication by $\al_1, \al_2$ increases $(\qdeg, \adeg)$-bidegree by $(2,0)$. 

\begin{figure}
\centering
\begin{center}
\begin{tabular}{  l | c | c } 
& label $1$ & label $2$ \\ 
\hline
$i$ odd & $1$ & $-1$ \\ 
\hline
$i$ even & $-1$  & $1$ \\ 
\end{tabular}
\end{center}
\caption{The contribution of the $i$-th anchor point on $S$ to $\adeg(S)$.} \label{fig:adeg} 
\end{figure}

\begin{example}
Let $C$ consist of two non-contractible circles. The bidegree $(\qdeg, \adeg)$ of the four anchored surfaces in $\Fr(C)$ whose underlying surface consists of two disks each intersecting $L_-$ once are recorded in Figure \ref{fig:deg examples}. 

\begin{figure}
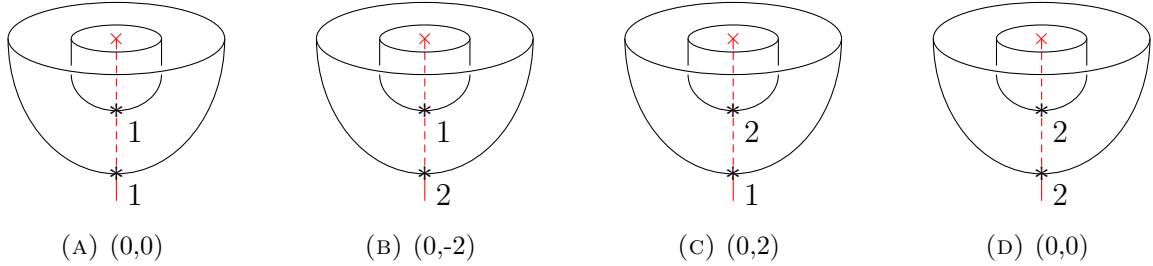

\centering
\subcaptionbox{ (0,0)
\label{fig:adeg ex1}}[.24\textwidth]
{\includestandalone{deg_ex1}
}
\subcaptionbox{ (0,-2)
\label{fig:adeg ex2}}[.24\textwidth]
{\includestandalone{deg_ex2}
}
\subcaptionbox{ (0,2)
\label{fig:adeg ex3}}[.24\textwidth]
{\includestandalone{deg_ex3}
}
\subcaptionbox{ (0,0) 
\label{fig:adeg ex4}}[.24\textwidth]
{\includestandalone{deg_ex4}
}
\caption{The $(\qdeg, \adeg)$ bidegrees of some anchored surfaces whose boundary consists of two non-contractible circles.}\label{fig:deg examples}
\end{figure}

\end{example}

\begin{lemma}
\label{lem:adeg is zero}
Let $S$ be an anchored surface. Then $\brak{S}=0$ or $\adeg(S) = 0$. 
\end{lemma}

\begin{proof}
If some component of $S$ has anchor points with different labels then $\brak{S} = 0$. Assume that all anchor points on any component of $S$ are labeled identically. We also assume that $S$ intersects $L$, otherwise $\adeg(S)=0$ is immediate. As usual, order the anchor points $p_1,\ldots, p_{m}$ from bottom to top. 

Take a generic half-plane $P$ in $\RR^3$ containing the anchor line $L$, so that $P\cap S$ consists of finitely many arcs (with boundary on $L$) and circles (disjoint from $L$). For any arc $a$ in $P\cap S$ with boundary $\d a = \{p_i, p_j\}$, necessarily $i$ and $j$ have opposite parities, and moreover $\l(p_i) = \l(p_j)$ by assumption. Therefore the total contribution of the anchor points $p_i$ and $p_j$ to $\adeg(S)$ is zero. Summing over all arcs in $P\cap S$ yields the statement of the lemma. 

\end{proof}

The subspace $\ker((,))\subset \Fr(C)$ respects this bigrading on $\Fr(C)$. Consequently, the bigrading descends to the state space $\brak{C}$.

Note that the relations \eqref{eq:neck cutting} and \eqref{eq:neck cutting line} are bi-homogeneous. Let $S\in \Fr(C)$ be a basis element of the form $S = S_1 \sqcup S_2$ where $S_1, S_2 \in \Fr(C)$ are anchored surfaces with $S_2$ closed. Then in $\brak{C}$ we have 
\begin{equation}
\label{eq:remove closed surface}
[S] = \brak{S_2} [S_1],  \ \  \brak{S_2}\in R_{\alpha}.
\end{equation}
Moreover, the relation \eqref{eq:remove closed surface} is bi-homogeneous. That it  is homogeneous with respect to $\qdeg$ follows from the fact that $\brak{S_2}\in \Ra$ is a polynomial of degree $\qdeg(S_2)$. Lemma \ref{lem:adeg is zero} ensures that $\adeg(S_2) = \adeg(\brak{S_2}) = 0$, so $\adeg(S) = \adeg(S_1)$. 

Given a bigraded module $M= \bigoplus_{(i,j) \in \Z^2} M_{i,j}$ over a commutative domain such that each $M_{i,j}$ has finite rank, define its \emph{graded rank} to be 
\[
\grank(M) = \sum_{i,j} \rank(M_{i,j}) q^i a^j. 
\]

\begin{lemma}
\label{lem:state space for circle}
Let $C\subset \P$ be a single circle. Then the state space $\brak{C}$ is a free $\Ra$-module of rank $2$. Moreover, we have
 \[
\grank(\brak{C}) = 
\begin{cases}
q+q^{-1} & \text{ if } C \text{ is contractible} \\
a+a^{-1} & \text{ if } C \text{ is non-contractible}. 
\end{cases}
\]
\end{lemma}

\begin{proof}

\begin{figure}
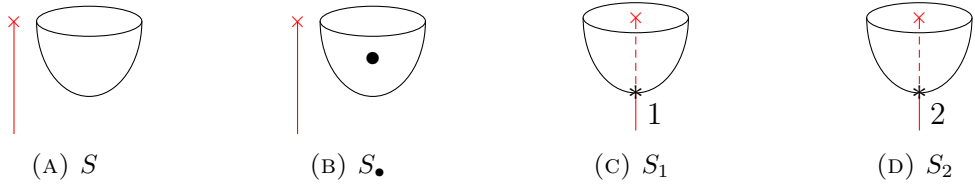

\centering
\subcaptionbox{ $S$
\label{fig:cup}}[.22\linewidth]
{\includestandalone{cup}
}
\subcaptionbox{ $S_\bullet$
\label{fig:dotted cup}}[.22\linewidth]
{\includestandalone{dotted_cup}
}
\subcaptionbox{ $S_1$
\label{fig:cup1}}[.22\linewidth]
{\includestandalone{cup1}
}
\subcaptionbox{ $S_2$
\label{fig:cup2}}[.22\linewidth]
{\includestandalone{cup2}
}
\caption{Basis elements for the state space of a single circle $C$. The first two surfaces form a basis if $C$ is contractible, and the last two form a basis if $C$ is non-contractible.}\label{fig:basis for circle}
\end{figure}

We consider two cases. If $C$ is contractible, then by applying the neck-cutting relation \eqref{eq:neck cutting} near $C$ and evaluating closed anchored surfaces as in equation \eqref{eq:remove closed surface}, we see that $\brak{C}$ is spanned by the two elements $S$ and $S_\bullet$ shown in Figures \ref{fig:cup} and \ref{fig:dotted cup}. Bidegrees of $S$ and $S_\bullet$ are $(-1,0)$ and $(1,0)$, respectively. Computing the matrix of the bilinear form \eqref{eq:bilinear form} for these elements yields
\[
\begin{pmatrix}
\b{S}S & \b{S}S_\bullet \\
\b{S_\bullet}S & \b{S_\bullet}S_\bullet 
\end{pmatrix} = 
\begin{pmatrix}
0 & 1 \\
1 & E_1
\end{pmatrix},
\]
which is invertible, thus $S, S_\bullet$ constitute a basis for $\brak{C}$. 

Now suppose $C$ is non-contractible. Applying the neck-cutting relation \eqref{eq:neck cutting line} near $C$ and evaluating closed anchored surfaces shows that the two elements $S_1, S_2$ depicted in Figures \ref{fig:cup1} and \ref{fig:cup2} span $\brak{C}$. Bidegrees of $S_1$ and $S_2$ are $(0,1)$ and $(0,-1)$, respectively. The matrix of the bilinear form is  
\[
\begin{pmatrix}
\b{S_1}S_1 & \b{S_1}S_2 \\
\b{S_2}S_1 & \b{S_2}S_2
\end{pmatrix} =
\begin{pmatrix}
1 & 0 \\
0 & 1
\end{pmatrix},
\]
hence $S_1, S_2$ are linearly independent and constitute a basis of $\brak{C}$. 
\end{proof}

\begin{theorem}
\label{thm:state space sl2}
Let $C\subset \P$ consist of $n$ contractible circles and $m$ non-contractible circles. Then the state space $\brak{C}$ is a free $R_\al$-module of rank $2^{n+m}$. Moreover, we have 
\[
\grank(\brak{C}) = (q+q^{-1})^n(a+a^{-1})^m. 
\]
\end{theorem}

\begin{proof}
Consider a $2^{n+m}$ element set $B(C)$ of basis vectors of $\Fr(C)$ consisting of surfaces $S$ satisfying 
\begin{itemize}
\item Each component of $S$ is a disk.
\item Each disk in $S$ with contractible boundary is disjoint from $L_-$ and carries either zero or one dot.
\item Each disk in $S$ with non-contractible boundary intersects $L_-$ exactly once, and its intersection point may be labeled by either $1$ or $2$. 
\end{itemize}
That $B(C)$ spans $\brak{S}$ follows from applying the two neck-cutting relations \eqref{eq:neck cutting} and \eqref{eq:neck cutting line} near the circles in $C$ and evaluating closed anchored surfaces. Linear independence of $B(C)$ and the statement regarding graded rank follow from the computations in Lemma \ref{lem:state space for circle}. 
\end{proof}

Elements of the basis $B(C)$ constructed above are \emph{standard generators}. For such a $\Sigma \in B(C)$ with $d$ dots and anchor points labeled $\l_1, \ldots, \l_m$, we have 
\begin{equation}
\label{eq:bidegree of generators}
\qdeg(\Sigma)=-n+2d, \hskip1em \adeg(\Sigma) = \sum_{i=1}^m (-1)^{i + \l(i)}. 
\end{equation}

Let $C_0, C_1 \subset \P$ be two collections of disjoint circles in the punctured plane. An \emph{anchored cobordism from $C_0$ to $C_1$} is a smoothly and properly embedded compact surface $S\subset \RR^2\times [0,1]$ with boundary $\d S = C_0 \sqcup C_1$, such that $C_i\subset \RR^2\times \{i\}$, $i=0,1$. Moreover, $S$ is required to intersect the arc $L_{[0,1]}:=(0,0) \times [0,1]$ transversely and come equipped with a labeling of these intersection points (called \emph{anchor points}), which is a map 
\[
\l = \l_S \: S\cap L_{[0,1]} \to \{1,2\}
\]
from the set of its anchor points to $\{1,2\}$. Anchored cobordisms are allowed to carry dots which can float on components but cannot jump to a different component. 

We compose anchored cobordisms in the usual manner. For anchored cobordisms $S_1 \:C_0\to C_1$, $S_2\: C_1\to C_2$, let $S_2 S_1 \: C_0 \to C_2$ denote the anchored cobordism obtained by gluing along the common boundary $C_1$ and re-scaling. Labels of anchor points of $S_2S_1$ are inherited from labels of $S_1$ and $S_2$.  

As above, if an anchored cobordism $S$ from $C_0$ to $C_1$ has $m$ anchor points and carries $d$ dots, define 
\[
\qdeg(S) = -\chi(S) + 2d + m.
\]
Let $\l(1), \ldots, \l(m)$ denote the labels of anchor points of $S$, ordered from bottom to top, and let $n$ be the number of non-contractible circles in $C_0$. Set 
\[
\adeg(S) =  (-1)^{n} \sum_{i=1}^m (-1)^{i+\l(i)}. 
\]

\begin{remark}
\label{rmk:cobordisms from empty set}
If $C_0=\varnothing$, then $S$ is a basis element of $\Fr(C_1)$, and moreover the two degrees $\qdeg(S)$, $\adeg(S)$ defined above for anchored cobordisms agree with the definitions in \eqref{eq:qdeg} and \eqref{eq:adeg} for elements of $\Fr(C_1)$. 
\end{remark}

An anchored cobordism $S$ from $C_0$ to $C_1$ induces an $\Ra$-linear map 
\[
S \: \Fr(C_0) \to \Fr(C_1)
\]
defined on the basis by gluing along the common boundary $C_0$. The definition of state spaces via universal construction immediately implies that we have an induced map
\begin{equation}
\label{eq:cobordism map}
\brak{S} \: \brak{C_0} \to \brak{C_1}.
\end{equation}

\begin{lemma}
\label{lem:bidegree}
Let $S_1 \: C_0\to C_1$ and $S_2 \: C_1\to C_2$ be anchored cobordisms. Then 
    \[
    \qdeg(S_2S_1) = \qdeg(S_2) + \qdeg(S_1), \hskip1em  \adeg(S_2S_1) = \adeg(S_2) + \adeg(S_1).
    \]
In particular, $\brak{S_1} \: \brak{C_0} \to \brak{C_1}$ is a map of bidegree $(\qdeg(S_1), \adeg(S_1))$. 

\end{lemma}

\begin{proof}
The first equality involving $\qdeg$ is straightforward. Let $n$ and $m$ denote the number of non-contractible circles in $C_0$ and $C_1$ respectively, and let $k$ denote the number of anchor points of $S_1$. We have 
\[
\adeg(S_2S_1) = \adeg(S_1) + (-1)^{n+m+k} \adeg(S_2).
\]
Note $n+m+k$ is even, since it is equal to the number of anchor points of the closed surface obtained by gluing disks to all boundary circles of $S_1$. 

The final statement concerning the bidegree of $\brak{S_1}$ follows from interpreting generators of $\brak{C_0}$ as anchored cobordisms $\varnothing \to C_0$, as in Remark \ref{rmk:cobordisms from empty set}. 
\end{proof}

\begin{definition}
\label{def:annular cobordism}
An \emph{annular cobordism} is an anchored cobordism $S\subset \RR^2\times [0,1]$ which is disjoint from the arc $L_{[0,1]}=(0,0) \times [0,1]$. An \emph{elementary} annular cobordism is one with a single non-degenerate critical point with respect to the height function $\RR^2\times [0,1] \to [0,1]$.
\end{definition} 

\begin{figure}
\centering 
\subcaptionbox{Type A \label{fig:type1}}[.22\linewidth]
{\includegraphics{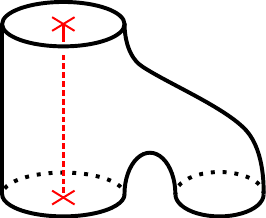}
}
\subcaptionbox{Type B \label{fig:type2}}[.22\linewidth]
{\includegraphics{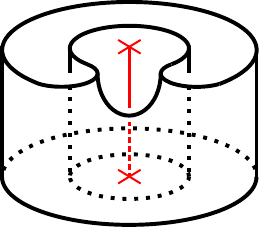}
}
\subcaptionbox{Type C \label{fig:type3}}[.22\linewidth]
{\includegraphics{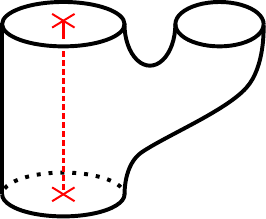}
}
\subcaptionbox{Type D \label{fig:type4}}[.22\linewidth]
{\includegraphics{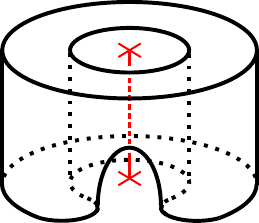}
}
\caption{Elementary saddles involving non-contractible circles.}\label{fig:elementary saddles}
\end{figure}

Elementary annular cobordisms consist of a union of a product cobordism with a cup, cap, or saddle. Every annular cobordism may be obtained by composing finitely many elementary ones. Cup and cap annular cobordisms always have contractible boundary. There are four types of elementary annular saddles involving at least one non-contractible circle, illustrated in Figure \ref{fig:elementary saddles}. In the next four examples we write down the maps assigned to these four cobordisms in the standard bases of state spaces, as defined in the proof of Theorem \ref{thm:state space sl2}. We also use the notation of shifted dots from \eqref{eq:shifted dot}. 

\begin{example}\label{ex:type1}
Figure \ref{fig:type1} map. The calculation for this map follows at once from the skein relation \eqref{eq:dot relations}. 
\begin{equation}
\begin{aligned}
     \includestandalone{type1_computation}   
\end{aligned}
\end{equation}
\end{example}

\begin{example}\label{ex:type2}
Figure \ref{fig:type2} map. This calculation follows easily from the skein relation \eqref{eq:moving cup off line}. 
\begin{equation}
\begin{aligned}
     \includestandalone{type2_computation}   
\end{aligned}
\end{equation}
\end{example}

\begin{example}\label{ex:type3}
Figure \ref{fig:type3} map. A convenient way to perform this calculation is to use neck-cutting with shifted dots \eqref{eq:neck cutting shifted dots} near the contractible circle and then simplify using the relations \eqref{eq:shifted dots relations}. 
\begin{equation}
\begin{aligned}
     \includestandalone{type3_computation}   
\end{aligned}
\end{equation}
\end{example}

\begin{example}\label{ex:type4}
Figure \ref{fig:type4} map. The neck-cutting relation \eqref{eq:neck cutting line} is helpful here. For the dotted cup we also use \eqref{eq:dot relations} to simplify further. 
\begin{equation}
\begin{aligned}
     \includestandalone{type4_computation}   
\end{aligned}
\end{equation}
\end{example}

Recall the involution $\tau$ of $R_{\alpha}$ that transposes $\alpha_1,\alpha_2$ and extend it to an antilinear involution, also denoted $\tau$, of the free $R_{\alpha}$-module $\Fr(C)$ as follows. Involution $\tau$ on $\Fr(C)$ sends a surface $S$ to the same surface with the labeling $\l$ of anchor points reversed and acts on linear combinations by   
\[ 
\tau\left(\sum_i \lambda_i S_i \right) = \sum_i \tau(\lambda_i)\tau(S_i).  
\] 
For a closed surface $S$ we have, by direct computation, $\brak{\tau(S)}= \tau(\brak{S})$, showing compatibility of the two involutions. If $S$, in addition, carries shifted dots, involution $\tau$ reverses their labels, so that $\tau(\circled{\mathsmaller{1}}) = \circled{\mathsmaller{2}}$ and $\tau(\circled{\mathsmaller{2}}) = \circled{\mathsmaller{1}}$. Involution $\tau$ descends to an involution, also denoted $\tau$, on $\brak{C}$. Annular degree is negated under $\tau$, $\adeg(\tau(S)) = - \adeg(S)$, for an anchored cobordism $S$.

\subsection{Annular link homology}

Let $\ACob$ denote the category whose objects consist of collections of finitely many disjoint simple closed curves in the punctured plane $\P$. A morphism from $C_0$ to $C_1$ in $\ACob$ is an anchored cobordism from $C_0$ to $C_1$, up to ambient isotopy fixing the boundary point-wise and mapping $L_{[0,1]}$ to itself. Let $\ACob'$ denote the subcategory of $\ACob$ with the same objects as $\ACob$ but whose morphisms are isotopy classes of annular cobordisms, disjoint from the anchor line $L$. The composition of annular cobordisms is again annular.

Let $\Ra\ggmod$ denote the category of bigraded $\Ra$-modules and homogeneous maps (of any bidegree) between them. We have a functor  
\[
\brak{-} \: \ACob \to \Ra\ggmod
\]
which sends a collection of circles $C\subset \P$ to the state space $\brak{C}$ and sends an anchored cobordism $S$ from $C_0$ to $C_1$ to the map $\brak{S} \: \brak{C_0} \to \brak{C_1}$ as in \eqref{eq:cobordism map}. By Lemma \ref{lem:bidegree}, $\brak{S}$ is a map of bidegree $(\qdeg(S), \adeg(S))$. We can restrict to the category of annular cobordisms to get a functor 
\[
\brak{-}' \: \ACob' \to \Ra \ggmod,
\]
which assigns to an annular cobordism $S$ a map $\brak{S}' = \brak{S}$ of bidegree $(\qdeg(S), 0)$. The restriction $\brak{-}'$ does not change the state space assigned to a collection of circles $C\subset \P$. 

On the other hand, a functor 
\[
\Ga \: \ACob' \to \Ra\ggmod
\]
was introduced in \cite{Akh}. We briefly recall $\Ga$ below.

Consider the algebra 
\[
\Aa = \Ra[X]/((X-\al_1)(X-\al_2)).
\]
It is a free $\Ra$-module with basis $\{1, X\}$. The trace $\epsilon_\al \: \Aa \to \Ra$ given by $1\mapsto 0$, $X\mapsto 1$ makes $\Aa$ into a Frobenius algebra, which defines a $(1+1)$-dimensional TQFT, a functor $\Fa$ from the category of dotted cobordisms to the category of $\Ra$-modules. A dot on a cobordism is interpreted as multiplication by $X\in \Aa$. Define a grading on $\Aa$ by setting 
\begin{equation}
\label{eq:A alpha grading}
\qdeg(1) = -1, \hskip1em \qdeg(X) = 1.
\end{equation}
With this grading, a cobordism $S$ with $d$ dots is assigned by $\Fa$ a map of degree $-\chi(S) + 2d$. Alternatively, the TQFT $\Fa$ is the result of applying the universal construction to the closed surface evaluation \eqref{eq:evaluation} when restricted to surfaces disjoint from $L$ and collections of contractible circles in $\P$. See \cite{KRfrobext} for further details about the Frobenius pair $(\Ra, \Aa)$. 

Let $C\subset \P$ be a collection of $n$ contractible and $m$ non-contractible circles. Define the bigraded $\Ra$-module $\Ga(C)$ as follows. As an $\Ra$-module, we set $\Ga(C) = \Fa(C) = \Aa^{\o (n+m)}$. Define the $\emph{annular grading}$, denoted $\adeg$, on $\Fa(C)$ as follows. 

Every tensor factor $\Aa$ corresponding to a contractible circle is concentrated in annular degree zero. Order the non-contractible circles in $C$ from outermost (furthest from the puncture) to innermost. Introduce the notation 
\begin{equation}
\label{eq:bases for essential circles}
  \begin{aligned}
  v_0 &= 1, & v_1 &= X- \al_1, \\
  v_0' &= 1, & v_1' &= X-\al_2.
  \end{aligned}  
\end{equation}
Both $\{v_0, v_1\} = \{1, X-\al_1\}$ and $\{v_0', v_1'\} = \{1, X-\al_2\}$ constitute an $\Ra$-basis for $\Aa$. Set 
\begin{equation}
\adeg(v_0) = \adeg(v_0') = -1, \hskip1em \adeg(v_1) =\adeg(v_1') = 1. 
\end{equation}
The annular grading on non-contractible circle is defined by assigning the homogeneous basis $\{v_0, v_1\}$ or $\{v_0', v_1'\}$ to the corresponding tensor factor of $\Aa$ in an alternating manner with respect to nesting in $\P$, with the convention that the outermost circle is assigned $\{v_0, v_1\}$. 

It is convenient to distinguish between the modules assigned to different types of circles in $\P$. Let $V_\al$ and $V_\al'$ denote the $\Ra$-modules $\Aa$ with bases $\{v_0, v_1\}$ and $\{v_0', v_1'\}$, respectively. The notation $\Aa$ will be reserved for the module assigned to a contractible circle, with basis $\{1, X\}$.

The $\Ra$-module $\Ga(C)$ also carries a quantum grading $\qdeg$ inherited from \eqref{eq:A alpha grading}. Define a modified quantum grading $\qdeg'$ on $\Ga(C)$ by
\begin{equation}
\label{eq:modified qdeg}
    \qdeg' = \qdeg - \adeg.
\end{equation}
We will consider $\Ga(C)$ as a bigraded $\Ra$-module with bigrading $(\qdeg', \adeg)$. Bidegrees are recorded in Figure \ref{fig:modified bidegrees}.

\begin{figure}
\centering
\begin{center}
\begin{tabular}{  l | c | c | c | c | c | c } 
& $1$ & $X$ & $v_0$ & $v_1$ & $v_0'$ & $v_1'$ \\ 
\hline
$\qdeg'$ & $-1$ & $1$ & $0$ & $0$ & $0$ & $0$\\ 
\hline
$\adeg$ & $0$  & $0$ & $-1$ & $1$ & $-1$ & $1$\\ 
\end{tabular}
\end{center}
\caption{The $(\qdeg', \adeg)$-bidegrees of relevant elements, where $\{1, X\}$ is a basis for a contractible circle and $\{v_0, v_1\}, \{v_0', v_1'\}$ are bases for non-contractible circles.} \label{fig:modified bidegrees} 
\end{figure}

\begin{remark}
 The modified quantum grading $\qdeg'$ appears elsewhere in the literature and is more natural in the context of annular link homology. In \cite{GLWsl2} this grading was denoted $j'$. Similarly, the annular link homology defined in \cite{BPW} carries the modified quantum grading. 
\end{remark}

We now define $\Ga$ on annular cobordisms. For an annular cobordism $S\subset \RR^2\times [0,1]$, if the boundary of $S$ is contractible in $\P$ then $\Ga(S)=\Fa(S)$, where $\Fa$ is the TQFT corresponding to the Frobenius algebra $\Aa$ as above. Formulas for the maps assigned by $\Ga$ to the four elementary cobordisms in Figure \ref{fig:elementary saddles} are recorded below. If other essential circles are present, then due to parity the formulas may be slightly different from those below. To obtain the full set of formulas, one interchanges  $v_0\leftrightarrow v_0'$, $v_1\leftrightarrow v_1'$, and $\al_1 \leftrightarrow \al_2$.

\noindent\begin{minipage}{.5\textwidth}
\begin{equation}\label{eq:equiv formula1}
  \begin{split}
   V_\al \o A_\al & \rar{\hyperref[fig:type1]{\operatorname{(A)}}} V_\al \\  v_0 \o 1 & \mapsto v_0 \\ v_1\o 1 &\mapsto v_1 \\
    v_0 \o X & \mapsto \al_1 v_0 \\
    v_1 \o X & \mapsto \al_2 v_1
  \end{split}
\end{equation}
\end{minipage}%
\begin{minipage}{.38\textwidth}
\begin{equation}\label{eq:equiv formula2}
  \begin{split}
   V_\al \o V_\al' & \rar{\hyperref[fig:type2]{\operatorname{(B)}}} A_\al \\  v_0 \o v_0' &\mapsto 0 \\
v_1 \o v_0' & \mapsto X - \al_1 \\
v_0 \o v_1' & \mapsto X - \al_2  \\
v_1 \o v_1' & \mapsto 0
  \end{split}
\end{equation}
\end{minipage}

\vskip1em

\noindent\begin{minipage}{.49\linewidth}
\begin{equation}\label{eq:equiv formula3}
  \begin{split}
   V_\al & \rar{\hyperref[fig:type3]{\operatorname{(C)}}} V_\al \o A_\al \\
v_0 & \mapsto v_0 \o (X-\al_2)\\
v_1 & \mapsto v_1 \o (X-\al_1)
  \end{split}
\end{equation}
\end{minipage} 
\begin{minipage}{.5\linewidth}
\begin{equation}\label{eq:equiv formula4}
  \begin{split}
   A_\al & \rar{\hyperref[fig:type4]{\operatorname{(D)}}} V_\al\o V_\al' \\
1 &\mapsto v_0 \o v_1' + v_1 \o v_0'\\
X & \mapsto  \al_1 v_0 \o v_1' + \al_2 v_1 \o v_0'
  \end{split}
\end{equation}
\end{minipage}
\vskip1em
\noindent

\begin{theorem}\label{thm_functor_iso} 
The functors $\brak{-}' \: \ACob' \to \Ra\ggmod$ and $\mathcal{G}_\al \: \ACob' \to \Ra\ggmod$ are naturally isomorphic via bidegree-preserving maps. 
\end{theorem}

\begin{proof}
Let $C\subset \P$ be a collection of circles. We will define an $\Ra$-linear, bidegree preserving isomorphism $\Phi_C \: \brak{C} \to \Ga(C)$ and show that it is natural with respect to annular cobordisms. 

Let $n$ and $m$ denote the number of contractible and non-contractible circles in $C$, respectively. Fix an ordering $1, \ldots, n$ on the contractible circles in $C$. The $\Ra$-module $\Ga(C)$ is free with basis given by elements of the form 
\[
y_1\o \cdots y_n \o z_1 \o \cdots \o z_m,
\]
where each $y_i$ is in $\{1, X\}$, specifying a basis element of the $i$-th contractible circle, and each $z_j$ is in either $\{v_0, v_1\}$ or $\{v_0', v_1'\}$, depending on nesting, specifying basis elements of the non-contractible circles. The ordering of factors $z_1\o \cdots \o z_m$ corresponding to non-contractible circles is from outermost to innermost as usual, so that the first factor $z_1$ labels the outermost non-contractible circle.

We now define the isomorphism $\Phi_C \: \brak{C} \to \Ga(C)$. Recall the standard basis $B = B(C)$ for $\brak{C}$ defined in the proof of Theorem \ref{thm:state space sl2}. For $\Sigma \in B$ with anchor points labeled $\l_1, \ldots, \l_m$, read from bottom to top, set 
\[
\Phi_C(\Sigma) = y_1\o \cdots y_n \o z_1 \o \cdots \o z_m,
\]
where $y_i = 1$ if the corresponding cup in $\Sigma$ is undotted and $y_i = X$ if the corresponding cup in $\Sigma$ is dotted. The generators $z_j$ of non-contractible circles are determined using the rule
\[
z_j =
\begin{cases}
v_1 & \text{ if } j \text{ is odd and } \l_j = 1 \\
v_0 & \text{ if } j \text{ is odd and } \l_j = 2 \\
v_0' & \text{ if } j \text{ is even and } \l_j = 1 \\
v_1' & \text{ if } j \text{ is even and } \l_j = 2 
\end{cases}
\]
See Figure \ref{fig:isomorphism ex} for an example of the assignment $\Phi_C$ when $n=1$, $m=2$. By comparing the bidegree formula \eqref{eq:bidegree of generators} for $\Sigma$ with the bidegree of $\Phi_C(\Sigma)$ (see Figure \ref{fig:modified bidegrees}), we see that $\Phi_C$ is a bidegree-preserving isomorphism. Recall that we use the modified quantum grading \eqref{eq:modified qdeg} for $\Ga(C)$. 

\begin{figure}
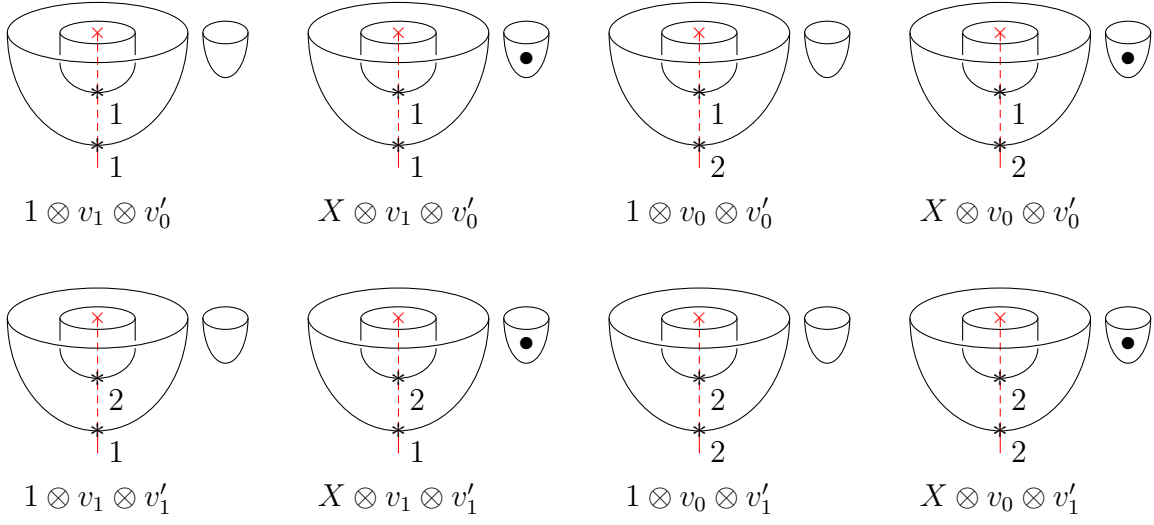

\centering
\includestandalone{basis_ex}
\caption{An example of the isomorphism $\Phi_C$ when $C$ consists of one contractible circle and two non-contractible circles. Basis elements $\Sigma$ of $\brak{C}$ are drawn with the corresponding basis element $\Phi_C(\Sigma) \in \Ga(C)$ written underneath.}\label{fig:isomorphism ex}
\end{figure}

Now let $S\:C_1 \to C_2$ be an annular cobordism. To complete the proof, we check that the square 
\[
\begin{tikzcd}
\brak{C_1} \ar[r, "{\Phi_{C_1}}"] \ar[d, "{\brak{S}}"'] & \Ga(C_1) \ar[d, "{\Ga(S)}"] \\
\brak{C_2} \ar[r, "{\Phi_{C_2}}"] & \Ga(C_2)
\end{tikzcd}
\]
commutes. If all the boundary circles of $S$ are contractible, then commutativity of the square is straightforward. Otherwise, if $S$ has at least one non-contractible boundary circle, it suffices to consider the case where $S$ is one of the elementary annular cobordisms depicted in Figure \ref{fig:elementary saddles}. Formulas for these maps were recorded in Examples \ref{ex:type1} -- \ref{ex:type4}. Comparing with the formulas \eqref{eq:equiv formula1} -- \eqref{eq:equiv formula4} completes the proof. 
\end{proof}

Let $\A : = S^1\times [0,1]$ denote the annulus. For an oriented link $L \subset \A \times [0,1]$ in the thickened annulus, a generic projection of $L$ onto $\A \times \{0\}$ yields a link diagram $D$ in the interior of $\A$. Identifying the interior of $\A$ with the punctured plane $\P$, we may form the cube of resolutions of $D$ in the usual way, with all smoothings drawn in $\P$. Diagrams representing isotopic annular links are related by Reidemeister moves away from the puncture. By standard arguments \cite{Kh, BNtangles}, applying the functor $\brak{-}'\: \ACob' \to \Ra \ggmod$ to the cube of resolutions yields a bigraded chain complex $\brak{D}$ whose chain homotopy type is an invariant of $L$. Theorem \ref{thm_functor_iso} implies that the resulting annular homology is isomorphic to that of \cite{Akh}.

%
%

\section{Unoriented \texorpdfstring{$SL(3)$}{SL(3)} anchored homology}
\label{sec:unoriented sl3}

We recall definitions and notations from~\cite{KRfoamev}, including that of  (unoriented) $SL(3)$ foams and refer the reader to~\cite[Section 2.1]{KRfoamev} for more details.

\begin{figure}
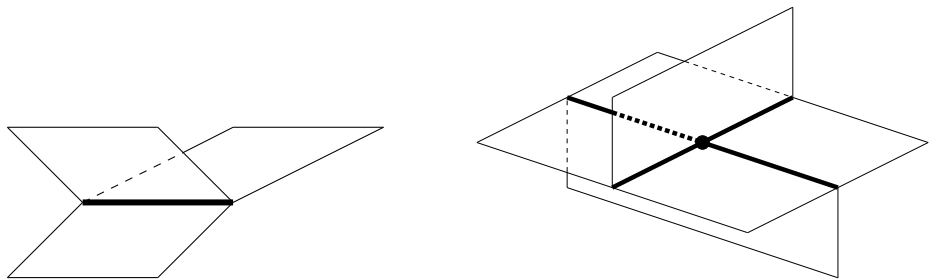

    \centering
\subcaptionbox{Seam points where three facets meet.\label{fig:seam line}}[.4\linewidth]
{\includestandalone{seam_line}
}
\subcaptionbox{A seam vertex where six facets meet.\label{fig:seam vertex}}[.4\linewidth]
{\includestandalone[scale=.6]{seam_vertex}
}
    \caption{Local model of a pre-foam near singular points. The singular graph $s(F)$ is drawn bold.}
    \label{fig:singular graph}
\end{figure}

\begin{definition}
A (closed) \emph{$SL(3)$ pre-foam} is a compact 2-dimensional CW complex equipped with a PL-structure such that each point has an open neighborhood that is either an open disk, the product of a tripod and an open interval (Figure \ref{fig:seam line}), or the cone over the $1$-skeleton of a tetrahedron (Figure \ref{fig:seam vertex}). Points of the first type are called \emph{regular}, those of the second are called \emph{seam points}, and those of the third are called \emph{seam vertices}. A (closed) \emph{$SL(3)$ foam} is a closed $SL(3)$ pre-foam together with a PL embedding into $\RR^3$. 
\end{definition}

We will simply write \emph{pre-foam} and \emph{foam} in place of \emph{closed} $SL(3)$ \emph{(pre-)foam}. For a pre-foam $F$, denote by $v(F)$ the set of seam vertices and by $s(F)$ the set of seam points and seam vertices. The subspace $s(F)$ is a $4$-valent graph which may contain closed loops. Connected components of $s(F) \setminus v(F)$ are called \emph{seams}. 

The subspace $F\setminus s(F)$ is a (not necessarily compact) surface, and a connected component of $F\setminus s(F)$ will be called a \emph{facet} of $F$. The (finite) set of facets of $F$ is denoted $f(F)$. Facets of pre-foams may be decorated by a finite number of dots, which are allowed to float freely on their facets but may not cross seams or enter seam vertices. 

\begin{figure}
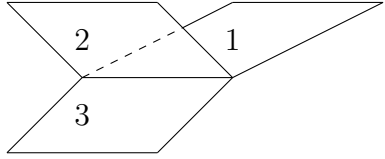

    \centering
    \includestandalone{coloring}
    \caption{The local model for a pre-admissible coloring near a seam point.}
    \label{fig:pre-admissible coloring}
\end{figure}

A \emph{coloring} of a pre-foam $F$ is a map 
\[
c \: f(F) \to \{1,2,3\}.
\]
That is, a coloring assigns $1, 2$ or $3$ to each facet of $F$. A coloring is called \emph{pre-admissible} if the three facets meeting at each seam of $F$ have distinct colors, see Figure \ref{fig:pre-admissible coloring}. For a pre-admissible coloring $c$ and $1\leq i, j \leq 3$, $i\not= j$, let $F_{ij}(c)$
denote the union of facets colored $i$ or $j$. The pre-admissibility condition guarantees that each $F_{ij}(c)$ is a closed surface (see \cite[Proposition 2.2]{KRfoamev}).

A coloring $c$ is called \emph{admissible} if each $F_{ij}(c)$ is orientable. For a foam $F$ (that is, a pre-foam embedded in $\RR^3$), every pre-admissible coloring is admissible, since $F_{ij}(c)$ is a closed surface in $\RR^3$. 

\subsection{Unoriented anchored \texorpdfstring{$SL(3)$}{SL(3)} foams and their evaluations} \label{subsec_unor_anch}

Fix a field $\kk$ of characteristic $2$. In this section the following commutative rings will be used. 
\begin{itemize}
\item $R_x' = \kk[x_1, x_2, x_3]$ is the ring of polynomials in three variables. 
\item $R_x=\kk[E_1, E_2, E_3]$ the subring of $R_x'$ that consists of symmetric polynomials in $x_1,x_2,x_3$, with generators $E_i$ being elementary symmetric polynomials: 
\begin{align*}
    E_1 &= x_1 + x_2 + x_3, \\
    E_2 &= x_1x_2 + x_1x_3 + x_2 x_3, \\
    E_3 &= x_1x_2x_3. 
\end{align*}
\item $R_x'' = R_x'[(x_1+x_2)^{-1}, (x_2+x_3)^{-1}, (x_1+x_3)^{-1}]$ is a localization of $R_x'$ given by inverting $x_i+x_j$, for $1\le i<j \le 3$. 
\item $\til{R}_x' = \kk[\sqrt{x_1}, \sqrt{x_2}, \sqrt{x_3}]$ is the extension of $R_x'$ obtained by introducing square roots of $x_1, x_2, x_3$. 
\item $\til{R}_x'' = \kk[\sqrt{x_1}, \sqrt{x_2}, \sqrt{x_3},(x_1+x_2)^{-1}, (x_2+x_3)^{-1}, (x_1+x_3)^{-1}]$ is a localization of $\til{R}_x'$ given by inverting $x_i+x_j$, for $1\le i<j \le 3$. 
\end{itemize}
All five of these rings are graded by setting $\deg(x_i) = 2$, $i=1,2,3$. Inclusions of the above rings are summarized in the following diagram. 

\begin{equation}
\label{eq:rings}
\begin{tikzcd}[row sep = .3em, column sep =.2em ]
& & \til{R}_x' & \subset & \til{R}_x'' \\
& & \cup & & \cup \\
R_x & \subset & R_x' & \subset & R_x''
\end{tikzcd}
\end{equation}

We follow the notation established in \cite{KRfoamev} for these rings with the additional subscript $x$ to distinguish from the notation in Section \ref{sec:anchored surfaces}.

\begin{definition}
\label{def:unoriented anchored foam}
An \emph{anchored $SL(3)$ foam} $F$ is an $SL(3)$ foam $F'\subset \RR^3$ that may intersect the line $L$ at finitely many points away from 
the singular graph $s(F')$ of $F'$. Thus each intersection point belongs to some facet $f$ of $F'$, and intersection of facets with $L$ are required to be transverse. Denote by $p(F) =F\cap L$ the set of intersection points (anchor points) of $F$. Intersection points carry labels in $\{1,2,3\}$; that is, $F$ comes equipped with a fixed map 
\[
\l \: p(F) \to \{1, 2, 3\}.
\]
It is convenient to order anchor points $p_1,\dots, p_m$ from  bottom to top, with labels
$\ell_i=\ell(p_i)$, $i=1,\dots, m$. 
\end{definition}

We now refine the notion of admissible coloring of a foam to that of admissible coloring of an anchored foam $F$. Consider an anchored foam $F$ with the underlying foam $F'$. A coloring $c\in \adm(F')$ induces a coloring of anchor points in $F'$, by assigning to each point the color of its facet. We say that $c$ is admissible if that's exactly the labeling of anchor points of $F$, that is, 
$\ell(p)=c(f)$ for each anchor point $p$ in a facet $f$, and then set $c(p)=\ell(p).$

In this way, the set of admissible colorings of $F'$ is in a bijection with the set of admissible colorings of anchored foams $F$ that become $F'$ upon forgetting the labeling of anchor points: 
\[  \adm(F') \ \cong \ \coprod_{F} \adm(F). 
\] 
Various constructions with $SL(3)$ foams in~\cite{KRfoamev} extend directly to anchored foams. In particular, bicolored surfaces $F_{ij}(c)$ are well-defined, associated to an admissible coloring $c$. We will also call an admissible coloring simply a coloring. We will use $i,j,k$ to denote the three elements of $\{1, 2,3\}$, not necessarily in that order. 

We refine \cite[Definition 2.9]{KRfoamev} for anchored foams. 

\begin{definition}
Let $F$ be an anchored foam, $c\in \adm(F)$ be an admissible coloring, and $\Sigma$ a connected component of $F_{ij}(c)$ which is disjoint from $L$. Define a coloring $c'$ of $F$ which swaps the colors $i$ and $j$ on facets of $\Sigma$, and leaves all other facets colored according to $c$. We say that $c$ and $c'$ are related by an \emph{$ij$-Kempe move along $\Sigma$}. Note that since $\Sigma$ has no anchor points, $c'$ is still an admissible coloring of $F$. 
\end{definition}

Kempe moves can be done on components $\Sigma$ of $F_{ij}(c)$ that intersect $L$ as well, but the resulting anchored foam $F_0$ is different from $F$ due to carrying different labels on anchor points on $\Sigma$. 

For $k\in \{1, 2, 3\}$, denote by $k', k''$ its two complementary elements, so that $\{k, k', k''\} = \{1, 2, 3\}$. Let $F$ be an anchored foam with labeling $\l$. Let $c\in \adm(F)$ be an admissible coloring. For an anchor point $p\in p(F)$ lying on a facet $f\in f(F)$, we set $c(p) = c(f)=\ell(p)$; that is, $c(p)$ is the color of the facet, according to $c$, on which $p$ lies, which equals $\l(p)$ since $c$ is admissible. For $1\leq i \leq 3$, let $d_i(c)$ denote the number of dots on facets colored $i$. For $1\leq i \neq j \leq 3$, let $F_{ij}(c)$ be the union of facets of $F$ colored $i$ or $j$. The space $F_{ij}(c)$ is a closed surface in $\RR^3$ and hence has even Euler characteristic. Set 
\begin{equation}
\label{eq:unoriented eval color}
\brak{F,c} = \frac{P(F,c)}{Q(F,c)},
\end{equation}
where
\begin{align}
P(F,c) &= \prod_{i=1}^3 x_i^{d_i(c)}\cdot \left( \prod_{p\in p(F)} (x_{c(p)} + x_{\l(p)'})(x_{c(p)} + x_{\l(p)''})   \right)^{1/2}, \label{eq:unoriented eval color numerator}\\
Q(F,c) &=  \prod_{1\leq i < j \leq 3} (x_i + x_j)^{\chi(F_{ij}(c))/2}. \label{eq:unoriented eval color denominator}
\end{align}
The product of the two terms under the square root, for a given anchor point $p$, is equal to 
\begin{align*}
(x_1+x_2)(x_1+x_3) \ \ &\mathrm{ if }\ \  c(p)=1, \\
(x_2+x_1)(x_2+x_3) \ \ &\mathrm{ if }\ \  c(p)=2, \\
(x_3+x_1)(x_3+x_2) \ \ &\mathrm{ if }\ \  c(p)=3 .
\end{align*} 
\begin{remark}
\label{rmk:square decoration}
This product is the inverse of the square decoration $\square$ in~\cite[Section 4.1]{KRfoamev}. The square decoration was used to study a separable version of the unoriented $SL(3)$ theory, with the discriminant  $\mathcal{D}=(x_1+x_2)(x_1+x_3)(x_2+x_3)$ inverted, which is a version of the Lee theory. Here, we use the defect line $L$ rather than freely floating square dots in~\cite[Section 4.1]{KRfoamev} in the opposite way, to add factors to the evaluation rather than divide by terms in the discriminant.  
\end{remark} 

\begin{remark}
\label{rmk:admissible colorings of anchored foams}
If $c$ is an admissible coloring of the underlying foam $F'$ of $F$ but not of the anchored foam $F$, then the evaluation \eqref{eq:unoriented eval color} is still defined and equal to zero:
\begin{equation}\label{eq_brak_is_0}
    \brak{F,c}=0, \ \  c\in\adm{F'}\setminus \adm{F}.
\end{equation}
This holds since, for some $p\in p(F)$, its color $c(p)$ differs from its label $\l(p)$, so that $x_{c(p)} + x_{c(p)}=0$ appears under the square root in \eqref{eq:unoriented eval color numerator} and $P(F,c)=0$. Thus, 
\begin{equation*}  
(x_{c(p)} + x_{\l(p)'})(x_{c(p)} + x_{\l(p)''})  = \begin{cases}
 (x_{\ell(p)} + x_{\l(p)'})(x_{\ell(p)} + x_{\l(p)''}) \ \ \mathrm{if} \ \ c(p)=\ell(p), \\
 0 \ \ \mathrm{otherwise}. 
\end{cases}
\end{equation*} 
\end{remark}

Define the \emph{evaluation of $F$} to be 
\begin{equation}
\label{eq:unoriented eval}
    \brak{F} = \sum_{c\in \adm(F)} \brak{F,c}.
\end{equation}
Alternatively, we can sum over the larger set of $c\in \adm(F')$, due to (\ref{eq_brak_is_0}). 

\vspace{0.1in} 

Let us explain the square root in equation \eqref{eq:unoriented eval color numerator}. The equality $\sqrt{x+y} =\sqrt{x} + \sqrt{y}$ holds in a commutative ring of characteristic $2$, so $\brak{F,c}$ is in the ring $\til{R}_x''$, see \eqref{eq:rings}. We will show in Proposition \ref{prop:eval has no square roots} that, in fact, no square roots appear, so that $\brak{F,c} \in R_x''$. Likewise, in Proposition \ref{prop:unoriented sl3 eval is polynomial} we show that $\brak{F} \in R_x'$. 

The evaluation \eqref{eq:unoriented eval} is multiplicative with respect to disjoint union and does not depend on a particular embedding of $F$ into $M=(\RR^3,L)$ as long as anchor points on $F$ and their labels are specified. 

If an anchored foam $F$ is a disjoint union of anchored foams $F_1 \sqcup \cdots \sqcup F_k$, then 
\[
\brak{F} = \brak{F_1} \cdots \brak{F_k}. 
\]
If $F$ is disjoint from $L$, then $\brak{F}$ is equal to the evaluation in \cite[Section 2.3]{KRfoamev}.

\begin{example}
\label{ex:sphere unoriented}
Let $F$ be a $2$-sphere $\SS^2$ with two anchor points and $d$ dots. Its evaluation is zero unless both points have the same label $i\in \{1, 2, 3\}$, in which case there is only admissible coloring $c$ which colors $F$ by $i$. Let $j, k\in \{1, 2, 3\}$ denote the complementary elements to $i$. The surfaces $F_{ij}(c), F_{ik}(c)$ are $2$-spheres, while $F_{jk}(c) = \varnothing$. Then the evaluation is 
\[
\brak{F} = \frac{x_i^d \left(( x_i + x_j)^2 (x_i+x_k)^2\right)^{1/2}}{(x_i + x_j)(x_i + x_k)} = x_i^d.
\]
\end{example}

\begin{example}
More generally, let $F$ be a genus $g$ surface carrying $d$ dots and $2n>0$ anchor points. It evaluates to zero unless all points are labeled by the same $i\in \{1,2,3\}$. In this case, letting $j,k\in \{1,2,3\}$ be the complementary elements to $i$, the evaluation is 
\[
\brak{F} = \frac{x_i^d \left( (x_i+x_j)(x_i+x_k) \right)^{n}}{\left( (x_i+x_j)(x_i+x_k) \right)^{1-g}}=x_i^d \left( (x_i+x_j)(x_i+x_k) \right)^{n+g-1} .
\]
\end{example}

\begin{example}
Consider the theta foam $F$ whose facets each intersect $L$ once, with anchor points labeled $i,j,k\in \{1,2,3\}$ and facets carrying $d_1, d_2$, and $d_3$ dots, as shown in \eqref{eq:theta foam}. 
\begin{equation}
    \begin{aligned}\label{eq:theta foam}
    \includestandalone{theta_foam}
    \end{aligned}
\end{equation}

In an admissible coloring of the underlying foam, the three facets must have distinct colors, so $\brak{F}=0$ if $i,j,k$ are not distinct. If $i,j,k$ are distinct, then there is one admissible coloring $c$ which colors the top, middle, and bottom facets, respectively, by $i,j$, and $k$. The surfaces $F_{ij}(c), F_{ik}(c), F_{jk}(c)$ are $2$-spheres, and the evaluation is
\[
\brak{F} = x_i^{d_1} x_j^{d_2} x_k^{d_3}. 
\]

\end{example}

\begin{remark}
Note that the evaluation of an anchored foam is in general not a symmetric function in $x_1, x_2, x_3$, whereas in \cite{KRfoamev} the evaluation is always an element of $R_x$. 
\end{remark}

Let us call a sequence $\ell \in \{1,2,3\}^m$ \emph{pre-admissible} if the following holds. Let $u_1,u_2,u_3$ be three non-zero elements of the abelian group $\Z/2\times \Z/2$. Sequence $\ell$ is \emph{pre-admissible} iff 
\begin{equation}\label{eq_usum} 
\sum_{i=1}^m u_{\ell(i)} =0 \in \Z/2\times \Z/2. 
\end{equation}

\begin{proposition}\label{prop_pre_adm} If an anchored foam $F$ has an admissible coloring, the sequence  $\ell$ of its anchor points is pre-admissible. 
\end{proposition} 
\begin{proof}
 Consider a generic intersection of $F$ with a half-plane in $\RR^3$ bounding $L$. This intersection is a trivalent graph $\Gamma$ in the half-plane. Coloring $c$ of $F$ induces a coloring $c'$ of edges of $\Gamma$ such that around each trivalent vertex of $\Gamma$ the colors of the three edges are distinct (Tait coloring). On the boundary points (one-valent vertices) of $\Gamma$ the coloring is given by labeling $\ell$. The sum on the left hand side of (\ref{eq_usum}) is zero since it can alternatively be written as the sum of triples of vectors $u_1+u_2+u_3=0$ over all trivalent vertices of $\Gamma$. Each inner edge of $\Gamma$, that bounds two trivalent vertices, contributes $u_i+u_i=0$ to the sum, where $i$ is the color of the edge. An edge with one or both endpoints on the boundary contributes the sum of $u_i$'s, over its boundary points. 
\end{proof}

For an anchored foam $F$ and $1\leq i \leq 3$, let $\anch(i)$ denote the number of anchor points of $F$ with label $i$ (the dependence on $F$ is omitted). 

\begin{proposition}
\label{prop:eval has no square roots}
 For an anchored foam $F$ and an admissible coloring $c$, we have $\brak{F,c}\in R_x''$.
\end{proposition}
\begin{proof}
Recall the rings $R_x''$ and $\til{R}_x''$ defined in \eqref{eq:rings}. It's clear that $\brak{F,c}$ belongs to the larger ring $\til{R}_x''$.

The expression in \eqref{eq:unoriented eval color} under the square root is equal to 
\[
(x_1 + x_2)^{\anch(1) + \anch(2)} (x_2+ x_3)^{\anch(2) + \anch(3)} (x_1 + x_3)^{\anch(1) + \anch(3)}.
\]

For $1\leq i < j \leq 3$, the integer $\anch(i) + \anch(j)$ is even since it is equal to the number of intersection points of the closed surface $F_{ij}(c)$ with $L$, see also Proposition~\ref{prop_pre_adm}. Consequently, taking the square root produces integral exponent of $x_i+x_j$, implying that $\brak{F,c}$ is in $R_x''$. 
\end{proof} 

Using the above notation, the square root term in \eqref{eq:unoriented eval color numerator} is equal to 
\begin{equation}
\til{Q}(F,c) := \prod_{1\leq i < j \leq 3} (x_i+x_j)^{(\anch(i) + \anch(j))/2},
\end{equation}
so that formula \eqref{eq:unoriented eval color} can be rewritten as 
\begin{equation}
    \brak{F,c} = \prod_{i=1}^3 x_i^{d_i(c)} \prod_{1\leq i < j \leq 3} (x_i+x_j)^{(\anch(i) + \anch(j) - \chi(F_{ij}(c)))/2}.
\end{equation}

\begin{proposition}
\label{prop:unoriented sl3 eval is polynomial}
 For an anchored foam $F$, we have $\brak{F}\in R_x'=\kk[x_1,x_2,x_3].$
\end{proposition}

\begin{proof}
    Proof of Theorem 2.17 in~\cite{KRfoamev} extends with minor changes to this case. Note that the evaluation is no longer a symmetric function. We must show that positive powers of $x_i + x_j$, $1\leq i < j \leq 3$, do not appear in the denominator of $\brak{F}$. Let us specialize to $i=1,j=2$. Denominators $x_1+x_2$ in the evaluations $\brak{F,c}$ may appear only from the components of $F_{12}(c)$ that are two-spheres. If a 2-sphere does  not intersect $L$, the proof in~\cite{KRfoamev} works in this case as well. Suppose a 2-sphere component $\Sigma$ of $F_{12}(c)$ intersects $L$ in $\anch(1)$ points colored $1$ and $\anch(2)$ points colored $2$ (necessarily in the corresponding facets of $F$ carrying those colors under $c$). These points contribute 
    \[ (x_1+x_2)^{\anch(1) + \anch(2)}(x_1+x_3)^{\anch(1)}(x_2+x_3)^{\anch(2)} 
    \] 
    to the expression under the square root, and $\anch(1)+\anch(2)\ge 2$, allowing to cancel the denominator term $x_1+x_2$ that $\Sigma$ contributes. Summing over all admissible colorings and otherwise following the arguments in~\cite[Theorem 2.17]{KRfoamev} implies the result. 

\end{proof}

\begin{remark}
 Contributions of anchor points to  the evaluation $\brak{F,c}$ can  be  interpreted as follows. Consider polynomial $f(x)=(x-x_1)(x-x_2)(x-x_3)\in R_x'[x]$. Then 
 \[ f'(x)=(x-x_2)(x-x_3) + (x-x_1)(x-x_3)+ (x-x_1)(x-x_2)
 \]
 and 
 \begin{eqnarray*}
     f'(x_1) &  = & (x_1-x_2)(x_1-x_3), \\
     f'(x_2) &  = & (x_2-x_1)(x_2-x_3), \\
     f'(x_3) &  = & (x_3-x_1)(x_3-x_2).  
 \end{eqnarray*}
 Contribution of an anchor point $p$ with a label $i=\ell(p)$ to the evaluations $\brak{F,c}$ and $\brak{F}$ is then $\sqrt{f'(x_i)}$, the square root of the derivative of $f$ at the root $x_i$ of the polynomial $f$. In characteristic two signs do not matter, but this observation hints how to extend the evaluation to characteristic $0$. 
\end{remark}

Since the labels $i_1,\dots, i_m$ of anchor points are fixed in a given $F$, these marked points contribute the same term, 
\[ \sqrtLF \ := \ \left(\prod_{r=1}^m f'(x_{i_r})\right)^{1/2}
\] 
and we have 
\begin{equation}
    \brak{F,c} = 
    \sqrtLF
    \cdot \brak{F',c},  \ \ \ \brak{F} = 
    \sqrtLF
    \cdot \brak{F'},
\end{equation}
where $F'$ is the foam $F$ viewed as a regular foam with anchored points and their labels ignored. When coloring $c$ of $F$ is not compatible with labels of anchor points, though, we should define $\sqrtLF=0$ to match the formula $\brak{F,c}=0$. 

Also notice that, switching to characteristic $0$ and from the matrix factorization viewpoint \cite{KRoz}, $f(x)=w'(x)$ is the derivative of the potential $w(x)=x^4/4 -E_1 x^3/3 + E_2 x^2/2 - E_3 x $, so that the contributions of anchor points are given by square roots of the second derivative $\sqrt{w''(x_i)}$ at critical points of $w$, analogous to the square root of the Hessian factor that appears, for example, in the steepest descent method formulas. 

%

\subsection{Skein relations}\label{subsec_skein} 

In this subsection we record several local relations satisfied by the evaluation of anchored $SL(3)$ foams. We start with the following proposition concerning the relations in \cite[Section 2.5]{KRfoamev}, which should be understood as occurring away from the anchor line $L$.

\begin{proposition}
\label{prop:KR local relations}
The twelve local relations in Propositions 2.22 -- 2.33 from \cite{KRfoamev} hold. 
\end{proposition}

\begin{proof}
The arguments in \cite{KRfoamev} apply without modification. 
\end{proof}

We will use shifted dots in this section, as in \eqref{eq:shifted dot}. For $1\leq i \leq 3$, we allow anchored foams to carry decorations of the form $\circled{i} = \bullet + x_i$ on a facet. They are required to be disjoint from $L$, float freely on their facets, but cannot move past seams or seam vertices. 
\begin{equation}
    \begin{aligned}
    \includestandalone{shifted_dot_unoriented_sl3}
    \end{aligned}
\end{equation}
For an anchored foam $F$ carrying $\circled{i}$ on some facet $f\in f(F)$, any coloring $c\in \adm(F)$ which colors $f$ by $i$ evaluates to zero, $\brak{F,c} = 0$. An anchor point labeled $i$ has the same effect as placing $\sqrt{\circled{j}\ \circled{k}}=\sqrt{\circled{i'}\ \circled{i''}}$ on the facet on which it lies (recall our conventions that $\{1,2,3\}=\{i,j,k\}=\{i,i',i''\}$). See also equation \eqref{eq:cup off line unoriented} and the discussion in Section \ref{sec:Remark on Lee's theory}. 

We also have relations involving the anchor line. 

\begin{lemma}
\label{lem:unoriented sl3 relations}
The local relations \eqref{eq:neck cutting line unoriented}, \eqref{eq:dot relation unoriented}, \eqref{eq:handle sum of pts}, \eqref{eq:cup off line unoriented},  and \eqref{eq:line past seam} hold.
\begin{equation}
\label{eq:neck cutting line unoriented} 
    \begin{aligned}
    \includestandalone{neck_cutting_line_unoriented_sl3}
    \end{aligned}
\end{equation}

\begin{equation}
\label{eq:dot relation unoriented} 
    \begin{aligned}
    \includestandalone{dot_relation_unoriented}
    \end{aligned}
\end{equation}

\begin{equation}
\label{eq:handle sum of pts} 
    \begin{aligned}
    \includestandalone{handle_sum_of_pts}
    \end{aligned}
\end{equation}

\begin{equation}
\label{eq:cup off line unoriented} 
    \begin{aligned}
    \includestandalone{moving_cup_off_line_unoriented_sl3}
    \end{aligned}
\end{equation}

\begin{equation}
\label{eq:line past seam} 
    \begin{aligned}
    \includestandalone{line_past_seam_unoriented}
    \end{aligned}
\end{equation}
\end{lemma}

In the last two equations, $\{i,j,k\}=\{1,2,3\}$.

\begin{proof}
Let us verify equation \eqref{eq:neck cutting line unoriented}; the other three relations are easier to check and the proof is left to the reader. Denote by $F$ the anchored foam on the left-hand side, and by $G^1, G^2, G^3$ the three foams on the right-hand side, with the superscript corresponding to the labels of the depicted anchor points. For $1\leq i \leq 3$, let $\adm_i(F)$ be the set of admissible colorings of $F$ in which the depicted tube is colored by $i$. Admissible colorings of $G^i$ must color the two disks by $i$, so there is a natural bijection $\adm_i(F) \cong \adm(G^i)$. 

For $c\in \adm_i(F)$, let $c'\in  \adm(G^i)$ denote the corresponding coloring. We will show that
\[
\brak{F,c} = \brak{G^i, c'},
\]
which completes the proof. 

 The anchored foam $G^i$ carries two more anchor points, both labeled $i$, than $F$ does, while the dot placement for $G^i$ and $F$ is the same, so
\[
P(G,c') = (x_i + x_j)(x_i + x_k) P(F,c), 
\]
where $\{i,j,k\}=\{1,2,3\}$. 
On the other hand, 
\[
\chi(G^i_{ij}(c')) = \chi(F_{ij}(c)) + 2, \hskip1em \chi(G^i_{ik}(c')) = \chi(F_{ik}(c))+ 2, \hskip1em \chi(G^i_{jk}(c')) = \chi(F_{jk}(c)),
\]
which yields 
\[
Q(G,c') = (x_i + x_j)(x_i+x_k) Q(F,c).
\]
Thus $\brak{F,c} = \brak{Q,c'}$ as desired. Summing over all admissible colorings of $F$ we get 
\[ \brak{F} = \brak{G^1}+\brak{G^2}+\brak{G^3},
\]
completing the proof. 
\end{proof}


\subsection{State spaces}
\label{sec:state spaces unoriented}
We generalize the notion of webs and cobordisms between them from \cite[Section 3.1]{KRfoamev} in the presence of the anchor line $L$.

\begin{definition}
\label{def:unoriented anchored foam with boundary}
A \emph{web} is a trivalent graph $\Gamma$ which is PL-embedded into the punctured plane $\P = \RR^2 \setminus \{(0,0)\}$. We allow webs to have closed loops with no vertices. A \emph{anchored foam with boundary} $V$ is the intersection of a closed anchored foam $F\subset \RR^3$ with a thickened plane $\RR^2 \times [0,1]$ such that $F\cap (\P\times \{i\})$, $i=0,1$ is a web (in particular, $F$ is disjoint from the two points $(0,0,0)$ and $(0,0,1)$), and dots of $F$ are disjoint from $\RR^2\times \{i\}$ and from $L$. Foams with boundary are considered equivalent if there is an orientation-preserving homeomorphism of $\RR^2\times [0,1]$ taking one to the other which fixes the boundary of $\RR^2\times [0,1]$ pointwise and maps the line segment $L_{[0,1]} := \{(0,0)\} \times [0,1]$ to itself.

For a foam with boundary $V$, let 
\[
p(V) = V \cap L_{[0,1]}
\]
denote its intersection points with the anchor line, called \emph{anchor points}. Each anchor point is required to carry a label in $\{1,2,3\}$. 
\end{definition}

We view $V$ as a cobordism from the web $\d_0 V := V \cap (\RR^2 \times \{0\})$ to the web $\d_1 V := V \cap (\RR^2 \times \{1\})$. A closed foam is then a cobordism from the empty web to itself. We will often refer to foams with boundary simply as \emph{foams} when the meaning is clear from context. Composition $WV$ of foams $V, W$ with $\d_1 V = \d_0 W$ is defined in the natural way. We obtain a category $\AFoam$ of webs and anchored foams.

The category $\AFoam$ has a contravariant involution $\omega$ which is the identity on webs and which sends a foam to its reflection about $\RR^2\times \{1/2\}$, preserving the labels of anchor points. As for closed foams, denote by $s(V)$ and $v(V)$ the singular graph and singular vertices, respectively, of a foam with boundary $V$. Define the \emph{degree} of $V$ to be 
\begin{equation}
\label{eq:deg of foam}
\deg(V) = 2 \big( \lr{d(V)}  +  \lr{p(V)} - \chi(V) \big) - \chi(s(V)) ,
\end{equation}
where $d(V)$ is the set of dots on $V$.

The definition of admissible colorings extend naturally to anchored foams with boundary. An admissible coloring induces a Tait coloring on the boundary webs. If a foam with boundary $V$ has an admissible coloring $c$, then by \cite[Remark 2.8]{KRfoamev},
\begin{equation}
    \deg(V) = 2 \lr{d(V)}  + 2 \lr{p(V)} - \big( \chi(V_{12}(c)) + \chi(V_{13}(c)) + \chi(V_{23}(c)) \big).
\end{equation}
It follows that for a closed foam $F$, its evaluation $\brak{F} \in R_x'$ is a homogeneous polynomial of degree $\deg(F)$. 

\begin{lemma}
For composable foams $V$ and $W$, we have 
\[
\deg(WV) = \deg(W) + \deg(V).
\]
\end{lemma}
\begin{proof}
This follows from \cite[Proposition 3.1]{KRfoamev} and  $\lr{p(WV)} = \lr{p(W)} + \lr{p(V)}$.
\end{proof}

We now define state spaces for webs via universal construction and the evaluation formula \eqref{eq:unoriented eval}. For a web $\Gamma$, let 
\[
\Fr(\Gamma)
\]
denote the free $R_x'$-module generated by all anchored foams $V$ from the empty web to $\Gamma$. Define a bilinear form 
\[
(-,-) \: \Fr(\Gamma) \times \Fr(\Gamma) \to R_x'
\]
by $(V,W) = \brak{\omega(V)W}$. This bilinear form is symmetric since $\brak{F} = \brak{\omega(F)}$ for any closed anchored foam $F$. Define the state space $\brak{\Gamma} := \Fr(\Gamma)/\ker((-,-))$ to be the quotient of $\Fr(\Gamma)$ by the kernel 
\[
\ker((-,-)) = \{ x\in \Fr(\Gamma) \mid (x,y) = 0 \text{ for all } y\in \Fr(\Gamma) \}
\]
of the bilinear form. Note that $(-,-)$ is degree-preserving, so its kernel and the state space $\brak{\Gamma}$ are graded $R_x'$-modules. 

An anchored foam $V \: \Gamma_0 \to \Gamma_1$ naturally induces a map
\[
\brak{V} \: \brak{\Gamma_0} \to \brak{\Gamma_1}
\]
of degree $\deg(V)$, defined by sending the equivalence class of a basis element $U\in \Fr(\Gamma_0)$ to the class of the composition $VU$. This is functorial with respect to composition of anchored foams, $\brak{WV} = \brak{W} \brak{V}$ for composable anchored foams with boundary $V$ and $W$. 

\begin{remark}
\label{rmk:no Tait colorings}
For a web $\Gamma$ and basis elements $V_1, V_2\in \Fr(\Gamma)$, an admissible coloring of the closed foam $\omega(V_2)V_1$ induces a Tait coloring of $\Gamma$. Thus $\brak{\Gamma} = 0$ if $\Gamma$ has no Tait colorings; see also \cite[Proposition 3.16]{KRfoamev}. 
\end{remark}

\begin{figure}
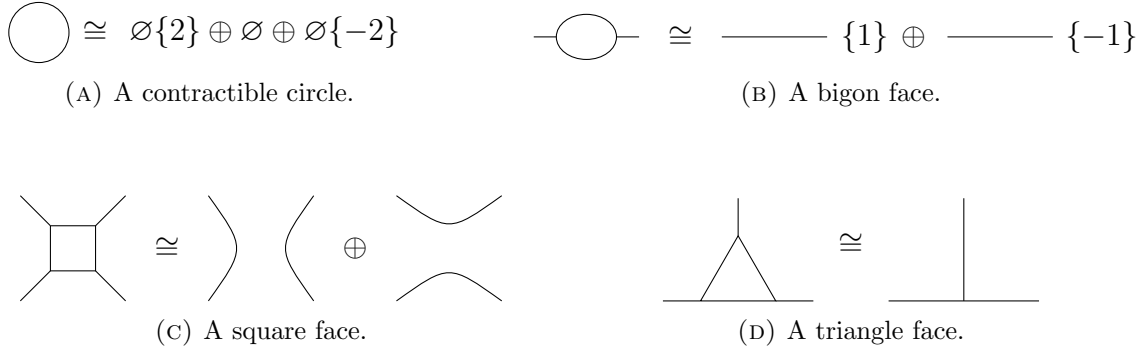

\centering 
\subcaptionbox{A contractible circle.}[.4\linewidth]
{\includestandalone{local_iso_1}
} \quad \quad 
\subcaptionbox{A bigon face.}[.5\linewidth]
{\includestandalone{local_iso_2}
}\\ \vskip6ex
\subcaptionbox{A square face.}[.4\linewidth]
{\includestandalone{local_iso_3}
}\quad \quad \quad 
\subcaptionbox{A triangle face.}[.4\linewidth]
{\includestandalone{local_iso_4}
}
\caption{Direct sum decompositions from \cite[Section 3.3]{KRfoamev}, where the depicted regions do not contain the puncture.}\label{fig:local contractible isos}
\end{figure}

\begin{proposition}
\label{prop:local isos}
The local\footnote{Here local means that the webs involved in the isomorphisms are identical outside of a disk which is disjoint from the puncture, and in this disk they are related as in the figures accompanying the statements of the propositions.} isomorphisms in Propositions 3.12 -- 3.15 from \cite{KRfoamev}, also shown in Figure \ref{fig:local contractible isos}, hold. 
\end{proposition}

\begin{proof}
Proposition \ref{prop:KR local relations} guarantees that the explicit isomorphisms defined in \cite{KRfoamev} hold in the anchored setting as well.
\end{proof}

\begin{proposition}
\label{prop:noncontractible circle}
Let $\Gamma \subset \P$ be a web with a non-contractible circle $C$ which bounds a disk in $\RR^2\setminus \Gamma$, and let $\Gamma' = \Gamma \setminus C$ be the web obtained by removing $C$. Then there is an isomorphism 
\[
\brak{\Gamma}  \cong \brak{\Gamma'} \oplus \brak{\Gamma'} \oplus \brak{\Gamma'}
\]
given by the maps shown in \eqref{eq:noncontractible iso}.
\begin{equation}
\label{eq:noncontractible iso}
    \begin{aligned}
    \includestandalone{noncontractible_iso}
    \end{aligned}
\end{equation}
\end{proposition}

\begin{proof}
This follows from Example \ref{ex:sphere unoriented} and the relation \eqref{eq:neck cutting line unoriented}. Note that there are no grading shifts in the three copies of $\brak{\Gamma'}.$ 
\end{proof}

It is an interesting and nontrivial problem to identify the state spaces $\brak{\Gamma}$. In the construction in \cite{KRfoamev} without the anchor line, state spaces can be simplified using the relations in \cite[Section 3.3]{KRfoamev}, see Figure \ref{fig:local contractible isos}. In particular, bipartite webs always contain a contractible circle, bigon, or square, so the state space in the bipartite case is a free module of graded rank equal to the Kuperberg bracket~\cite{Kup}, normalized as in \cite{Khsl3}; see also \cite[Proposition 3.17, Proposition 4.15]{KRfoamev}. The simplest web which cannot be simplified using the relations in Figure~\ref{fig:local contractible isos} and for which the state space is unknown is the dodecahedral graph, as explored in \cite{Boozer,KRconicalfoams}, and, on the gauge theory side, in~\cite{KM2,KM1,KM3}. 

One may also ask to identify state spaces in the presence of the anchor line and the modified evaluation considered in this paper. Propositions \ref{prop:local isos} and \ref{prop:noncontractible circle} give some ways to simplify state spaces. In general, we are not able to decompose the bigon, square, and triangle regions in Figure \ref{fig:local contractible isos} if they contain the puncture. An extended evaluation, obtained by introducing additional types of intersection points of $L$ with a foam, is discussed in Section \ref{sec:new intersection points}. The following lemma addresses reducibility of smallest webs. 

\begin{lemma}
\label{lem:graphs} 
Let $\Gamma \subset \RR^2$ be a connected, planar, trivalent graph with no edges connecting a vertex to itself \footnote{A graph with such an edge has trivial state space, see Remark \ref{rmk:no Tait colorings}.}.  

\begin{enumerate}
    \item  If $\Gamma$ is bipartite, then $\Gamma$ has at least two bounded faces with at most four edges each. 
    
    \item If at most one of the bounded faces of $\Gamma$ has fewer than five edges, then $\Gamma$ has at least eight vertices. 
\end{enumerate}
\end{lemma}

\begin{proof}
Let $v,e, f$ denote the number of vertices, edges, and faces (including the unbounded face) of $\Gamma$, respectively. Label the faces $1, \ldots, f$, and for $1\leq i \leq f$, let $r_i$ denote the number of edges comprising the boundary of the $i$-th face. We have 
\begin{equation}
\label{eq:faces and vertices formula}
    \sum_{i=1}^f r_i = 2e = 3v,
\end{equation} 
where the second equality holds since $\Gamma$ is trivalent. 

We first prove statement $(1)$. Since $\Gamma$ is bipartite, each $r_i$ is even. Suppose for the sake of contradiction that at most one bounded face of $\Gamma$ has four or fewer edges. Then equation \eqref{eq:faces and vertices formula} implies
\[
\sum_{i=1}^f r_i > 6(f-2),
\]
so $12 > 6f-3v$. On the other hand, an Euler characteristic computation gives 
\[
12 = 6(f-e+v) = 6f-3v,
\]
which is a contradiction. 

Let us now address statement $(2)$. From equation \eqref{eq:faces and vertices formula} we obtain
\[
3v \geq 5(f-2) + 4 = 5f -6
\]
since, by assumption, there are $(f-2)$ faces with at least $5$ edges each, and the remaining two faces each have at least two edges. This together with an Euler characteristic computation gives $f\geq 6$, and it follows that $v\geq 8$. 
\end{proof}

\begin{corollary}
Let $\Gamma \subset \P$ be a bipartite web. Then $\brak{\Gamma}$ is a free $R_x'$-module of rank equal to the number of Tait colorings of $\Gamma$.
\end{corollary}

\begin{proof}
By statement (1) of Lemma \ref{lem:graphs}, any such web has either an innermost non-contractible circle or a region, not containing the puncture, which either bounds a closed loop, or is a bigon or square face. Thus state space can be reduced using Propositions \ref{prop:local isos} and \ref{prop:noncontractible circle}, and since the resulting web remains bipartite we can continue the procedure. 
\end{proof}

\begin{figure}
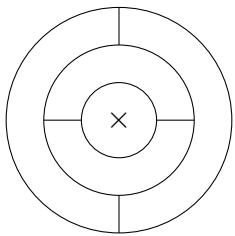

    \centering
    \includestandalone{graph}
    \caption{The simplest non-reducible web in the punctured plane.}
    \label{fig:non reducible web}
\end{figure}

It is natural to ask what is the simplest web for which the state space cannot be reduced using Propositions~\ref{prop:local isos} and~\ref{prop:noncontractible circle}. By Statement (2) of Lemma \ref{lem:graphs}, such a web has at least eight vertices. The web shown in Figure~\ref{fig:non reducible web} has precisely eight vertices and cannot be simplified using our local relations. 
We have not identified the state space of this web, but it can be approached via the 4-periodic  (and, in general, non-exact) complex described in~\cite[Section 4.3]{KRfoamev}. It can be applied along any of the four edges of Figure~\ref{fig:non reducible web} web near either the marked or the infinite point. One of the other three webs in the complex contains a loop  and has trivial homology, but additional computations are needed to identify the state space due to non-exactness of the complex.   

\vspace{0.1in} 

An annular graph $\Gamma \subset \P$ is called \emph{reducible} if its state space can be reduced to a sum of those for the empty annular graph by recursively applying relations (A)-(D) in Figure~\ref{fig:local contractible isos} and relation in Proposition~\ref{prop:noncontractible circle}. It may make sense to also allow reductions to annular graphs without Tait colorings (including graphs with loops), since such graphs have trivial state spaces. 

A reducible annular graph allows an identification of its state space with a suitable free graded $R_x$-module by recursively applying the above state sum decompositions.  
As a special case, we have the following decomposition formula for collections of simple closed curves in an annulus. 

\begin{proposition}
Let $\Gamma \subset \P$ consist of $n$ contractible circles and $m$ non-contractible circles. Then the state space $\brak{\Gamma}$ is a free $R_x'$-module of graded rank $3^m(q^2+1+q^{-2})^n$. 
\end{proposition}

In particular, for a reducible $\Gamma$, the graded rank of the free $R_x'$-module $\brak{\Gamma}$ can be computed recursively.

Anchored foams and state spaces carry an additional $(\Z/2 \times \Z/2)$-grading as follows. Recall that $u_1, u_2, u_3$ denote the nonzero elements of $\Z/2 \times \Z/2$. For a foam $V$ with (possibly empty) boundary, define
\[
\adeg(V) = \sum_{p\in p(V)} u_{\l(p)}.
\]
We call $\adeg$ the \emph{annular degree}. Clearly $\adeg$ is additive under disjoint union and composition.

The annular degree extends to a $(\Z/2\times \Z/2)$-grading on $\Fr(\Gamma)$, for a web $\Gamma \subset \P$,  by setting the ground ring $R_x'$ to be concentrated in annular degree zero. Proposition \ref{prop_pre_adm} implies that $\brak{F} = 0$ or $\adeg(\brak{F}) = 0$ for any closed foam $F$. It follows that $(-,-)$ preserves annular degree, so $\adeg$ descends to a $(\Z/2 \times \Z/2)$-grading on the state space $\brak{\Gamma}$. The annular grading is the unoriented version of the grading on state spaces  of annular oriented 
webs by the integral weight lattice of $\mathfrak{sl}_3$, see  Section~\ref{sec:oriented sl3 state spaces}, even though the action of the latter is lacking on the equivariant annular state spaces.

\vspace{0.1in} 

In~\cite[Section 4]{KRfoamev} the authors consider localization of the unoriented $SL(3)$ theory given by  inverting the discriminant $\mcD=(x_1+x_2)(x_1+x_3)(x_2+x_3)$. This localization results in a significant simplification of the theory, making it separable, so to speak. In particular, a suitable 4-term sequence of web state spaces in~\cite[Section 4.3]{KRfoamev} is exact. 

This localization easily extends to the annular case. The corresponding 4-term sequence is exact in the annular case as well. The ground ring for that theory is $R'_\mcD:=\kk[x_1,x_2,x_3,\mcD^{-1}]$, with $\kk$ a characteristic two field. The analogue of~\cite[Proposition 4.13]{KRfoamev} holds: the localized state space of an annular web $\Gamma$ is a  projective $R'_\mcD$-module of rank equal to the number of Tait colorings of $\Gamma$. The latter is the number of edge colorings of $\Gamma$ into three colors so that at each vertex the colors are distinct. Proof of this result in~\cite{KRfoamev} easily adapts to the annular case, with the modification that the region around the marked point can be inductively simplified, if necessary, by reducing to the other three terms in the exact sequence, until it has a single edge (a loop around the marked point).

\vspace{0.1in} 

%
%

\subsection{Remark on Lee's theory}
\label{sec:Remark on Lee's theory}
Recall the function 
\begin{equation}
    f(x)=(x+x_1)(x+x_2)(x+x_3)= x^3 + E_1 x^2 + E_2 x + E_3
\end{equation} 
(in characteristic  $2$ signs do not matter)
with coefficients in the ring $R_x$ and roots in $R'_x\supset R_x$. One can form the quotient ring 
$A:= R'_x[x]/(f(x))$, naturally isomorphic to the homology of a contractible circle in our theory. 
Let 
\begin{equation}
\mcD=(x_1+x_2)(x_1+x_3)(x_2+x_3) = E_1E_2+E_3
\end{equation}
be the discriminant. Consider the localization 
\begin{equation}
    R'_{\mcD} := R'_x[\mcD^{-1}], \ \  
    A_{\mcD}  := R'_{\mcD}\otimes_{R'_x} A. 
\end{equation}
Introduce idempotents $e_1,e_2,e_3\in A_{\mcD}$: 
\begin{equation}
    e_i \ := \ \frac{(x+x_j)(x+x_k)}{(x_i+x_j)(x_i+x_k)}, \ \ \{i,j,k\} =\{1,2,3\}. 
\end{equation}
We have 
\begin{equation}
    1 = e_1 + e_2 + e_3 , \ \ e_i e_j = \delta_{i,j} e_i. 
\end{equation}
These idempotents decompose ring $A_{\mcD}$ into the direct product 
\begin{equation}
    A_{\mcD}\cong R'_{\mcD}e_1 \times R'_{\mcD}e_2 \times R'_{\mcD}e_3 \cong 
    R'_{\mcD}\times R'_{\mcD} \times R'_{\mcD} .
\end{equation}

An idempotent $e_i$ can be visualized as floating on a facet of a foam $F$, in the localized theory. These idempotents allow to decompose an evaluation of a foam $F$ with $n$ facets into $3^n$  terms by summing over all ways to place each of these three idempotents onto facets of $F$. Each term is straightforward to compute and equals zero unless the idempotents define a Tait coloring (an admissible coloring) of $F$. 

Idempotent $e_i$ bears  a close relation to an anchor point labeled $i$. The anchor point $p$ on a facet $f$ contributes the term $\sqrt{f'(x_{c(f)})}=\sqrt{(x_{c(f)}+x_j)(x_{c(f)}+x_k)}$ to the evaluation $\brak{F,c}$. The square of this term is either $0$ (if $i\not=c(f)$) or the denominator of $e_i$, if $i=c(f)$, for any coloring $c$ of $F$.

Comparing $e_i$ and an anchor point $p$ labeled $i$, when coloring $c$ associates color $c(f)\not=i$ to the facet $f$ carrying $e_i$ or $p$, both evaluations are zero. When $c(f)=i$, the idempotented dot $e_i$ contributes $1$ to the evaluation, while the anchor point contributes$\sqrt{f'(x_{i})}$. Denominator of $e_i$ is  $f'(x_i)$.

\vspace{0.1in} 

One can try to unify $e_i$ and anchor points $p$ by considering anchor lines and circles $L$ in $\RR^3$ possibly intersecting a foam $F$. 
Intersection points (anchor points) carry labels $i\in \{1,2,3\}$ and a circle anchor points labeled $i$ is the idempotent $e_i$. 
Then a ``small" circle intersecting a facet $f$ at two points, both labeled $i$, can also be converted into $e_i$. 
Notice that once $e_i$ are allowed, integrality is lost and an evaluation of such a foam may contain denominators which are products of $x_i+x_j$.   

\vspace{0.1in} 

For a different generalization, instead of a single line $L\subset \RR^3$ consider a 1-manifold $L$ property embedded in $\RR^3$, say a finite union of lines and circles, possibly knotted.  All anchor points (intersection points with $L$) on a foam $F$ carry labels, with the usual contribution to the evaluation, as in formula (\ref{eq:unoriented eval color numerator}). Integrality Theorem~\ref{thm:oriented sl3 eval is polynomial} still holds for such generalized evaluation. In particular, given $k$ points on a plane, one can define various state spaces for webs $\Gamma$ embedded in the plane and disjoint from these marked points.  Also note that for $k\geq 2$ punctures, bipartite graphs are in general not reducible, which makes it harder to understand corresponding state spaces in the oriented $SL(3)$ case.

\vspace{0.1in} 

\begin{remark}
A handle next to but disjoint from an anchor line can be written as a sum of three lower genus terms intersecting the line, see equation (\ref{eq:handle sum of pts}), 
which follows from the formula
\[ m\circ \Delta(1) = (x_1+x_2)(x_1+x_3)+(x_1+x_2)(x_2+x_3) +  (x_1+x_3)(x_2+x_3) = f'(x_1)+ f'(x_2)+f'(x_3). 
\]
\end{remark}

\vspace{0.1in} 

\subsection{Unlabeled anchor points and bigon decomposition}
\label{sec:new intersection points}

Direct sum decompositions for webs $\Gamma$ containing a bigon, triangle, or square face which do not contain the puncture are given in Proposition \ref{prop:local isos}. On the other hand, Proposition \ref{prop:noncontractible circle} describes how to simplify a web containing an innermost non-contractible circle. In order to have direct sum decompositions for more general regions containing the puncture we introduce additional types of intersections of the anchor line $L$ with a foam and modify the evaluation $\brak{-}$.

\begin{figure}
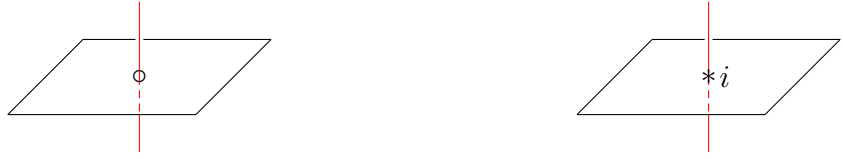

\centering    
\subcaptionbox{A Type 1 anchor point marked $\circ$ and carrying no label. \label{fig:type 1 point}}[.45\linewidth]
{\includestandalone{type1_intersection} 
}
\subcaptionbox{A Type 2 anchor point marked $*$ with label $i\in \{1,2,3\}$. \label{fig:type 2 point}}[.45\linewidth]
{\includestandalone{type2_intersection}
}
\caption{Type 1 and Type 2 anchor points.}
\label{fig:type 1 and 2 anchor points}
\end{figure}

In addition to anchor points, which carry labels in $\{1,2,3\}$ as in Definition \ref{def:unoriented anchored foam}, we allow finitely many transverse intersections of $L$ with a foam $F$ away from the singular graph $s(F)$, and we do not require labels. We will call the usual (labeled) anchor points \emph{Type 2}, and the new (unlabeled) anchor points \emph{Type 1}. In the figures, we denote Type 2 anchor points by an asterisk $*$ as usual, along with a label in $\{1, 2, 3\}$, and Type 1 anchor points will be indicated by a small unshaded circle $\circ$. Figure \ref{fig:type 1 and 2 anchor points} illustrates the convention. Let $p_1(F)$ and $p_2(F)$ denote the set of Type 1 and Type 2 anchor points, respectively (using the notation in Section \ref{subsec_unor_anch}, $p(F) = p_2(F)$). The definition of admissible coloring remains the same.

We modify the evaluation in the presence of Type 1 points as follows.  Let $c\in \adm(F)$. For $p\in p_1(F)$ lying on some facet $f\in f(F)$, let $c(p) := c(f)$ denote the coloring of the facet on which $p$ lies. Also recall that for $i\in \{1, 2, 3\}$, we write $i', i''$ and $j,k$ to denote the two complementary elements, so that $\{1,2,3\} = \{i, j, k\} = \{i, i', i''\}$.

Define 
\begin{align}
\til{Q}_\circ(F,c) &= \prod_{p\in p_1(F)} \sqrt{x_{c(p)'} + x_{c(p)''} }, \label{eq:type 1 contribution} \\
P_\circ(F,c) &= P(F,c) \cdot \til{Q}_\circ(F,c) , \\
\brak{F,c}_\circ &= \frac{P_\circ(F,c)}{Q(F,c)},\\
\brak{F}_\circ &= \sum_{c\in \adm(F)} \brak{F,c}_\circ. 
\end{align}
where $P(F,c)$ and $Q(F,c)$ are as defined in \eqref{eq:unoriented eval color numerator} and \eqref{eq:unoriented eval color denominator}. In other words, a Type 1 point $p$ on an $i$-colored facet contributes a factor of $\sqrt{x_j + x_k}$ to the evaluation $\brak{F,c}_\circ$. 

\begin{remark}
Type 1 intersection points are related to the triangle decoration from \cite[Section 4.1]{KRfoamev}. Precisely, the contribution of a Type 1 point $p$ to the square root in \eqref{eq:type 1 contribution} equals the inverse of placing a triangle decoration on the facet where $p$ lies. See relation \eqref{eq:cup off line type 1}, as well as Remark \ref{rmk:square decoration} for a related discussion.
\end{remark}

Note that Type 1 intersection point contributes half the degree of a Type 2 point to the degree of the evaluation and, thus, to the degree of a cobordism represented by a foam with boundary.  
\vspace{0.1in}

\begin{example}
\label{ex:sphere type 1}
Consider a $2$-sphere $F$ carrying $d$ dots and intersecting $L$ in two Type 1 anchor points, as shown in \eqref{eq:sphere type 1}.

\begin{equation}
\label{eq:sphere type 1}
    \begin{aligned}
    \includestandalone{sphere_type_1}
    \end{aligned}
\end{equation}

For $1\leq i \leq 3$, let $c_i \in \adm(F)$ color $F$ by $i$. Then 
\begin{align*}
    \brak{F,c_i}_\circ &= \frac{x_i^d(x_j + x_k)}{(x_i+ x_j)(x_i + x_k)}\ ,\\
    \brak{F}_\circ &= \brak{F,c_1}_\circ + \brak{F,c_2}_\circ + \brak{F,c_3}_\circ \\
    & = \frac{ x_1^d (x_2 + x_3)^2 + x_2^d (x_1 + x_3)^2 + x_3^d (x_1 + x_2)^2 }{ (x_1 + x_2)(x_1 + x_3)(x_2 + x_3) } \\
    & = \frac{ x_1^d (x_2^2 + x_3^2) + x_2^d (x_1^2 + x_3^2) + x_3^d (x_1^2 + x_2^2)}{ (x_1 + x_2)(x_1 + x_3)(x_2 + x_3) } .
\end{align*}
Thus $\brak{F}_\circ = 0$ if $d=0, 2$, and $\brak{F}_\circ = 1$ if $d=1$. For $d\geq 3$, the last expression above equals the ratio of the antisymmetrizer with exponent $(d,2,0)$ and antisymmetrizer with exponent $(2,1,0)$ (up to adding signs, which does not matter in characteristic $2$). Thus $\brak{F}_\circ$ equals the Schur function $s_\lambda(x_1, x_2, x_3)$ for the partition $\lambda = (d-2, 1,0)$ when $d\geq 3$.

\end{example}

\begin{example}
\label{ex:sphere mixed types}
Consider a $2$-sphere $F$ carrying $d$ dots and intersecting $L$ in one Type 1 anchor point and one Type 2 anchor point, as shown in \eqref{eq:sphere mixed types}.
\begin{equation}
\label{eq:sphere mixed types}
    \begin{aligned}
    \includestandalone{sphere_mixed_types}
    \end{aligned}
\end{equation}
Then $F$ has one admissible coloring, and 
\[
\brak{F}_\circ = \frac{x_i^d  \sqrt{(x_i + x_j)(x_i+x_k)(x_j+x_k)}}{(x_i + x_j)(x_i+x_k)} = \frac{x_i^d \sqrt{x_j+x_k}}{\sqrt{(x_i + x_j)(x_i+x_k)}}.
\]
\end{example}

From Example \ref{ex:sphere mixed types} we see that the evaluation $\brak{F}_\circ$ in general has denominators and square roots, so we can only conclude that   
\[\brak{F}_\circ \in\til{R}_{\circ} := \kk[x_1, x_2, x_3,(x_1+x_2)^{-1/2}, (x_2+x_3)^{-1/2}, (x_1+x_3)^{-1/2}].
\] 
Note that $\til{R}_{\circ}$ is a subring of $\til{R}''_x$, see Section~\ref{subsec_unor_anch} and diagram (\ref{eq:rings}).

We use $\til{R}_{\circ}$ as the ground ring of the theory. Evaluations of closed anchored foams $F$ with two types of anchor points belong to this ring. We define the state space $\brak{\Gamma}_{\circ}$ of a trivalent graph $\Gamma \subset \P$ using this evaluation and following the general recipe of Section~\ref{sec:state spaces unoriented}. The state space is a graded $\til{R}_{\circ}$-module, but, due to the presence of invertible elements $(x_i+x_j)^{1/2}$ of degree $1$, grading carries little information, and for many purposes one can downsize and consider the degree zero part $\brak{\Gamma}^0_{\circ}$ of the state space, which is a module over the degree $0$ subring $\til{R}^0_{\circ}$ of $\til{R}_{\circ}$. 

This theory  is functorial and foams with top and bottom boundary and anchor points of those two different types induce maps between the corresponding state spaces. Various  direct sum decompositions  that hold for the unoriented $SL(3)$ theory $\brak{-}$ hold for this theory as well.  

We also have local relations involving Type 1 intersection points.

\begin{lemma}
\label{lem:type1 relations}
The local relations\footnote{To clarify relation \eqref{eq:type1 relation1}: the first term on the right-hand side of the equality has a Type 1 anchor point on each of two front-facing half-bubbles, while the second term has a Type 1 anchor point on each of the two back-facing half-bubbles.} \eqref{eq:cup off line type 1}, \eqref{eq:type1 relation1}, \eqref{eq:type1 relation2}, and \eqref{eq:type1 relation3} hold for the theory $\brak{-}_\circ$.

\begin{equation}
\label{eq:cup off line type 1}
    \begin{aligned}
    \includestandalone{moving_cup_off_line_type1}
    \end{aligned}
\end{equation}

\begin{equation}
\begin{aligned}
    \label{eq:type1 relation1}
    \includestandalone{type1_relation1}
\end{aligned}
\end{equation}

\begin{equation}
\begin{aligned}
    \label{eq:type1 relation2}
    \includestandalone{type1_relation2}
\end{aligned}
\end{equation}

\begin{equation}
\begin{aligned}
    \label{eq:type1 relation3}
    \includestandalone{type1_relation3}
\end{aligned}
\end{equation}
\end{lemma}

\begin{proof}
Relation \eqref{eq:cup off line type 1} is straightforward and left to the reader. Let us verify relation \eqref{eq:type1 relation1}. Denote by $F$ the foam on the left-hand side of the equality, and denote by $F^1$ and $F^2$ the two foams on the right-hand side. There is a natural identification $\adm(F^1) = \adm(F^2)$. 

Let $c\in \adm(F^1)$ be a coloring in which the front two half-bubble facets are differently colored, say the top front half-bubble is colored $j$, the bottom front half-bubble is colored $k$, and the remaining ``big'' facet is colored $i$. Continue to denote by $c\in \adm(F^2)$ the corresponding coloring of $F^2$. The top Type 1 intersection point of $F^1$ contributes $\sqrt{x_i+ x_k}$ to $\brak{F^1,c}$ and the bottom Type 1 intersection point of $F^1$ contributes $\sqrt{x_i + x_j}$, while the contributions of these points to $\brak{F^2,c}$ are reversed. Thus in characteristic two we have 
\[
\brak{F^1, c} + \brak{F^2,c} = 0. 
\]

Next, the admissible colorings of $F$ are in natural bijection with the admissible colorings of $F^1$ (and of $F^2$) in which the front half-bubbles of $F^1$ are colored the same. Let $c\in \adm(F)$, and let $c'\in \adm(F^1) \cong \adm(F^2)$ denote the corresponding colorings. Suppose that $c'$ colors the front half-bubbles of $F^1$ by $j$, the ``big'' facet by $i$, and the back half-bubbles by $k$. Then
\[
\brak{F^1, c'} = \frac{x_i + x_k}{x_j + x_k} \brak{F,c}\ \text{ and }\ \brak{F^2, c'} = \frac{x_i+x_j}{x_j + x_k} \brak{F,c},
\]
from which we obtain
\[
\brak{F,c} = \brak{F^1, c'} + \brak{F^2, c'}, 
\]
which completes the proof of relation \eqref{eq:type1 relation1}.

We now address the relation \eqref{eq:type1 relation2}. Let $G$ denote the foam on the left-hand side of the equation, and let $G'$ denote the foam on the right-hand side. Let $c\in \adm(G)$, and assume $c$ colors the ``big'' facet of $G$ by $i$, the front bubble by $j$, and the back bubble by $k$. Let $c'\in \adm(G)$ denote the coloring which is identical to $c$ except the front and back bubbles are colored by $k$ and $j$, respectively. Let $c''\in \adm(G')$ denote the coloring of $G'$ in which the depicted facet is colored $i$, and the remaining facets are colored according to $c$ (equivalently, $c'$). We claim that 
\[
\brak{G,c} + \brak{G,c'} = \brak{G',c''},
\]
which completes the proof. To verify the above equality, observe that 
\[
\brak{G,c} = \frac{x_i + x_k}{x_j + x_k} \brak{G',c''} \ \text{ and }\ \brak{G,c'} = \frac{x_i+ x_j}{x_j + x_k} \brak{G',c''}. 
\]

The proof of relation \eqref{eq:type1 relation3} is similar and left to the reader. 
\end{proof}

\vspace{0.1in} 

The previous lemma allows us to simplify the state space $\brak{\Gamma}_{\circ}$ assigned to a web $\Gamma\subset \P$ with a bigon region containing the puncture. 

\begin{proposition}
The two maps shown in Figure \ref{fig:bigon maps} are mutually inverse isomorphisms between state spaces of graphs in the theory $\brak{-}_{\circ}.$

\begin{figure}
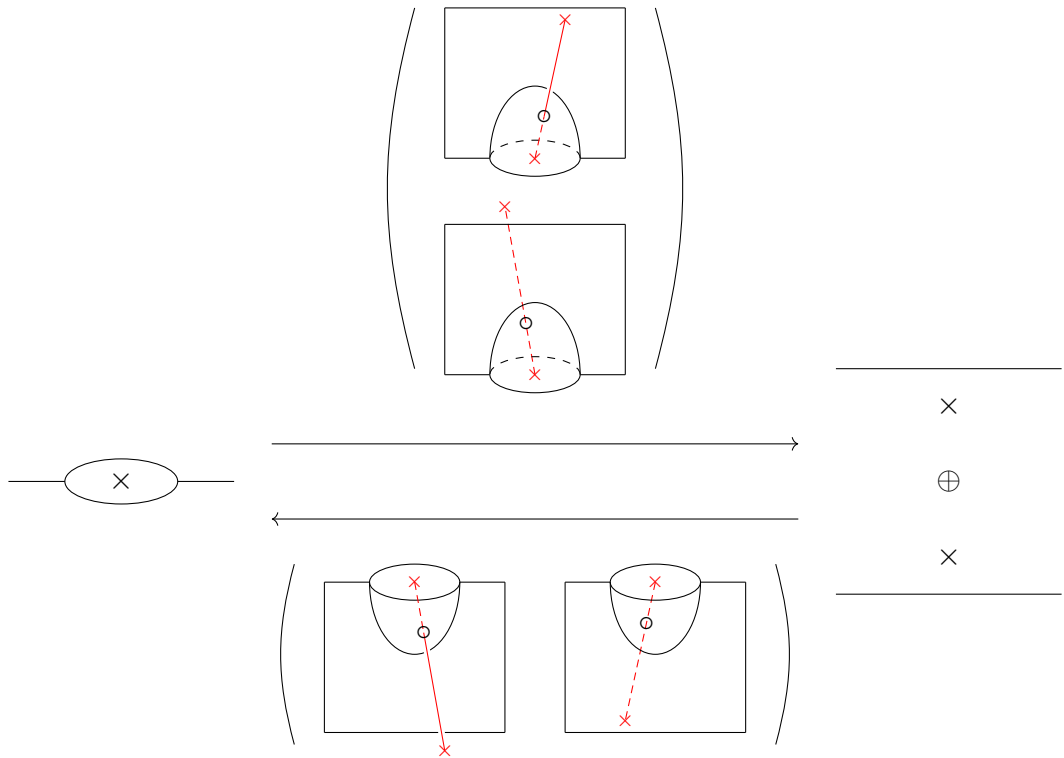

    \centering
    \includestandalone{bigon_maps}
    \caption{Isomorphisms which simplify a bigon region containing the puncture, for the theory $\brak{-}_{\circ}$. In the top map, the top foam has a Type 1 point on the front half-bubble, and the bottom foam has a Type 1 point on the back half-bubble. In the bottom map, the first foam has a Type 1 point on the front half-bubble, and the second foam has a Type 1 point on the back half-bubble.}
    \label{fig:bigon maps}
\end{figure}
\end{proposition}

\begin{proof}
This follows from the relations in Lemma \ref{lem:type1 relations}. 
\end{proof}

%
%

\section{Oriented \texorpdfstring{$SL(3)$}{SL(3)} anchored homology}\label{sec_oriented_hom} 

In this section we recall oriented $SL(3)$ foams, which were introduced in \cite{Khsl3} in the context of $sl(3)$ link homology. An equivariant analogue was defined in \cite{MV}, see also \cite{MN, Clark, MPT, Mac, Robert} for various aspects of $SL(3)$ foams and link homology. In Section \ref{sec:oriented sl3 foams} we define an evaluation of oriented $SL(3)$ foams via colorings in the style of Robert-Wagner \cite{RWfoamev} and show in Theorem \ref{thm:universal sl3 via coloring} that our evaluation agrees with that of \cite{MV}. In Section \ref{sec:oriented anchored sl3 foams} we deform the evaluation in the presence of the anchor line $L$. In Theorem \ref{thm:oriented sl3 eval is polynomial} we show that our evaluation is always a polynomial.  

To avoid introducing new notation, in this section we will reuse the notation for various rings from Section \ref{sec:unoriented sl3}. 

\begin{itemize}
\item $R_x' = \Z[x_1, x_2, x_3]$ is the ring of polynomials in three variables. 
\item $R_x=\Z[E_1, E_2, E_3]$ the subring of $R_x'$ that consists of symmetric polynomials in $x_1,x_2,x_3$, with generators $E_i$ being elementary symmetric polynomials: 
\begin{align*}
    E_1 &= x_1 + x_2 + x_3, \\
    E_2 &= x_1x_2 + x_1x_3 + x_2 x_3, \\
    E_3 &= x_1x_2x_3. 
\end{align*}
\item $R_x'' = R_x'[(x_1-x_2)^{-1}, (x_2-x_3)^{-1}, (x_1-x_3)^{-1}]$ is a localization of $R_x'$ given by inverting $x_i-x_j$, for $1\le i<j \le 3$. 
\item $\til{R}_x' = R_x'[\sqrt{x_1-x_3}, \sqrt{x_2-x_3}, \sqrt{x_1-x_3}]$ is the extension of $R_x'$ obtained by introducing square roots of $\sqrt{x_i-x_j}$, for $1\leq i < j \leq 3$. 
\item $\til{R}_x'' = \til{R}_x'[(x_1-x_2)^{-1}, (x_2-x_3)^{-1}, (x_1-x_3)^{-1}]$ is a suitable localization of the ring $\til{R}_x'$. 
\end{itemize}
All five of these rings are graded by setting $\deg(x_1) = \deg(x_2) = \deg(x_3) = 2$. Inclusions of the above rings are summarized in the following diagram.

\begin{equation}
\label{eq:rings char 0}
\begin{tikzcd}[row sep = .3em, column sep =.2em ]
& & \til{R}_x' & \subset & \til{R}_x'' \\
& & \cup & & \cup \\
R_x & \subset & R_x' & \subset & R_x''
\end{tikzcd}
\end{equation}


\subsection{Oriented \texorpdfstring{$SL(3)$}{SL(3)} foams and their evaluations}
\label{sec:oriented sl3 foams}

We begin by recalling the definition of oriented $SL(3)$ foams from \cite[Section 3.2]{Khsl3}.

\begin{definition}
\label{def:oriented sl3 foam}
A (closed) \emph{oriented $SL(3)$ pre-foam} $F$ consists of the following data.
\begin{itemize}
    \item An orientable surface $F'$ with connected components $F_1, \ldots, F_k$ and a partition of the boundary components of $F'$ into triples. The underlying CW structure of $F$ is obtained by identifying the three circles in each triple. The image of the three circles in each triple becomes a single circle in $F$, called a \emph{singular circle}. The image of the surfaces $F_i$ are called \emph{facets}. Three facets meet at each singular circle. 
    
    \item For each singular circle $Z$, we fix a cyclic ordering of the three facets meeting at $Z$. There are two possible choices of cyclic ordering for each $Z$. 
    
    \item Each facet may carry some number of dots, which are allowed to float freely along the facet but cannot cross singular circles. 
\end{itemize}

A \emph{oriented $SL(3)$ foam} is a pre-foam as above equipped with an embedding into $\RR^3$, along with an orientation on each facet such that any two of the three facets meeting at each singular circle are incompatibly oriented, as shown in Figure \ref{fig:facet orientation}. Each singular circle $Z$ acquires an induced orientation, see Figure \ref{fig:induced orientation}. This induced orientation on $Z$ specifies a cyclic ordering of the three facets meeting at $Z$ by following the left-hand rule, Figure \ref{fig:cyclic ordering}, and we require this to match the cyclic ordering specified by the pre-foam $F$.
\end{definition}

\begin{figure}
\centering
\subcaptionbox{Orientations of three facets meeting at a singular circle.
\label{fig:facet orientation}}[.3\linewidth]
{\includestandalone[scale=.7]{facet_orientation}
}
\subcaptionbox{ The induced orientation of a singular circle.
\label{fig:induced orientation}}[.3\linewidth]
{\includestandalone[scale=.7]{induced_orientation}
}
\subcaptionbox{ The induced cyclic ordering.
\label{fig:cyclic ordering}}[.3\linewidth]
{\includestandalone[scale=.7]{cyclic_ordering}
}
\caption{}\label{fig:orientation conventions}
\end{figure}

Note that unlike unoriented foams considered in Section \ref{sec:unoriented sl3}, the oriented $SL(3)$ pre-foams in the present section do not contain singular vertices. When there is no risk of confusion between the foams introduced in the Definition \ref{def:oriented sl3 foam} and those of Section \ref{sec:unoriented sl3}, in this section we will simply write (pre-)foam rather than oriented $SL(3)$ (pre-)foam. 

For a pre-foam $F$, let $\Theta(F)$ denote the set of its singular circles and $\theta(F) = \lr{\Theta(F)}$ the number of singular circles. Each $Z\in \Theta(F)$ has a neighborhood homeomorphic to the product of a circle $\SS^1$ and a tripod. Let $f(F)$ denote the set of facets of $F$. We use the definitions of pre-admissible and admissible colorings of pre-foams and foams from Section \ref{sec:unoriented sl3} in the present situation. For a pre-foam $F$, $\adm(F)$ denotes the set of admissible colorings of $F$. Note that if $F$ is a foam, every pre-admissible coloring is also admissible.  

Fix a pre-foam $F$ and an admissible coloring $c\in \adm(F)$. For $1\leq i \neq j \leq 3$, bicolored surfaces $F_{ij}(c)$ consist of all facets colored $i$ or $j$; each $F_{ij}(c)$ is a closed, orientable surface. For $1\leq i \leq 3$, let $F_i(c)$ be the surface consisting of all facets of $F$ which are colored $i$ by $c$; the surface $F_i(c)$ is orientable and has $\theta(F)$ boundary components. Denote by $\b{F}_i(c)$ the closed surface obtained by gluing disks along boundary components of $F_i(c)$. We have 
\begin{equation}
\label{eq:chi of surfaces}
\begin{aligned}
\chi(\b{F}_i(c)) &= \chi(F_i(c)) + \theta(F),\ 1\leq i \leq 3 \\
\chi(F_{ij}(c)) &= \chi(F_i(c)) + \chi(F_j(c)),\ 1\leq i < j \leq 3.
\end{aligned}
\end{equation}

The three facets meeting at each singular circle are colored by $i,j,k$, where as before we use $i,j,k$ to denote the three elements of $\{1, 2, 3\}$. We now define quantities $\theta^{\pm}(c)$, $\theta^{\pm}_{ij}(c)$ associated with the set of singular circles $\Theta(F)$ and the admissible coloring $c$.

\begin{definition}
Let $F$ be a pre-foam with admissible coloring $c$, and let $1\leq i < j \leq 3$. A singular circle $Z\in \Theta(F)$ is \emph{positive} with respect to $(i,j)$ if the cyclic ordering of the colors of the three facets meeting at $Z$ is $(i\ k\ j)$. If $F$ is a foam, then an equivalent formulation is as follows: when looking along the orientation of $Z$ with the facet colored $k$ drawn below, the $i$-colored facet is to the left of the $j$-colored facet. Otherwise, we say $Z$ is negative with respect to $(i,j)$. See Figure \ref{fig:ij circle} for a pictorial definition. Let $\theta^+_{ij}(c)$ (resp. $\theta^{-}_{ij}(c)$) denote the number of positive (resp. negative) circles with respect to $(i,j)$. We have
\[
\theta^+_{ij}(F,c) + \theta^{-}_{ij}(F,c) = \theta(F).
\]

We say that a singular circle $Z$ is \emph{positive} with respect to $c$ if the colors of the three facets meeting at $Z$ are $(1\ 2\ 3)$ in the cyclic ordering, and otherwise $Z$ is negative, see Figures \ref{fig:positive circle} and \ref{fig:negative circle}. Let $\theta^+(F,c)$ (resp. $\theta^{-}(F,c)$) denote the number of positive (resp. negative) circles in $F$ with respect to $c$. We have 
\begin{equation}
\label{eq:positive and negative circles}
\theta^+(F,c) + \theta^-(F,c) = \theta(F). 
\end{equation}
We will often omit $F$ from the notation and simply write $\theta$, $\theta^{\pm}_{ij}(c)$, and $\theta^{\pm}(c)$.
\end{definition}

\begin{figure}
\centering
\subcaptionbox{A positive $(i,j)$-circle, where $i<j$.
\label{fig:ij circle}}[.3\linewidth]
{\includestandalone[scale=.7]{ij_circle}
}
\subcaptionbox{ A positive singular circle.
\label{fig:positive circle}}[.3\linewidth]
{\includestandalone[scale=.7]{positive_circle}
}
\subcaptionbox{ A negative singular circle.
\label{fig:negative circle}}[.3\linewidth]
{\includestandalone[scale=.7]{negative_circle}
}
\caption{}\label{fig:signed circles}
\end{figure}

We now define the evaluations $\brak{F,c}$ and $\brak{F}$. For a pre-foam $F$, $c\in \adm(F)$, and $1\leq i \leq 3$, let $d_i(c)$ denote the number of dots on facets colored $i$. Define 
\begin{align}
    P(F,c) &= \prod_{i=1}^3 x_i^{d_i(c)} \label{eq:oriented sl3 numerator}\\
    Q(F,c) &= \prod_{1\leq i < j \leq 3} (x_i - x_j)^{\chi(F_{ij}(c))/2} \label{eq:oriented sl3 denominator}\\
    s(F,c) &= \sum_{i=1}^3 i \chi(\b{F}_i(c))/2 + \sum_{1\leq i < j \leq 3 } \theta^+_{ij}(c). \label{eq:sl3 sign}
\end{align}
Set 
\begin{align}
    \brak{F,c} & = (-1)^{s(F,c)} \frac{P(F,c)}{Q(F,c)},\\
    \brak{F} & = \sum_{c\in \adm(F)} \brak{F,c}. \label{eq:oriented pre foam eval}
\end{align}
A priori, the evaluations $\brak{F,c}$ and $\brak{F}$ lie in the ring $R_x''$ (see diagram \eqref{eq:rings char 0}).  

In what follows, we use the symbol $\equiv$ to mean equality modulo $2$. Note that
\begin{equation}
    \sum_{i=1}^3 i \chi(\b{F}_i(c))/2 \equiv \frac{\chi(\b{F}_1(c)) + \chi(\b{F}_3(c)) }{2},
\end{equation}
since $\chi(\b{F}_2(c))$ is even. Moreover, from \eqref{eq:chi of surfaces} we obtain 
\begin{equation}
\label{eq:alternate formula}
\sum_{i=1}^3 i \chi(\b{F}_i(c))/2 \equiv \theta + \sum_{i=1}^3 i \chi(F_i(c))/2.
\end{equation}

\begin{lemma}
\label{lem:alternate sign expression}
For a pre-foam $F$ and $c\in \adm(F)$, we have 
\[
\sum_{1\leq i < j \leq 3} \theta^+_{ij}(c) \equiv \theta^+(c).
\]
It follows that 
\begin{equation}
\label{eq:alternate sign expression}
    s(F,c) \equiv \sum_{i=1}^3 i \chi({F}_i(c))/2 + \theta^-(c).
\end{equation}
\end{lemma}
\begin{proof}
Let $Z\in \Theta(F)$. Observe that if $Z$ is positive with respect to $c$, then it contributes only to $\theta_{13}^+(c)$. Likewise, if $Z$ is negative then it contributes to both $\theta_{12}^+(c)$ and $\theta_{23}^+(c)$ but not to $\theta_{13}^+(c)$, which verifies the first equality. The second equality follows from equations \eqref{eq:alternate formula} and \eqref{eq:positive and negative circles}.
\end{proof}

\begin{example}
\label{ex:oriented sl3 sphere}
Let $F$ be a $2$-sphere $\SS^2$ with $d$ dots. For $1\leq i \leq 3$, let $c_i\in \adm(F)$ color $F$ by $i$. We have 
\begin{align*}
    \brak{F} &= \brak{F,c_1} + \brak{F,c_2} + \brak{F,c_3} \\
    &= -\frac{ x_1^d}{(x_1-x_2)(x_1 - x_3)} + \frac{x_2^d}{(x_1-x_2)(x_2-x_3)} - \frac{x_3^d}{(x_1-x_3)(x_2-x_3)} \\
    &= \frac{ -x_1^d(x_2-x_3) + x_2^d(x_1-x_3) - x_3^d(x_1-x_2)}{(x_1-x_2)(x_2-x_3)(x_1 - x_3)} \\
    &= -s_{(d-2,0,0)}(x_1,x_2,x_3) = -h_{d-2}(x_1,x_2,x_3) = - \sum_{i+j+k=d-2} x_1^i x_2^j x_3^k,
\end{align*}
where $s_{(d-2,0,0)}(x_1,x_2,x_3)$ is the Schur function of the partition $(d-2,0,0)$, and $h_{d-2}(x_1,x_2,x_3)$ is the complete symmetric function of degree $d-2$. In particular $\brak{F}=0$ if $d=0$ or $d=1$, and $\brak{F}=-1$ if $d=2$. 
\end{example}

\begin{example}
\label{ex:oriented sl3 theta foam}
Let $F$ be the theta foam shown in \eqref{eq:oriented theta foam}.

\begin{equation}
\label{eq:oriented theta foam}
    \begin{aligned}
    \includestandalone{oriented_theta_foam}
    \end{aligned}
\end{equation}

Given any $c\in \adm(F)$, each capped-off surface $\b{F}_i(c)$ and each bicolored surface $F_{ij}(c)$ is a $2$-sphere. In particular,
\[
s(F,c) \equiv \theta^+(c).
\]
For $\sigma \in S_3$, let $c(\sigma)\in \adm(F)$ denote the coloring which colors the top facet by $\sigma(1)$, the middle facet by $\sigma(2)$, and the bottom facet by $\sigma(3)$. We have 
\[
\brak{F} = \sum_{\sigma\in S_3} \brak{F, c(\sigma)} = \frac{\sum_{\sigma\in S_3} (-1)^{\theta^+(c(\sigma))} x_{\sigma(1)}^{d_1} x_{\sigma(2)}^{d_2} x_{\sigma(3)}^{d_3}}{(x_1-x_2)(x_1-x_3)(x_2-x_3)},
\]
and moreover 
\[
\theta^+(c(\sigma)) \equiv \lr{\sigma},
\]
where $\lr{\sigma}$ is the length of $\sigma$. 

Therefore if $d_1\geq d_2\geq d_3$, we have 
\[
\brak{F} = s_{(d_1-2, d_2-1, d_3)}(x_1,x_2,x_3),
\]
the Schur function with partition $(d_1-2, d_2-1, d_3)$. In particular, $\brak{F}=0$ if $d_1, d_2, d_3$ are not distinct.  If $d_1, d_2, d_3$ are distinct and $d_1+d_2+d_3\leq 3$, then up to cyclic permutation there are two choices, for which the evaluation is recorded in \eqref{eq:theta foam evaluation}. 
\begin{equation}
\label{eq:theta foam evaluation}
    \begin{aligned}
    \includestandalone{theta_foam1} \hskip3em \includestandalone{theta_foam2}
    \end{aligned}
\end{equation}

\end{example}

The symmetric group $S_3$ naturally acts on $\adm(F)$ and on the five rings in the diagram \eqref{eq:rings char 0}. The following lemma is analogous to \cite[Lemma 2.16]{RWfoamev}. 

\begin{lemma}
\label{lem:oriented sl3 S_3 action}
Let $F$ be a pre-foam, $c\in \adm(F)$, and $\sigma\in S_3$. Then 
\[
\sigma  ( \brak{F,c}) = \brak{F, \sigma  (c)}.
\]
\end{lemma}

\begin{proof}

We may assume that $\sigma$ is a transposition $(i \ i+1)$ for $i=1,2$. We have 
\[
\sigma ( P(F,c) )= P(F,\sigma(c)), \hskip2em \sigma (  Q(F,c)) = (-1)^{\chi(F_{i(i+1)}(c))/2} Q(F,\sigma   (c)). 
\]

Let $k \in \{1,2,3\}\setminus \{i, i+1\}$. Note that a singular circle $Z$ is positive with respect to $c$ if and only if $Z$ is negative with respect to $\sigma(c)$, so 
\[
\theta^+(c) + \theta^+ (\sigma(c))= \theta = \theta^- (c) + \theta^- (\sigma(c)).
\]
Moreover, we have 
\[\ F_i(c) = F_{i+1}(\sigma(c)), \  F_{i+1}(c) = F_i(\sigma(c)), \ F_k(c) = F_k(\sigma(c)).
\]
Therefore
\begin{align*}
    s(F,c) - s(F,\sigma(c)) & =  \frac{\chi(F_{i+1}(c)) - \chi(F_i(c))}{2} + \theta^-(c) - \theta^-(\sigma(c)) \\
    & \equiv  \frac{\chi(F_{i+1}(c)) - \chi(F_i(c))}{2} + \theta \\
    &\equiv \frac{\chi(F_{i+1}(c)) + \chi(F_i(c))}{2} \\
    & \equiv \frac{\chi(F_{i(i+1)}(c))}{2},
\end{align*}
which completes the proof.
\end{proof}

\begin{corollary}
The evaluation $\brak{F}$ is a symmetric rational function.
\end{corollary}

Later we will prove that $\brak{F}$ is in fact a polynomial, see Corollary \ref{cor:oriented sl3 eval is polynomial}.

\begin{lemma}
\label{lem:Kempe move sign}
Let $i\in \{1,2\}$, let $F$ be a pre-foam, and let $c\in \adm(F)$ be an admissible coloring. Suppose $c'\in \adm(F)$ is obtained from $c$ by an $(1,2)$-Kempe move along a surface $\Sigma \subset F_{12}(c)$. Then 
\[
s(F,c) \equiv s(F,c') + \frac{\chi(\Sigma)}{2}.
\]
\end{lemma}

\begin{proof}
Note that this is analogous to \cite[Lemma 2.19]{RWfoamev}. Letting $\theta(\Sigma)$ denote the number of seam circles on $\Sigma$, we have 
\[
\theta^-(c) + \theta^-(c') \equiv \theta(\Sigma) \equiv \chi(F_1(c) \cap \Sigma).
\]
Note also that 
\begin{align*}
\chi(F_1(c)) - \chi(F_1(c')) & = \chi(F_1(c) \cap \Sigma) - \chi(F_{2}(c) \cap \Sigma),\\
\chi(F_{2}(c)) - \chi(F_{2}(c')) &= \chi(F_{2}(c) \cap \Sigma) - \chi(F_1(c) \cap \Sigma).
\end{align*}
We compute: 
\begin{align*}
    s(F,c) - s(F,c') & \equiv \frac{\chi(F_1(c)) - \chi(F_1(c'))}{2}  + \frac{2(\chi(F_{2}(c)) - \chi (F_{2}(c')))}{2} + \theta(\Sigma) \\
    & \equiv \frac{\chi(F_{2}(c) \cap \Sigma )  - \chi(F_1(c) \cap \Sigma)}{2} + \chi(F_1(c) \cap \Sigma) \\
    &\equiv \frac{\chi(\Sigma)}{2}.
\end{align*}
\end{proof}


\subsection{Oriented anchored \texorpdfstring{$SL(3)$}{SL(3)} foams and their evaluations}
\label{sec:oriented anchored sl3 foams}

In this section we introduce (oriented) anchored $SL(3)$ foams and their evaluations. 

\begin{definition}
An \emph{oriented anchored $SL(3)$ foam} $F$ is an oriented foam $F'\subset \RR^3$ that may intersect the anchor line $L$ at finitely many points away from the singular circles of $F'$, so that each intersection point belongs to some facet of $F'$, and moreover these intersections are required to be transverse. Denote by $p(F) = F'\cap L$ the set of intersection points (anchor points) of $F$. The anchor points carry labels in $\{1,2,3\}$; that is, $F$ comes equipped with a fixed map 
\[
\l \: p(F) \to \{1,2,3\}.
\]
\end{definition}

Fix an anchored foam $F$ and an admissible coloring $c$ of the underlying foam $F'$. Each anchor point $p\in p(F)$ lying on a facet $f$ inherits a color $c(p) := c(f)$. As in Section \ref{sec:unoriented sl3}, we say that $c$ is an admissible coloring of the anchored foam $F$ if for each $p\in p(F)$, the color of $p$ equals the label of $p$, that is, $c(p) = \l(p)$. Denote by $\adm(F)$ the set of admissible colorings of $F$.  

For $i\in \{1,2,3\}$, let $i', i''$ denote the complementary elements, so that $\{i,i', i''\} = \{1,2,3\}$. Define the evaluations 
\begin{align}
    \brak{F,c} &= (-1)^{s(F,c)} \frac{P(F,c)}{Q(F,c)}\left( \prod_{p\in p(F)} (-1)^{c(p)-1} (x_{c(p)} - x_{\l(p)'})(x_{c(p)} - x_{\l(p)''}) \right)^{1/2} , \label{eq:oriented anchored sl3 eval color}\\
    \brak{F} &= \sum_{c\in \adm(F)} \brak{F,c}, \label{eq:oriented anchored sl3 eval}
\end{align}
where $P(F,c)$, $Q(F,c)$ and $s(F,c)$ are as defined in \eqref{eq:oriented sl3 numerator}, \eqref{eq:oriented sl3 denominator}, and \eqref{eq:sl3 sign}, respectively.

Let us explain the square root in equation \eqref{eq:oriented anchored sl3 eval color}. We have $c(p) = \l(p)$ for every anchor point $p\in p(F)$. If $p$ is labeled $i$, then it contributes 
\[
(-1)^{i-1} (x_i - x_j)(x_i - x_k) 
\]
to the product under the square root. More concretely, the product of the two terms under the square root, for a fixed anchor point $p$, is equal to 
\begin{align*}
(x_1-x_2)(x_1-x_3) \ \ &\mathrm{ if }\ \  c(p)=1, \\
(x_1-x_2)(x_2-x_3) \ \ &\mathrm{ if }\ \  c(p)=2, \\
(x_1-x_3)(x_2-x_3) \ \ &\mathrm{ if }\ \  c(p)=3 .
\end{align*} 
Let $\anch(i)$ be the number of anchor points $p$ with $c(p)=i$. Then for $1\le i < j \le 3$ the sum $\anch(i)+\anch(j)$ is even, which follows from Proposition~\ref{prop_pre_adm}. 

We define the square root  as the product 
\begin{equation}
\label{eq:anchor point contribution}
   \widetilde{Q}(F,c) \ := \  \prod_{1\le i<j\le 3} (x_i-x_j)^{(\anch(i)+\anch(j))/2} 
\end{equation}
and rewrite formula (\ref{eq:oriented anchored sl3 eval color}) as 
\begin{equation}\label{eq_no_roots}
\begin{split}
    \brak{F,c} \ &:= \ (-1)^{s(F,c)} \frac{P(F,c)\widetilde{Q}(F,c)}{Q(F,c)} \\
    &= 
    (-1)^{s(F,c)}P(F,c)  \prod_{1\le i<j\le 3} (x_i-x_j)^{(\anch(i)+\anch(j)-\chi(F_{ij}(c)))/2}.
    \end{split}
\end{equation}

Note that $\widetilde{Q}(F,c)$ depends only on the labels of anchor points and not on the coloring $c$, as long as $c$ respects labels of anchor points (otherwise, the evaluation is $0$). Consequently, it can also be denoted $\widetilde{Q}(F)$. Alternatively, it may be useful to allow more general colorings $c$, with $\widetilde{Q}(F,c)=0 $ for $c$ not compatible with the labels of anchor points. 

\vspace{0.1in} 

Recall diagram \eqref{eq:rings char 0} and the surrounding discussion for notations of various rings. The above formula implies the following proposition. 

\begin{proposition}
The evaluation $\brak{F,c}$ is an element of $R_x''$.
\end{proposition}

\begin{remark}
As discussed in Remark \ref{rmk:admissible colorings of anchored foams}, if $c$ is an admissible coloring of the underlying foam $F'$ but not of the anchored foam $F$, then the evaluation \eqref{eq:oriented anchored sl3 eval color} is still well-defined and equal to zero. Even if we don't restrict the notion of admissible colorings of an anchored foam to those which color anchor points according to their labels, additional terms in the evaluation will each be $0$, not contributing anything. 
\end{remark}

\begin{example}
\label{ex:oriented anchored sl3 sphere}
Let $F$ be a $2$-sphere $\SS^2$ carrying $d$ dots and intersecting $L$ twice. Then $\brak{F}=0$ unless both anchor points are labeled by $i \in \{1,2,3\}$. In this case, there is one admissible coloring $c$ which colors $F$ by $i$. We see that $s(F,c) \equiv i$, and the evaluation is
\[
\brak{F} = (-1)^i x_i^d.
\]
\end{example}

\begin{example}
\label{ex:oriented anchored sl3 theta foam}
Consider the theta foam $F$ whose facets each intersect $L$ exactly once, shown in \eqref{eq:oriented anchored theta foam}. There is one admissible coloring $c$, and we have 
\[
\brak{F} = \brak{F,c}  =
\begin{cases}
x_i^{d_1} x_j^{d_2} x_k^{d_3} & \text{ if } (i,j,k) = (1,3,2) \text{ or a cyclic permutation}, \\
-x_i^{d_1} x_j^{d_2} x_k^{d_3}  & \text{ if } (i,j,k) = (1,2,3) \text{ or a cyclic permutation}.
\end{cases}
\]

\begin{equation}
\label{eq:oriented anchored theta foam}
    \begin{aligned}
    \includestandalone{theta_foam_oriented}
    \end{aligned}
\end{equation}
\end{example}

The symmetric group $S_3$ acts on all five of the rings in diagram \eqref{eq:rings char 0}. Recall also that $S_3$ acts on the set of admissible colorings of an un-anchored foam (i.e., those considered in Section \ref{sec:oriented sl3 foams}). However, for an anchored foam $F$, $c\in \adm(F)$, and $\sigma\in S_3$, the coloring $\sigma(c)$ is in general not admissible for $F$. 

Consider instead the anchored foam $\sigma(F)$ defined as follows. The underlying foam of $\sigma(F)$ agrees with the underlying foam of $F$. If anchor points of $F$ are labeled by $\l \: p(F) \to \{1,2,3\}$, then the anchor points of $\sigma(F)$ are labeled by $\sigma(l) \: p \mapsto \sigma(\l(p))$. Note that $\sigma$ provides a bijection $\adm(F) \cong \adm(\sigma(F))$ via $c\mapsto \sigma(c)$. The following lemma says that the evaluations $\brak{F}$ and $\brak{\sigma(F)}$ differ by a sign, and moreover the sign depends only on $\sigma$ and on labels of anchor points of $F$.

\begin{lemma}
\label{lem:anchored oriented S_3 action}
For an anchored foam $F$, $c\in \adm(F)$, and $\sigma\in S_3$, we have 
\[
\sigma\left(\brak{F,c}\right) = (-1)^{\eps(F,\sigma)}  \brak{\sigma(F), \sigma(c)},
\]
where 
\begin{equation}\label{eq_sigma_eps}
\eps(F,\sigma) = \sum_{\substack{1\leq i < j \leq 3, \\ \sigma(i) > \sigma(j)} } \frac{\anch(i) + \anch(j)}{2}.
\end{equation}
It follows that
\[
\sigma\left(\brak{F}\right) = (-1)^{\eps(F,\sigma)}  \brak{\sigma(F)}.
\]
\end{lemma}

\begin{proof}
By Lemma \ref{lem:oriented sl3 S_3 action}, we have 
\[
\sigma \left( (-1)^{s(F,c)}
\frac{P(F,c)}{Q(F,c)}
\right) = (-1)^{s(\sigma(F), \sigma(c))} \frac{P(\sigma(F), \sigma (c))}{Q(\sigma(F), \sigma(c))}.
\]
It is clear that 
\[
\sigma( \til{Q}(F) ) = (-1)^{\eps(F,\sigma)} \til{Q}(\sigma(F)),
\]
and the first equality follows. For the second equality, we have
\begin{align*}
    \sigma\left(\brak{F}\right) &= \sum_{c\in \adm(F)} \sigma \left(  \brak{F,c} \right) \\
    &= (-1)^{\eps(F,\sigma)}  \sum_{c\in \adm(F)} \brak{\sigma(F), \sigma(c)} \\
    &= (-1)^{\eps(F,\sigma)} \brak{\sigma(F)}.
\end{align*}
\end{proof}

For $1\leq i \neq j \leq 3$, consider the ring
\[
R_{ij}'' := R_x'[(x_i-x_k)^{-1}, (x_j-x_k)^{-1}].
\]
Each $R_{ij}''$ is a subring of $R_x''$. A permutation $\sigma\in S_3$ sends $R_{ij}''$ isomorphically onto $R_{\sigma(i) \sigma(j)}''$. 

We are now ready for the main result of this section. 

\begin{theorem}
\label{thm:oriented sl3 eval is polynomial}
The evaluation $\brak{F}$ of an anchored foam is an element of $R_x'$, the polynomial ring in variables $x_1, x_2, x_3$. 
\end{theorem}

\begin{proof}
The proof is similar to that of \cite[Theorem 2.17]{KRfoamev} and \cite[Proposition 2.18]{RWfoamev}. 
By Lemma \ref{lem:anchored oriented S_3 action}, it suffices to show that $\brak{F} \in R_{12}''$ for any anchored foam $F$. This is because we may take a permutation $\sigma\in S_3$ sending $1$ to $i$ and $2$ to $j$, and consider the anchored foam $\sigma^{-1}(F)$. Then $\brak{\sigma^{-1}(F)} \in R_{12}''$ implies that 
\[
 \pm \brak{F} = \pm \brak{\sigma(\sigma^{-1}(F))} =\pm \sigma\left( \brak{\sigma^{-1}(F)}\right)   \in R_{ij}'',
\]
where the first equality comes from Lemma \ref{lem:anchored oriented S_3 action}. It follows that 
\[
\brak{F} \in R_{12}'' \cap R_{23}'' \cap R_{13}'' = R_x'.
\]

Let us show that $\brak{F} \in R_{12}''$. Partition $\adm(F)$ into equivalence classes as follows. For $c \in \adm(F)$, the class $C_c$ containing $c$ consists of colorings obtained from $c$ by performing a sequence of $(1,2)$ Kempe moves along surfaces in $F_{12}(c)$ which are disjoint from $L$. If $F_{12}(c)$ has $n$ connected components, $k\geq 0$ of which are disjoint from $L$, then $C_c$ consists of $2^k$ elements. We will show that 
\[
\sum_{c'\in C_c} \brak{F,c'} \in R_{12}'',
\]
which will conclude the proof. 

Write $\Sigma := F_{12}(c)$ as a disjoint union
\[
\Sigma = \Sigma' \cup \Sigma_1 \cup \cdots \cup \Sigma_k,
\]
where each $\Sigma_a$, $a=1, \ldots, k$ is connected and disjoint from $L$, and where each component of $\Sigma'$ intersects $L$. For $i=1,2$ and $a=1, \ldots, k$, let $t_i(a)$ denote the number of dots on $i$-colored facets (according to $c$) of $\Sigma_a$, and let $t_3$ denote the number of dots on $3$-colored facets (according to $c$) of $F$. We claim that 
\begin{equation}
\label{eq:sum formula}
\sum_{c'\in C_c} \brak{F,c'} = \frac{x_3^{t_3} \cdot \prod_{a=1}^k
\left(
x_1^{t_1(a)}x_2^{t_2(a)} + (-1)^{\chi(\Sigma_a)/2} x_2^{t_1(a)}x_1^{t_2(a)} \left(\frac{(x_1-x_3)}{(x_2 - x_3)}^{\l_{\Sigma_a}(c)/2} \right)
\right) 
\cdot \widetilde{Q}(F)
}
{(x_1 - x_2)^{\chi(\Sigma)/2} ( x_1 - x_3)^{\chi(F_{13}(c))/2} (x_2 - x_3)^{\chi(F_{23}(c))/2}
}
\end{equation}
where 
\begin{itemize}
    \item $\l_{\Sigma_a}(c) \in 2\Z$ is an even integer such that 
\[
\chi(F_{13}(c')) = \chi(F_{13}(c)) - \l_{\Sigma_a}(c) 
\ \text{ and }\ \chi(F_{23}(c')) = \chi(F_{23}(c)) + \l_{\Sigma_a}(c)
\]
for the coloring $c'\in C_c$ which is obtained from $c$ by a $(1,2)$ Kempe move along $\Sigma_a$. See \cite[Lemma 2.12 (3)]{KRfoamev} for details regarding this integer.

\item $\widetilde{Q}(F)$ is the contribution from the anchor points of $F$, equation \eqref{eq:anchor point contribution}. 
\end{itemize}
To verify the claimed equality, expand the product to  obtain $2^k$ terms, each of which corresponds to one of the $2^k$ colorings in $C_c$. That the sign is correct follows from Lemma \ref{lem:Kempe move sign}.

Finally, we argue that $(x_1 - x_2)^{\chi(\Sigma)/2}$ divides the numerator of \eqref{eq:sum formula}. Positive contributions to $\chi(\Sigma)$ come from $2$-sphere components of $\Sigma$. Each $\Sigma_a$ which is a $2$-sphere contributes one to the exponent $\chi(\Sigma)/2$. On the other hand, the corresponding factor in the product in the numerator of \eqref{eq:sum formula} is divisible by $x_1 - x_2$. The remaining positive contributions to $\chi(\Sigma)/2$ come from $2$-sphere components of $\Sigma'$. Such a component $\Sigma_0$ contains at least two anchor points, each labeled $1$ or $2$, so the contribution from $\Sigma_0$ can be cancelled with terms in $\widetilde{Q}(F)$. 
\end{proof}

\begin{corollary}
\label{cor:oriented sl3 eval is polynomial}
If $F$ is a pre-foam or a foam which is disjoint from $L$, then $\brak{F} \in R_x$, the ring of symmetric polynomials in $x_1, x_2, x_3$.
\end{corollary}

\begin{proof}
This follows from Lemma \ref{lem:anchored oriented S_3 action} and Theorem \ref{thm:oriented sl3 eval is polynomial}. 
\end{proof}

\subsection{Skein relations}

In this section we record several local relations involving oriented anchored $SL(3)$ foams.

\begin{lemma}
\label{lem:oriented sl3 relations}
The local relations \eqref{eq:three dots}, \eqref{eq:neck cutting oriented sl3}, \eqref{eq:bigon identity}, and \eqref{eq:square identity} hold for anchored foams. Seam lines are drawn in bold in relation \eqref{eq:square identity} to clarify the picture. 
\begin{equation}
\label{eq:three dots}
    \begin{aligned}
    \includestandalone{three_dots}
    \end{aligned}
\end{equation}

\begin{equation}
\label{eq:neck cutting oriented sl3}
    \begin{aligned}
    \includestandalone{neck_cutting_oriented_sl3}
    \end{aligned}
\end{equation}

\begin{equation}
\label{eq:bigon identity}
    \begin{aligned}
    \includestandalone{bigon_identity}
    \end{aligned}
\end{equation}

\begin{equation}
\label{eq:square identity}
    \begin{aligned}
    \includestandalone{square_identity}
    \end{aligned}
\end{equation}

\end{lemma}

\begin{proof}

Proofs of these four relations are similar to Propositions 2.33, 2.22, 2.23, and 2.24 in \cite{KRfoamev}, respectively, with the caveat that we must keep track of the sign \eqref{eq:sl3 sign}. Moreover, $S_3$ symmetry is used in \cite{KRfoamev} to simplify the calculations. Anchor points and their labels are the same for the foams depicted in each of these four relations, so Lemma \ref{lem:anchored oriented S_3 action} implies that we may use $S_3$ symmetry in a similar manner.  

We verify relation \eqref{eq:neck cutting oriented sl3} and leave the remaining three relations to the reader. Let $F$ denote the foam appearing on the left-hand side of the equality. The six foams on the right-hand side are identical except for placement of dots. We denote them by $G^1, \ldots, G^6$, so that the relation reads
\[
\brak{F} = -\left( \brak{G^1} + \brak{G^2} + \brak{G^3}\right)  +E_1 \left( \brak{G^4} + \brak{G^5} \right) - E_2 \brak{G^6} .
\]
Admissible colorings of $G^1, \ldots, G^6$ are in canonical bijection. For $c\in \adm(G^1)$, let 
\[
\brak{G,c}:= -\left( \brak{G^1,c} + \brak{G^2,c} + \brak{G^3,c}\right)  +E_1 \left( \brak{G^4,c} + \brak{G^5,c} \right) - E_2 \brak{G^6,c} .
\]
There are two types of colorings of $G^1$: those which color the two depicted disks the same, and those which color them differently. Those of the first type are in canonical bijection with colorings of $F$. 

Suppose $c\in \adm(G^1)$ colors both disks the same color, say $i$, and denote by $c\in \adm(G^2) \cong \cdots \cong \adm(G^6)$ and $c'\in \adm(F)$ the corresponding colorings. We will show that $\brak{F,c'} = \brak{G,c}$. We may assume $i=1$. Then
\[
\brak{G^1,c} = \brak{G^2,c} = \brak{G^3,c} = x_1^2 \brak{G^6,c}, \ \ \brak{G^4,c} = \brak{G^5,c} = x_1 \brak{G^6,c}, 
\]
which yields
\[
\brak{G,c} = -3x_1^2 \brak{G^6,c} + 2E_1 x_1 \brak{G^6,c} - E_2 \brak{G^6,c} =-(x_1-x_2)(x_1-x_3) \brak{G^6,c}. 
\]
To compare this with $\brak{F,c'}$, observe that
\[
\chi(F_1(c')) + 2= \chi(G^6_1(c)), \ \ \chi(F_2(c')) = \chi(G^6_2(c)), \ \ \chi(F_3(c')) = \chi(G^6_3(c)),
\]
which implies $s(F,c') \equiv s(G,c) + 1$. Moreover, we have
\[
\chi(F_{12}(c')) + 2 = \chi(G^6_{12}(c)), \ \  \chi(F_{13}(c')) + 2 = \chi(G^6_{13}(c)), \ \ \chi(F_{23}(c')) = \chi(G^6_{23}(c)).
\]
Therefore 
\[
\brak{G^6,c} = -\frac{\brak{F,c'}}{(x_1-x_2)(x_1-x_3)},
\]
which verifies $\brak{F,c'}  = \brak{G,c}$.

To complete the proof, suppose that $c$ colors the top depicted disk by $i$ and the bottom disk by $j$, with $i\neq j$. We have 
\[
\brak{G^1, c} = x_i^2 \brak{G^6, c}, \ \  \brak{G^2,c} = x_ix_j \brak{G^6,c}, \ \   \brak{G^3,c} = x_j^2 \brak{G^3,c},
\]
\[
\brak{G^4,c} = x_i \brak{G^6,c}, \ \   \brak{G^5,c} = x_j \brak{G^6,c} .
\]
Therefore $\brak{G,c}= 0$, which concludes the proof.
\end{proof}

\begin{lemma}
\label{lem:bubble removal}
Let $F$ be an anchored foam. Denote by $F_{n,m}$ the anchored foam obtained from $F$ by adding a bubble (disjoint from $L$) to some facet in $F$, with the two new facets carrying $n$ and $m$ dots respectively, such that the facet with $n$ dots directly precedes the facet with $m$ dots in the cyclic ordering. Let $F_n$ denote the foam obtained from $F$ by adding $n$ dots to the same facet. This is shown in \eqref{eq:bubble removal}. Then 
\begin{align*}
    \brak{F_{n,n}} &= 0, \\
    \brak{F_{1,0}} &= - \brak{F_{0,1}} = \brak{F}, \\
    \brak{F_{2,0}} &= - \brak{F_{0,2}}  = E_1 \brak{F} - \brak{F_1}. 
\end{align*}

\begin{equation}
    \label{eq:bubble removal}
    \begin{aligned}
    \includestandalone{bubble_removal}
    \end{aligned}
\end{equation}
\end{lemma}

\begin{remark}
The relations in Lemmas \ref{lem:oriented sl3 relations} and \ref{lem:bubble removal} also hold for pre-foams. 
\end{remark}

Similar to the $SL(2)$ and unoriented $SL(3)$ setting, for oriented $SL(3)$ foams we allow shifted dots $\circled{i} = \bullet - x_i$ ($1\leq i \leq 3)$ on a facet. 
\begin{equation*}
    \begin{aligned}
    \includestandalone{shifted_dot_oriented_sl3}
    \end{aligned}
\end{equation*}
They must be disjoint from $L$ and are allowed to float freely on their facets but cannot cross seam lines.

\begin{lemma}
\label{lem:oriented anchored sl3 neck cutting}
The local relations shown in \eqref{eq:oriented anchored neck cutting}, \eqref{eq:cup off line oriented sl3}, and \eqref{eq:line past seam oriented} hold. 

\begin{equation}
\label{eq:oriented anchored neck cutting}
    \begin{aligned}
    \includestandalone{neck_cutting_line_oriented}
    \end{aligned}
\end{equation}

\begin{equation}
\label{eq:cup off line oriented sl3}
    \begin{aligned}
    \includestandalone{moving_cup_off_line_oriented_sl3}
    \end{aligned}
\end{equation}

\begin{equation}
    \label{eq:line past seam oriented}
    \begin{aligned}
    \includestandalone{line_past_seam_oriented}
    \end{aligned}
\end{equation}
\end{lemma}

In the last equation we assume $j<k$. 

\begin{proof}
We verify equation \eqref{eq:oriented anchored neck cutting} and leave the remaining relations to the reader. The argument is similar to that of relation \eqref{eq:neck cutting line unoriented} in Lemma \ref{lem:unoriented sl3 relations}, so we will be brief. Let $F$ denote the foam on the left-hand side, and let $G^1, G^2, G^3$ denote the three foams on the right-hand side, with superscript corresponding to labels of the anchor points. For $1\leq i \leq 3$, let $\adm_i(F)$ consist of all admissible colorings of $F$ which color the depicted tube by $i$. There is a natural bijection $\adm_i(F) \cong \adm(G^i)$. 

Given $c\in \adm_i(F)$, let $c'\in \adm(G^i)$ denote the corresponding coloring. Arguing as in the proof of Lemma \ref{lem:unoriented sl3 relations}, we obtain 
\[
\brak{F,c} = \pm \brak{G^i, c'}.
\]
It remains to show that the above sign is equal to $(-1)^i$. We have
\[
\chi(F_j(c)) = \chi(G_j^i(c')), \ \ \chi(F_k(c)) = \chi(G_k^j(c')),\ \ \chi(F_i(c)) = \chi(G_i^i(c')) - 2, \ \ \theta^\pm(F,c) = \theta^\pm(G^i,c'),
\]
so $s(F,c) \equiv s(G^i,c') + i$ as needed. 

\end{proof}

\subsection{State spaces}
\label{sec:oriented sl3 state spaces}

In this section we define state spaces associated to oriented $SL(3)$ webs. Much of this is analogous to notions in Section \ref{sec:state spaces unoriented}.

\begin{figure}
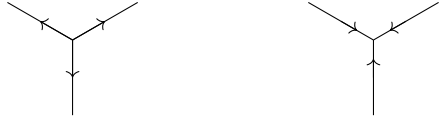

    \centering
    \includestandalone{oriented_webs}
    \caption{The orientations at each trivalent vertex of an oriented $SL(3)$ web must be either all outgoing or all incoming.}
    \label{fig:oriented webs}
\end{figure}

\begin{definition}
An \emph{oriented $SL(3)$ web} is a planar trivalent graph $\Gamma \subset \P$ in the punctured plane, which may have closed loops with no vertices. Moreover, edges and loops of $\Gamma$ carry orientations such that each vertex is either a source or a sink, as shown in Figure \ref{fig:oriented webs}. In this section we will simply write \emph{web} rather than oriented $SL(3)$ web. 
\end{definition}

The definition of an anchored foam with boundary in the oriented setting is analogous to that of Definition \ref{def:unoriented anchored foam with boundary}. The singular graph of a foam with boundary $V$ is a union of finitely many arcs (with boundary in $\RR^2 \times \{0,1\}$) and circles (disjoint from $\RR^2 \times \{0,1\}$). Intersection points of $V$ with $L_{[0,1]}$ (anchor points) must be disjoint from the singular graph and carry labels in $\{1,2,3\}$. Facets of $V$ are required to carry orientations satisfying the convention in Figure \ref{fig:facet orientation} near singular points. As usual, we will use the left-hand rule to specify these orientations and cyclic orderings by orienting each singular circle and arc, as shown in Figures \ref{fig:induced orientation} and \ref{fig:cyclic ordering}. 

\begin{figure}
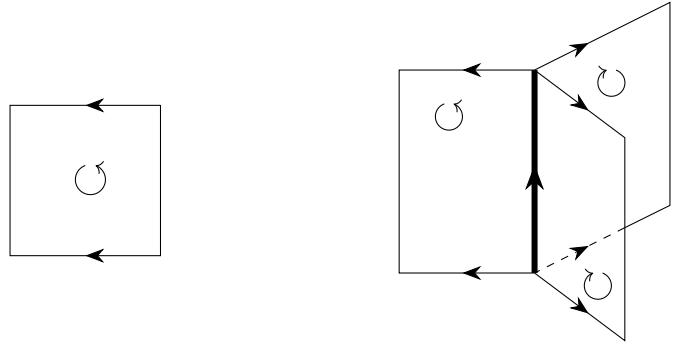

\centering
\includestandalone{oriented_foam_induced_orientation} \hskip7em
\includestandalone{oriented_foam_induced_orientation2}
\caption{Our convention for the induced orientation on the webs $\d_0 V$ (bottom) and $\d_1 V$ (top). 
}\label{fig:foam with boundary induced orientation}
\end{figure}

As in Section \ref{sec:state spaces unoriented}, let $\d_i V := V \cap (\RR^2 \times \{0\} )$ for $i=0,1$. The orientation of facets of $V$ induces an orientation on $\d_0 V$ and $\d_1 V$ via the convention in Figure \ref{fig:foam with boundary induced orientation}. We view $V$ as a cobordism from the oriented web $\d_0 V$ to the oriented web $\d_1 V$. Composition $WV$ of foams $V, W$ with $\d_1 V = \d_0 W$ is defined in the natural way. 

Denote by $p(V) = V \cap L_{[0,1]}$ the set of anchor points of $V$ and by $\lr{d(V)}$ the number of dots. The degree of $V$ is defined to be 
\begin{equation}
\label{eq:oriented foam degree}
\deg(V) =  2 \left( \lr{d(V)} + \lr{p(V)} - \chi(V) \right)+ \chi(\d V) .
\end{equation}
Degree is clearly additive under composition and is compatible with the grading on $R_x'$, in the sense that if $V$ is a closed foam, then $\deg(V) = \deg (\brak{V})$. 

As in Definition \ref{def:annular cobordism}, by an \emph{annular foam} we mean a foam (with boundary) which is disjoint from $L$. The composition of two annular foams is again annular.

There is an involution $\omega$ defined by reflecting a foam with boundary through $\RR^2\times \{1/2\}$. We have $\d_1 V = \d_0 (\omega(V))$ and $\d_0 V = \d_1 (\omega(V))$ for any foam with boundary $V$. Given a web $\Gamma \subset \P$, let $\Fr(\Gamma)$ denote the free $R_x'$-module generated by foams with boundary $V$ from the empty web to $\Gamma$ (that is, $\d_0 V = \varnothing$, $\d_1 V = \Gamma$). Define a bilinear form 
\[
(-,-) \: \Fr(\Gamma) \times \Fr(\Gamma) \to R_x'
\]
by $(V,W) = \omega(V)W$. This bilinear form is symmetric since $\brak{F} = \brak{\omega(F)}$ for any closed foam $F$. The state space $\brak{\Gamma}$ is the quotient of $\Fr(\Gamma)$ by the kernel 
\[
\ker((-,-)) = \{ x\in \Fr(\Gamma) \mid (x,y) = 0 \text{ for all } y\in \Fr(\Gamma) \}
\]
of the bilinear form, 
\[
\brak{\Gamma} := \Fr(\Gamma)/ \ker((-,-)). 
\]
The state space $\brak{\Gamma}$ inherits the grading from $\Fr(\Gamma)$ since $(-,-)$ is degree-preserving. A foam with boundary $V$ from $\Gamma_0 $ to $\Gamma_1$ naturally induces a map
\[
\brak{V} \: \brak{\Gamma_0} \to \brak{\Gamma_1}
\]
of degree $\deg(V)$, defined by sending the equivalence class of a basis element $U\in \Fr(\Gamma_0)$ to the equivalence class of $VU \in \Fr(\Gamma_1)$.  This assignment is functorial with respect to composition of foams, $\brak{WV} = \brak{W} \brak{V}$ for composable $V, W$.

\begin{figure}
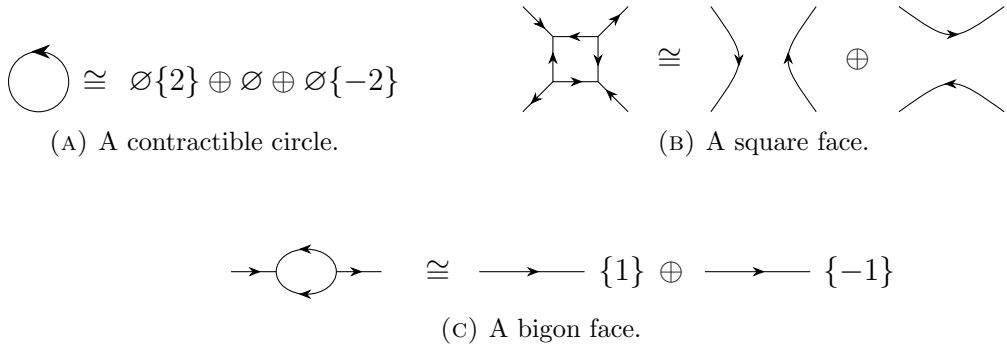

\centering 
\subcaptionbox{A contractible circle.}[.3\linewidth]
{\includestandalone{oriented_contractible_circle}
}  \hskip2em
\subcaptionbox{A square face.}[.5\linewidth]
{\includestandalone{oriented_square}
}\\ \vskip6ex
\subcaptionbox{A bigon face.}[.5\linewidth]
{\includestandalone{oriented_bigon}
}
\caption{Local relations for state spaces of oriented $SL(3)$ webs, where the depicted regions do not contain the puncture.}\label{fig:oriented local relations}
\end{figure}

\begin{lemma}
\label{lem:oriented sl3 local isos}
The three local isomorphisms shown in Figure \ref{fig:oriented local relations} hold. 
\end{lemma}

\begin{proof}
The arguments for relations (a), (b), and (c) of the figure are analogous to Propositions 7, 9, and 8, respectively, of \cite{Khsl3}. The relevant relations are Lemmas \ref{lem:oriented sl3 relations} and \ref{lem:bubble removal}. 
\end{proof}

\begin{proposition}
\label{prop:oriented sl3 noncontractible circle}
Let $\Gamma \subset \P$ be a web with a non-contractible circle $C$ which bounds a disk in $\RR^2\setminus \Gamma$, and let $\Gamma'=\Gamma \setminus C$ be the web obtained by removing $C$. Then there is an isomorphism 
\[
\brak{\Gamma} \cong \brak{\Gamma'} \oplus \brak{\Gamma'} \oplus \brak{\Gamma'}
\]
given by the maps shown in \eqref{eq:noncontractible iso oriented} (orientation of the circle is omitted). 

\begin{equation}
    \label{eq:noncontractible iso oriented}
    \begin{aligned}
    \includestandalone{noncontractible_iso_oriented}
    \end{aligned}
\end{equation}
\end{proposition}

\begin{proof}
This follows from Example \ref{ex:oriented anchored sl3 sphere} and the neck-cutting relation \eqref{eq:oriented anchored neck cutting}. 
\end{proof}

\begin{theorem}
\label{thm:oriented state spaces}
For any web $\Gamma \subset \P$, the state space $\brak{\Gamma}$ is a free graded $R_x'$-module of rank equal to the number of Tait colorings of $\Gamma$. Moreover, if $\Gamma$ is contractible, then the graded rank of $\brak{\Gamma}$ equals the Kuberberg polynomial \cite{Kup} of $\Gamma$, normalized as in \cite[Section 2]{Khsl3}. 
\end{theorem}

\begin{proof}
Lemma \ref{lem:graphs} (1) guarantees that we can reduce $\brak{\Gamma}$ to a direct sum of empty webs by recursively applying the local isomorphisms in Lemma \ref{lem:oriented sl3 local isos} and Proposition \ref{prop:oriented sl3 noncontractible circle}. It is then clear that the rank equals the number of Tait colorings. 

If $\Gamma$ is contractible then $\brak{\Gamma}$ can be simplified using only the isomorphisms in Lemma \ref{lem:oriented sl3 local isos}. Upon taking graded ranks, these isomorphisms recover the recursive relations for computing the Kuperberg polynomial. 
\end{proof}

Theorem \ref{thm:oriented state spaces} does not address the graded rank of state spaces of non-contractible webs. These may be computed recursively. As a special case, if $\Gamma$ consists of $n$ contractible and $m$ non-contractible circles, then $\brak{\Gamma}$ is free of graded rank $3^m (q^2 + 1 + q^{-2})^n$. 

\vspace{0.1in} 

Given a web $\Gamma \subset \P$, we can forget the puncture and the anchor line $L$ and apply the universal construction to the evaluation \eqref{eq:oriented pre foam eval}. Precisely, let $\Fr(\Gamma)_{\rm forget}$ denote the free $R_x$-module generated by all foams with boundary $\Gamma$ (forgetting the anchor line).  By Corollary \ref{cor:oriented sl3 eval is polynomial}, we can define the bilinear form $(-,-) \: \Fr(\Gamma)_{\forget} \times \Fr(\Gamma)_{\forget} \to R_x$ and the corresponding state space $\brak{\Gamma}_{\forget}$ in the usual way. Thus we obtain state spaces for webs in $\RR^2$, functorial with respect to foams in $\RR^2\times [0,1]$. These state spaces and maps induced by foams are graded via equation \eqref{eq:oriented foam degree}, where $\lr{p(V)} =  0$.

\begin{proposition}
For a contractible web $\Gamma \subset \P$, there is a degree-preserving isomorphism
\[
\brak{\Gamma} \cong \brak{\Gamma}_{\forget},
\]
natural with respect to foams with contractible boundary and which are disjoint from $L$. 
\end{proposition}

\begin{proof}
This follows from Theorem \ref{thm:oriented state spaces}. 
\end{proof}

On the other hand, Mackaay-Vaz \cite{MV} define an evaluation $\brak{-}_{\MV}$ for oriented $SL(3)$ pre-foams and use it to define an equivariant (also called \emph{universal}) version of the $sl(3)$ link homology introduced in \cite{Khsl3}. They work over the ground ring $\Z[a,b,c]$ and associate a state space $\brak{\Gamma}_{\MV}$ to each web $\Gamma \subset \RR^2$ via the universal construction applied to their pre-foam evaluation $\brak{-}_{\MV}$. To compare with our situation, identify $\Z[a,b,c]$ with the ring $R_x = \Z[E_1, E_2, E_3]$ of symmetric functions in $x_1, x_2, x_3$ via a ring isomorphism $\varphi$ defined by $\varphi(a)= E_1$, $\varphi(b)= -E_2$, $\varphi(c)= E_3$.

\begin{theorem}
\label{thm:universal sl3 via coloring}
For any closed pre-foam $F$, we have 
\[
\brak{F} = \varphi\left({\brak{F}_{\MV}}\right).
\]
It follows that there are isomorphisms $\brak{\Gamma}_{\forget} \cong \brak{\Gamma}_{\MV}\o_{\Z[a,b,c]} R_x$ for any web $\Gamma\subset \RR^2$, natural with respect to maps induced by foams with boundary. 
\end{theorem}

\begin{proof}
The evaluation $\brak{-}_{\MV}$ is defined by applying the local relations (3D), (CN), (S), and $(\Theta)$ in \cite[Section 2.1]{MV} to reduce any foam to an element of $\Z[a,b,c]$. Under the change of variables $a\mapsto E_1$, $b\mapsto -E_2$, $c \mapsto E_3$, these four relations hold for our evaluation $\brak{-}$ by relation \eqref{eq:three dots}, relation \eqref{eq:neck cutting oriented sl3}, Example \ref{ex:oriented sl3 sphere}, and Example \ref{ex:oriented sl3 theta foam}. The statement follows. 
\end{proof}

As in the $SL(2)$ and unoriented $SL(3)$ setting considered earlier in the paper, we can define an additional grading on oriented $SL(3)$ foams and state spaces. Define the abelian group 
\begin{equation}
    \Lambda = \Z  w_1 \oplus  \Z w_2\oplus \Z w_3  /(w_1+w_2+w_3), 
\end{equation}
on three generators and one relation. $\Lambda$ is a free abelian group of rank  two.

\begin{figure}
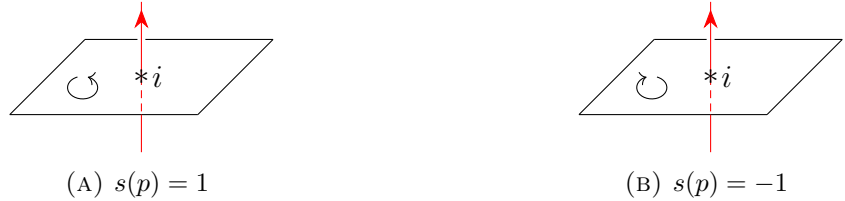

\centering    
\subcaptionbox{$s(p) = 1$ \label{fig:positive intersection}}[.45\linewidth]
{\includestandalone{positive_intersection} 
}
\subcaptionbox{$s(p)=-1$ \label{fig:negative intersection}}[.45\linewidth]
{\includestandalone{negative_intersection}
}
\caption{The oriented intersection number between a facet and $L$.}
\label{fig:oriented intersection}
\end{figure}

Orient the anchor line $L$ from bottom to top. For an anchored foam $V$ with boundary and $p\in p(V)$ an anchor point lying on some facet $f$, let $s(p) \in \{\pm 1\}$ denote the oriented intersection number between $f$ and $L$ ($s(p)$ does not depend on the label of $p$),  see Figure \ref{fig:oriented intersection} for the convention. Define the \emph{annular degree of $V$} to be 
\begin{equation}
    \label{eq:oriented adeg}
    \adeg(V) = \sum_{p\in p(V} s(p)w_{\l(p)} \in \Lambda.  
\end{equation}

\begin{proposition}\label{prop_deg_0}
If $F$ is a closed anchored foam with an admissible coloring $c$, then $\adeg(F) = 0$.
\end{proposition}
\begin{proof}
The proof is similar to that of Proposition \ref{prop_pre_adm}. The intersection of $F$ with a generic half-plane that bounds $L$ is an oriented web $\Gamma$ with boundary points on $L$. An admissible coloring $c$ of $F$ induces a Tait coloring of $\Gamma$. The boundary points (one-valent vertices) of $\Gamma$ are colored according to their label. The sum in \eqref{eq:oriented adeg} may be rewritten as the sum of terms $\pm (w_1 + w_2 + w_3) =0$ over all trivalent vertices of $\Gamma$, where the sign is $+1$ if all edges are incoming and $-1$ if all edges are outgoing. Each $i$-colored inner edge $e$ of $\Gamma$ bounds two trivalent vertices and contributes $\pm (w_i - w_i)=0$ since $e$ is oriented towards one of its boundary vertices and away from the other. The remaining edges, with one or both endpoints on $L$, contribute precisely $\adeg(F)$. 
\end{proof}

Let $\Gamma \subset \P$ be an (annular oriented) $SL(3)$ web. An anchored foam $F\subset \RR^3_-$ with $\partial F=\Gamma$ has a well-defined degree $\adeg(F)\in \Lambda$ via (\ref{eq:oriented adeg}). Furthermore, we equip the coefficient ring  $R_x'$ with an $\Lambda$-grading, with all elements of degree $0$. This makes free $R_x'$-module $\Fr(\Gamma)$ into an $\Lambda$-graded module, and 
Proposition~\ref{prop_deg_0} implies that the kernel of the bilinear form on $\Fr(\Gamma)$ is $\Lambda$-graded as well. Consequently, the grading descends to an $\Lambda$-grading on the state
space $\brak{\Gamma}$. A foam $V$ with boundary  induces a map $\brak{V} \: \brak{-\partial_0\Gamma} \to \brak{\partial_1\Gamma}$ which changes $\adeg$ by $\adeg(V)$. If $V$ has no anchor points, it induces an annular degree $0$ map between the state spaces of its boundaries. 
The state space of a contractible web is concentrated in annular degree zero. 

$\Lambda$-grading on $\brak{\Gamma}$ is the analogue of grading on finite-dimensional $SL(3)$ representations by the weight lattice. In fact, in the non-equivariant version of our construction, where all $x_i$'s are set to $0$ upon closed foam evaluation (and state spaces are defined accordingly, over a ground field rather than the ring $R_x'$), the state space $\brak{\Gamma}$ is naturally an $\mathfrak{sl}_3$-representation. We also refer the reader to  Queffelec-Rose~\cite{QR} for the construction of sutured annular $\mathfrak{sl}_n$-homology, with state spaces of annular webs carrying an $\mathfrak{sl}_n$-action. 
In the equivariant case, it is not clear how to define an $\mathfrak{sl}_3$-action or what's the substitute for it. 

\vspace{0.1in}

Denote by $\AFoamor$ the category whose objects consist of oriented $SL(3)$ webs in $\P$ and whose morphisms are $R_x'$-linear combinations of anchored  cobordisms between webs. Morphism spaces in this category are triply graded via $(\qdeg, \adeg)$. The state space construction assembles into a functor
\[
\brak{-} \: \AFoamor \to R_x'\gggmod
\]
landing in the category of triply-graded $R_x'$-modules.

This functor respects the trigradings on the hom spaces in the two categories.
Restricting to the subcategory of annular cobordisms and their linear combinations, the induced maps have annular degree $0$.

\subsection{Annular \texorpdfstring{$SL(3)$}{SL(3)}-link homology}

Let $L\subset \A \times [0,1]$ be a link in the thickened annulus. Projecting onto $\A \times \{0\} = \A$ and identifying the interior of $\A$ with the punctured plane $\P$, we obtain a link diagram $D\subset \P$. Following \cite{Khsl3, MV}, form the cube of resolutions of $D$, with each resolution web drawn in the punctured plane. Applying the functor $\brak{-} \: \AFoamor \to R_x' \gggmod $ yields a chain complex $C(D)$ of $\Z \oplus \Lambda$-graded $R_x'$-modules, with the $\Z$-grading given by $\deg$, equation \eqref{eq:oriented foam degree},  and the $\Lambda$-grading given by $\adeg$, equation \eqref{eq:oriented adeg}. Degree shifts in the cube of resolutions are applied only to the $\Z$-degree $\deg$. Diagrams in $\P$ representing isotopic annular links are related by Reidemeister moves away from the puncture. Proofs of Reidemeister invariance in \cite{MV} are local, and all local relations (away from $L$) on foams in \cite{MV} also hold for our evaluation $\brak{-}$ by \eqref{eq:neck cutting oriented sl3}, Example \ref{ex:oriented sl3 sphere}, and Example \ref{ex:oriented sl3 theta foam}. It follows that the chain homotopy class of $C(D)$ is an invariant of the annular link $L$. We define \emph{equivariant annular $SL(3)$ homology} as cohomology groups $H(C(D))$. 

Moreover, foams between webs appearing in the cube of resolutions are disjoint from $L$. Thus the differential preserves annular degree throughout the complex. Consequently, equivariant annular $sl(3)$ link homology carries a homological grading as well as an internal $\Z \oplus \Lambda$-grading $(\deg, \adeg)$. Cohomology groups  $H(C(D))$ are trigraded $R'_x$-modules.

\bibliographystyle{alpha}
\bibliography{references}
\end{document}